\documentclass[12pt]{article}
\usepackage{amssymb,amsmath,latexsym}
\allowdisplaybreaks

\begin{document}

\begin{doublespace}

\def\1{{\bf 1}}
\def\ind{{\bf 1}}
\def\nn{\nonumber}
\newcommand{\I}{\mathbf{1}}

\def\sA {{\cal A}} \def\sB {{\cal B}} \def\sC {{\cal C}}
\def\sD {{\cal D}} \def\sE {{\cal E}} \def\sF {{\cal F}}
\def\sG {{\cal G}} \def\sH {{\cal H}} \def\sI {{\cal I}}
\def\sJ {{\cal J}} \def\sK {{\cal K}} \def\sL {{\cal L}}
\def\sM {{\cal M}} \def\sN {{\cal N}} \def\sO {{\cal O}}
\def\sP {{\cal P}} \def\sQ {{\cal Q}} \def\sR {{\cal R}}
\def\sS {{\cal S}} \def\sT {{\cal T}} \def\sU {{\cal U}}
\def\sV {{\cal V}} \def\sW {{\cal W}} \def\X {{\mathfrak X}}
\def\sY {{\cal Y}} \def\sZ {{\cal Z}}

\def\bA {{\mathbb A}} \def\bB {{\mathbb B}} \def\bC {{\mathbb C}}
\def\bD {{\mathbb D}} \def\bE {{\mathbb E}} \def\bF {{\mathbb F}}
\def\bG {{\mathbb G}} \def\bH {{\mathbb H}} \def\bI {{\mathbb I}}
\def\bJ {{\mathbb J}} \def\bK {{\mathbb K}} \def\bL {{\mathbb L}}
\def\bM {{\mathbb M}} \def\bN {{\mathbb N}} \def\bO {{\mathbb O}}
\def\bP {{\mathbb P}} \def\bQ {{\mathbb Q}} \def\bR {{\mathbb R}}
\def\bS {{\mathbb S}} \def\bT {{\mathbb T}} \def\bU {{\mathbb U}}
\def\bV {{\mathbb V}} \def\bW {{\mathbb W}} \def\bX {{\mathbb X}}
\def\bY {{\mathbb Y}} \def\bZ {{\mathbb Z}}
\def\R {{\mathbb R}} \def\RR {{\mathbb R}} \def\H {{\mathbb H}}
\def\n{{\bf n}} \def\Z {{\mathbb Z}}

\newcommand{\expr}[1]{\left( #1 \right)}
\newcommand{\cl}[1]{\overline{#1}}
\newtheorem{thm}{Theorem}[section]
\newtheorem{lemma}[thm]{Lemma}
\newtheorem{defn}[thm]{Definition}
\newtheorem{prop}[thm]{Proposition}
\newtheorem{corollary}[thm]{Corollary}
\newtheorem{remark}[thm]{Remark}
\newtheorem{example}[thm]{Example}
\numberwithin{equation}{section}
\def\ee{\varepsilon}
\def\qed{{\hfill $\Box$ \bigskip}}
\def\NN{{\mathcal N}}
\def\AA{{\mathcal A}}
\def\MM{{\mathcal M}}
\def\BB{{\mathcal B}}
\def\CC{{\mathcal C}}
\def\LL{{\mathcal L}}
\def\DD{{\mathcal D}}
\def\FF{{\mathcal F}}
\def\EE{{\mathcal E}}
\def\QQ{{\mathcal Q}}
\def\SS{{\mathcal S}}
\def\RR{{\mathbb R}}
\def\R{{\mathbb R}}
\def\L{{\bf L}}
\def\K{{\bf K}}
\def\S{{\bf S}}
\def\A{{\bf A}}
\def\E{{\mathbb E}}
\def\F{{\bf F}}
\def\P{{\mathbb P}}
\def\N{{\mathbb N}}
\def\eps{\varepsilon}
\def\wh{\widehat}
\def\wt{\widetilde}
\def\pf{\noindent{\bf Proof.} }
\def\pff{\noindent{\bf Proof} }
\def\cp{\mathrm{Cap}}

\title{\Large \bf Accessibility, Martin boundary and minimal thinness for Feller processes in metric measure spaces}

\author{{\bf Panki Kim}\thanks{This work was  supported by the National Research Foundation of
Korea(NRF) grant funded by the Korea government(MSIP) (No. NRF-2015R1A4A1041675)
},
\quad {\bf Renming Song\thanks{Research supported in part by a grant from
the Simons Foundation (208236)}} \quad and
\quad {\bf Zoran Vondra\v{c}ek}
\thanks{Research supported in part by the Croatian Science Foundation under the project 3526}
}

\date{}

\maketitle

\begin{abstract}
In this paper we study the Martin boundary at infinity for a large class of purely discontinuous Feller processes in metric measure spaces. We show that if $\infty$ is accessible from an open set $D$, then there is only one Martin boundary point of $D$ associated with it, and this point is minimal. We also prove the analogous result for finite boundary points. As a consequence, we show that minimal thinness of a set is
a local property.
\end{abstract}



\section{Introduction and setup}\label{s:intro}

The Martin kernel and Martin boundary of an open set with respect to a transient strong Markov process
were introduced in \cite{KW} with the goal of representing non-negative harmonic functions (with respect to the
underlying process) as an integral of the Martin kernel against a finite measure on the (minimal) Martin boundary.
The identification of the Martin boundary for purely discontinuous Markov processes began in late nineties when it
was shown in \cite{B99, CS98}
that for the isotropic $\alpha$-stable process the Martin boundary of a bounded Lipschitz domain coincides with its Euclidean boundary.  Soon after, the result was extended in \cite{SW99} to the so-called $\kappa$-fat open sets. These results were subsequently extended in two directions: 
to more general processes and
to general open sets.

In the first direction, the Martin boundary of bounded $\kappa$-fat open sets was  studied in \cite{KSV1} for a class of subordinate Brownian motions and then in \cite{KSV14b} for some symmetric L\'evy processes. In both papers the Martin boundary was identified with the Euclidean boundary. In fact, the latter paper gives a local result: if an open set $D\subset \R^d$ is
$\kappa$-fat at $z_0\in\partial D$, then there is exactly one (minimal) Martin boundary point associated to $z_0$.
A related result is the identification of the Martin boundary at infinity of an unbounded open set with a single point provided the set is $\kappa$-fat at infinity, see \cite{KSV14}. In all of these papers  an appropriate boundary Harnack principle for non-negative harmonic functions played a major role.

In the second direction, the boundary Harnack principle and the Martin kernel for arbitrary open sets  were studied in
\cite{BKK} for isotropic $\alpha$-stable processes. The authors of \cite{BKK} introduced
the concepts of accessible and inaccessible boundary points
and proved a result that leads to the identification of the finite Martin boundary of an arbitrary bounded open with its
Euclidean boundary. It was also proved in \cite{BKK} that a finite Martin boundary point is minimal if and only if the
corresponding Euclidean boundary point is accessible.
By use of the Kelvin transform, they were able to identify the infinite part of the Martin boundary as well.

The main goal of this paper is to generalize results of \cite{BKK, KSV14}
to more general processes. Inspired by the paper \cite{BKuK} we will work with a 
class of purely discontinuous Feller processes
in duality in a measure metric space $\X$.
The jumps of these processes are assumed to be quite regular as precisely described in Assumptions \textbf{C}, \textbf{C1} and \textbf{C2} below.
Most of L\'evy processes fall into our framework, see Section \ref{s:exam} for details.
Our main results can be roughly stated as follows: let $D$ be an open set
in $\X$. If $z_0\in \partial D$ (the boundary of $D$ in the original topology of $\X$) is accessible,
then there is exactly one Martin boundary point associated with $z_0$. In case $D\subset \X$ is bounded and all its boundary points are accessible,
the Martin boundary and the minimal Martin boundary of $D$ are identified with $\partial D$.
In case of unbounded open set such that infinity is accessible, we identify the Martin boundary at infinity with a single point.

Another goal of this paper is to show that minimal thinness is a local property.
We will use our results on the Martin boundary to show that
under certain geometric assumptions, if $E\subset D\subset \X$ are open sets with a common boundary point $z_0$
which is accessible from both $E$ and $D$,
then $F\subset E$ is minimally thin at $z_0$ in $E$ if and only if $F$ is minimally thin at $z_0$ in $D$.

We now provide a precise description of the process and the assumptions it satisfies, introduce all necessary notation, state the results and explain the methods of proofs.

Let $(\X, d, m)$ be a metric measure space with a countable base
 such that all bounded closed sets are compact and the measure $m$ has full support. For $x\in \X$ and $r>0$, let $B(x,r)$ denote the ball centered at $x$ with radius $r$.
 Let $R_0\in (0,\infty]$ be the localization radius such that $\X\setminus B(x,2r)\neq \emptyset$ for all $x\in \X$ and all $r<R_0$.

Let $X=(X_t, \sF_t, \P_x)$ be a Hunt process on $\X$. We will assume the following

\smallskip
\noindent
\textbf{Assumption A:}
$X$ is a Hunt process admitting a strong dual process $\widehat{X}$ with respect to the measure $m$ and $\widehat{X}$ is also a Hunt process.
The transition semigroups $(P_t)$ and $(\widehat{P}_t)$ of $X$ and $\widehat{X}$ are both Feller and strongly Feller.
Every semi-polar set of $X$ is polar.

For the definition of Hunt processes see p. 45 of \cite{BG}, and for the definition of
a strong dual, see Definition VI.(1.2) on p.225 of \cite{BG}. For the definition of Feller
processes see pp. 49--50 
of \cite{CW}, and  for the definition of strong Feller processes see p. 129 of \cite{CW}. For the definitions of polar and semi-polar sets, see
Definition II.(3.1) on p.79 of \cite{BG}.

\smallskip
In the sequel, all objects related to the dual process $\widehat{X}$ will be denoted by a hat.
Recall that a set is polar (semi-polar, respectively) for $X$ if and only if it is polar (semi-polar, respectively)
for $\wh X$ (see VI. (1.19) in \cite{BG}).
Under assumption {\bf A} the process $X$ admits a (possibly infinite) Green function $G(x,y)$ serving as a density of the occupation measure: $G(x,A):=\E_x \int_0^{\infty}\ind_{(X_t\in A)} dt =\int_A G(x,y)m(dy)$. Moreover, $G(x,y)=\widehat{G}(y,x)$ for all $x,y\in \X$, 
cf. VI.1 in \cite{BG} for details. 
Further, let $D$ be an open subset of $\X$ and $\tau_D=\inf\{t>0:\, X_t\notin D\}$ the exit time from $D$. The killed process $X^D$ is defined by $X_t^D=X_t$ if $t<\tau_D$ and $X_t^D=\partial$
where $\partial$ is an extra point added to $\X$.
 The killed process $\widehat{X}^D$ is defined analogously. 
By Hunt's switching identity (Theorem 1.16 in  \cite{BG}), 
 it holds that $\E_x[G(X_{\tau_D},y)]=
\widehat{\E}_y[\widehat{G}(\widehat{X}_{\widehat{\tau}_D},x)]$
 for all $x,y \in \X$ which implies that $X^D$ and $\widehat{X}^D$ 
are in duality, see  p.43 in \cite{G71}.
Again by VI.1 in \cite{BG}, $X^D$ admits a unique (possibly infinite) 
Green function (potential kernel) $G_D(x,y)$ such that for every non-negative Borel function $f$,
$$
G_D f(x):=\E_x \int_0^{\tau_D}f(X_t)dt=\int_D G_D(x,y)f(y) \, m(dy)\,  ,
$$
and $G_D(x,y)=\widehat{G}_D(y,x)$, $x,y\in D$, with $\widehat{G}_D(y,x)$ the Green function of $\widehat{X}^D$.
It is assumed throughout the paper that $G_D(x,y)=0$ for $(x,y)\in (D\times D)^c$.
We also note that the killed process $X^D$ is strongly Feller, see e.g.~the
first part of the proof of Theorem on pp.~68--69 in \cite{Chu}.
From now on, we will always assume that $D$ is Greenian, that is, the Green function
$G_D(x,y)$ is finite for all $x,y\in D$, $x\neq y$. Under this assumption,
the killed process $X^D$ is
 transient in the sense that there exists a non-negative Borel function $f$ on $D$ such that $0<G_D f<\infty$ (and the same is true for $\widehat{X}$).

Recall that $z\in \partial D$ is said to be regular with respect to $X$ if $\P_z(\tau_D =0)=1$ and irregular otherwise.
We will denote the set of regular (respectively irregular) points of $\partial D$ with respect to $X$
by $D^{\mathrm{reg}}$ (respectively $D^{\mathrm{irr}}$).
$\widehat{D}^{\mathrm{reg}}$ (respectively $\widehat{D}^{\mathrm{irr}}$) stands
for the sets of regular (respectively irregular) points of $\partial D$ with respect to $\widehat{X}$ respectively.
It is well known that 
$D^{\mathrm{irr}}$ and $\widehat{D}^{\mathrm{irr}}$ 
are semipolar, hence polar under {\bf A}.

The process $X$, being a Hunt process, admits a L\'evy system $(J,H)$ where $J(x,dy)$
is a kernel on the state space $\X$ (called the L\'evy kernel of $\X$),
and $H=(H_t)_{t\ge 0}$ is a positive continuous additive functional of $X$. We assume that $H_t=t$
so that for every function $f:\X\times \X \to [0,\infty)$ vanishing on the diagonal and every stopping time $T$,
$$
\E_x \sum_{0<s\le T} f(X_{s-}, X_s)=\E_x \int_0^T f(X_s,y)J(X_s,dy) ds\, .
$$
By using $\tau_D$ in the displayed formula above and taking $f(x,y)=\1_D(x)\1_A(y)$ with $A\subset \overline{D}^c$, we get that
\begin{equation}\label{e:exit-distribution}
\P_x(X_{\tau_D}\in A, \tau_D<\zeta)=\E_x \int_0^{\tau_D} J(X_s, A) ds= \int_D G_D(x,y)J(y,A)m(dy)\, ,
\end{equation}
where $\zeta$ is the life time of $X$.
Similar formulae hold for the dual process $\wh{X}$ and $\wh{J}(x,dy)m(dx)=J(y,dx)m(dy)$.

\medskip
\noindent
\textbf{Assumption C:} The L\'evy kernels of $X$ and $\wh{X}$ are of the form $J(x,dy)=j(x,y)m(dy)$, $\wh{J}(x,dy)=\widehat{j}(x,y)m(dy)$, where $j(x,y)=\widehat{j}(y,x)>0$ for all $x,y\in \X$, $x\neq y$.

The next two related assumptions control the decay of the density $j$.

\noindent
\textbf{Assumption C1}$(z_0, R)$:
Let $z_0\in \X$ and $R\le R_0$. For all $0<r_1<r_2<R$, there exists a constant $c=c(z_0, r_2/r_1)>0$ such that for all $x\in B(z_0,r_1)$ and all $y\in \X\setminus B(z_0, r_2)$,
$$
 c^{-1} j(z_0,y)\le j(x,y)\le c j(z_0,y), \qquad  c^{-1} \widehat{j}(z_0,y)\le \widehat{j}(x,y)\le c \widehat{j}(z_0,y).
$$

In the next assumption we require that the localization radius $R_0=\infty$.

\noindent
\textbf{Assumption C2}$(z_0, R)$:
Let $z_0\in \X$ and $R>0$. For all $R\le r_1<r_2$, there exists a constant $c=c(z_0, r_2/r_1)>0$ such that for all $x\in B(z_0,r_1)$ and all $y\in \X\setminus B(z_0, r_2)$,
$$
 c^{-1} j(z_0,y)\le j(x,y)\le c j(z_0,y), \qquad  c^{-1} \widehat{j}(z_0,y)\le \widehat{j}(x,y)\le c \widehat{j}(z_0,y).
$$

\medskip
We \emph{define}  the Poisson kernel of an open set $D\in \X$ by
\begin{equation}\label{e:pk-def}
P_D(x,z)=\int_D G_D(x,y)j(y,z) m(dy), \qquad x\in D, z\in D^c .
\end{equation}
By \eqref{e:exit-distribution}, we see that $P_D(x,\cdot)$ is the density of the exit distribution of $X$ from $D$
restricted to $\overline{D}^c$:
$$
\P_x(X_{\tau_D}\in A, 
\tau_D <\zeta
)=\int_A P_D(x,z) m(dz), \qquad A\subset \overline{D}^c .
$$

Recall that $f:\X\to [0,\infty)$ is regular harmonic in $D$ with respect to $X$ if
$$
f(x)=\E_x[f(X_{\tau_D}),
\tau_D <\zeta
]\, , \quad \textrm{for all }x\in D\,  ,
$$
and it is harmonic in $D$ with respect to $X$ if for every relatively compact open $U\subset \overline{U}\subset D$
$$
f(x)=\E_x[f(X_{\tau_U}), \tau_U <\zeta]\, , \quad \textrm{for all }x\in U\,  .
$$
Throughout the paper we will adopt the convention that $X_{\zeta}=\partial$ and $f(\partial)=0$ for every function $f$. Thus we will drop $\tau_D<\zeta$ in expressions similar to the right-hand side in the penultimate display.
A function $f:\X\to [0,\infty)$ harmonic in $D$ with respect to $X^D$ if for every relatively compact open $U\subset \overline{U}\subset D$
$$
f(x)=\E_x[f(X^D_{\tau_U})]\, , \quad \textrm{for all }x\in U\,  .
$$
It follows from the Hunt switching formula that for every $y\in D$ and any open neighborhood $U$ of $y$, $G_D(\cdot, y)$ is regular harmonic in $D\setminus \overline{U}$. In particular, $G_D(\cdot, y)$ is harmonic in  $D\setminus \{y\}$.

The next pair of assumptions is about an approximate factorization of harmonic functions. This approximate factorization is
a crucial tool in proving the oscillation reduction. The first one is
an approximate factorization of harmonic functions at a finite boundary point.

\noindent
\textbf{Assumption F1}$(z_0, R)$:
Let $z_0\in \X$ and $R\le R_0$. For any $\frac{1}{2} < a < 1$, there exists $C(a)=C(z_0, R, a)\ge 1$ such that for
every $r\in (0, R)$, every open set $D\subset B(z_0,r)$, every non-negative function $f$ on $\X$
which is  regular harmonic in $D$ with respect to $X$ and vanishes in
$B(z_0,r)\cap(\overline{D}^c\cup D^{\mathrm{reg}})$,
and all $x\in D\cap B(z_0,r/8)$,
\begin{eqnarray}\label{e:af-1}
\lefteqn{C(a)^{-1}\E_x[\tau_{D}] \int_{\overline{B}(z_0,ar/2)^c}  j(z_0,y) f(y)m(dy)  } \nonumber \\
&&\le f(x) \le C(a)\E_x[\tau_{D}]\int_{\overline{B}(z_0,ar/2)^c}  j(z_0,y) f(y)m(dy).
\end{eqnarray}

In the second assumption we require that the localization radius $R_0=\infty$.

\noindent
\textbf{Assumption F2}$(z_0, R)$:
Let $z_0\in \X$ and $R>0$. For any $1 < a  < 2$, there exists  $C(a)=C(z_0, R, a)\ge 1$ such that for every $r\ge R$,
every unbounded open set $D\subset \overline{B}(z_0,r)^c$,
every non-negative function $f$ on $\X$ which is  regular harmonic in $D$
with respect to $X$ and vanishes on
$\overline{B}(z_0, r)^c\cap( \overline{D}^c\cup D^{\mathrm{reg}})$,
and all $x\in D\cap \overline{B}(z_0,8r)^c$,
\begin{eqnarray}\label{e:af-2}
\lefteqn{C(a)^{-1}\, P_{D}   (x,z_0) \int_{B(z_0, 2ar)}  f(z)m(dz)} \nonumber \\
&&\le f(x) \le
C(a)\, P_{D}   (x,z_0) \int_{B(z_0, 2ar)}  f(z)m(dz).
\end{eqnarray}

The approximate factorization of harmonic functions stated in {\bf F1} and {\bf F2} can be proved under somewhat
stronger assumptions than the Assumptions {\bf A}, {\bf B}, {\bf C} and {\bf D} in \cite{BKuK}.
This is done in the companion paper \cite{KSVp1}.

\medskip

Let $D\subset \X$ be an open set. A point $z\in \partial D$ is called accessible from $D$ with respect to $X$ if
\begin{equation}\label{e:accessible-finite}
P_D(x,z)=\int_D G_D(x,w)j(w,z) m(dw) = \infty \quad \text{ for all } x \in D\, ,
\end{equation}
and inaccessible otherwise.

In case $D$ is unbounded we say that $\infty$ is accessible from $D$ with respect to $X$ if
\begin{equation}\label{e:accessible-finite2}
\E_x \tau_D =\int_D G_D(x,w) m(dw)=  \infty \quad \text{ for all } x \in D
\end{equation}
 and inaccessible otherwise.
The concepts of accessible and inaccessible points were introduced in \cite{BKuK}.

For $D\subset \X$, let $\partial_M D$ denote the Martin boundary of $D$ with respect to $X^D$ in the sense of Kunita-Watanabe \cite{KW}, see Section 3 for more details. 
A point $w\in \partial_M D$ is said to be minimal if the Martin kernel $M_D(\cdot, w)$ is a minimal
harmonic function with respect to $X^D$. We will use $\partial_m D$ to denote the
minimal Martin boundary of $D$ with respect to $X^D$.
A point $w\in \partial_M D$ is said to be a \emph{finite Martin boundary point}
if there exists a bounded (with respect to the metric $d$) sequence $(y_n)_{n\ge 1}\subset D$
converging to $w$ in the Martin topology.
A point $w\in \partial_M D$ is said to be an \emph{infinite Martin boundary point}
if there exists an unbounded (with respect to the metric $d$) sequence $(y_n)_{n\ge 1}\subset D$
converging to $w$ in the Martin topology. We note that these two definitions do not rule out the possibility that a point $w\in \partial_M D$ is at the same time finite and infinite Martin boundary point. We will show in Corollary \ref{c:finite-not-infinite}(a) that under appropriate and natural assumptions this cannot happen. 
A point $w\in \partial_MD$ is said to be associated
with $z_0\in \partial D$ if there is a sequence $(y_n)_{n\ge 1}\subset D$ converging to $w$
in the Martin topology and to $z_0$ in the topology of $\X$. The set of Martin
boundary points associated with $z_0$ is denoted by $\partial_M^{z_0} D$.
A point $w\in \partial_MD$ is said to be associated with $\infty$ if $w$ is an infinite Martin boundary point.
The set of Martin boundary points associated with $\infty$ is denoted by $\partial_M^{\infty} D$.
$\partial^f_M D$ and $\partial^f_m D$ will be used to denote the finite part of the Martin
boundary and minimal boundary respectively.
Note that $\partial_M^{\infty} D$ is the set of infinite Martin boundary points.

Now we can state the first main result of the paper. We will always assume that Assumptions \textbf{A} and \textbf{C} hold true.

\begin{thm}\label{t:main-mb0}
Let $D\subset \X$ be an open set.
(a) Suppose that $z_0\in \partial D$.
Assume that there exists $R\le R_0$ such that
\textbf{C1}$(z_0, R)$  holds, and that
$\widehat{X}$ satisfies \textbf{F1}$(z_0, R)$.
If $z_0$
is accessible from $D$ with respect to $X$,
then there is only one Martin boundary point associated with $z_0$.

\noindent
(b) Suppose that, in addition to the assumptions in (a),
 for all $r \in (0, R]$,
\begin{align}\label{e:nGassup1}
\sup_{x\in D\cap B(z_0,r/2)} \sup_{y \in \X \setminus B(z_0, r)} \max(G_D(x, y), \widehat{G}_D(x, y))=:c(r) <\infty\, .
\end{align}
Then the Martin kernel $M_D(\cdot, z_0)$ is harmonic with respect to $X^D$.

\noindent
(c) Suppose, in addition, that $X$ satisfies \textbf{F1}$(z_0, R)$, that
\begin{align}\label{e:zreg}
\lim_{D\ni x\to z}G_D(x,y)=0 \quad \text{for every } z\in D^{\mathrm{reg}} \text{ and every }y\in D,
\end{align}
  and that, if $D$ is unbounded then for  $r \in (0, R]$,
\begin{align}\label{e:nGassup2}
\lim_{x \to \infty} G_D(x, y)=0 \quad \text{for all } y \in D \cap B(z_0, r).
\end{align}
Then the corresponding Martin boundary point is minimal.
\end{thm}

\begin{corollary}\label{c:main-mb0} Suppose that every point $z_0\in \partial D$ is accessible from $D$ with respect to $X$, and that the assumptions of Theorem \ref{t:main-mb0}(c) are satisfied for all $z_0\in \partial D$ (with $c(r)$ in \eqref{e:nGassup1} independent of $z_0$).

\noindent
(a)  The finite part of the Martin boundary $\partial_M D$ and the minimal Martin boundary $\partial_m D$ can be identified
with $\partial D$.

\noindent
(b) If $D$ is bounded, then $\partial D$ and $\partial_M D$ are homeomorphic.

\noindent
(c) Let $D$ be bounded. For any non-negative function $u$ which is harmonic with respect to $X^D$, there exists a unique finite measure $\mu$ on $\partial D$ such that
$$
u(x)=\int_{\partial D}M_D(x, z)\mu(dz), \qquad x\in D.
$$
\end{corollary}

\begin{thm}\label{t:main-mb1}
(a) Suppose that $R_0=\infty$, $D$ is an unbounded open subset of $\X$, and $\infty$ is accessible from $D$ with respect to $X$.
If there is a point $z_0\in\X$ and $R>0$ such that
\textbf{C2}$(z_0, R)$ holds,
and $\widehat{X}$ satisfies \textbf{F2}$(z_0, R)$,
then there is only one Martin boundary point associated with $\infty$.

\noindent
(b) Suppose that, in addition to the assumptions in (a),  for all $r\ge R$
\begin{align}\label{e:nGassup1-infty}
\sup_{x\in D\cap B(z_0,r/2)} \sup_{y \in \X \setminus B(z_0, r)} \max(G_D(x, y), \widehat{G}_D(x, y))=:c(r) <\infty\, .
\end{align}
Then the Martin kernel associated with $\infty$ is harmonic with respect to $X^D$.

\noindent
(c) Suppose, in addition,   that $X$ satisfies \textbf{F2}$(z_0, R)$, that
\eqref{e:zreg} holds, and that
\begin{align}\label{e:nGassup2i}
\lim_{x \to \infty} G_D(x, y)=0 \quad \text{for all } y \in D.
\end{align}
Then the Martin boundary point associated with $\infty$ is minimal.
\end{thm}

\begin{corollary}\label{c:finite-not-infinite} Let $R_0=\infty$ and $D\subset \X$ be unbounded. Suppose that every point $z_0\in \partial D$ is accessible from $D$ with respect to $X$, that $\infty$ is accessible from $D$ with respect to $X$, that the assumptions of Theorem \ref{t:main-mb0}(c) are satisfied for all $z\in \partial D$ (with $c(r)$ in \eqref{e:nGassup1} independent of $z$) and that the assumptions of Theorem \ref{t:main-mb1}(c) are satisfied. Then 

\noindent
(a) $\partial_M^f D\cap \partial_M^{\infty}D=\emptyset$.

\noindent 
(b) The Martin boundary $\partial_M D$ is homeomorphic with the one-point compactification of $\partial D$.

\noindent
(c) For any non-negative function $u$ which is harmonic with respect to $X^D$, there exists a unique finite measure $\mu$ on $\partial D$ and $\mu_{\infty}\ge 0$ such that
$$
u(x)=\int_{\partial D}M_D(x, z)\mu(dz)+M_D(x,\infty)\mu_{\infty}\, , \qquad x\in D,
$$
where $M_D(\cdot, \infty)$ denotes the Martin kernel associated with $\infty$.
\end{corollary}

The preliminary version of the results of this paper (and the forthcoming paper \cite{KSVp3}) was presented at
the 11th Workshop on Markov Processes and Related Topics held in Shanghai Jiaotong University from June 27 to June 30 2015,
and at the International Conference on Stochastic Analysis and Related Topics held in Wuhan University from August 3 to
August 8 2015.
In the recent preprint \cite{JK}, Juszczyszyn and Kwa\'snicki independently considered similar problems as those in
Corollary \ref{c:main-mb0}  for  bounded $D$. Our main motivation for the current paper was to investigate the Martin boundary at infinity. The investigation starts with the result stating that there is only one Martin boundary point associated with $\infty$ which should be understood as a local result about the Martin boundary in the sense that no other information about the remaining part of the boundary is required. This motivated our approach in studying the finite part of Martin boundary through the local approach -- if $z_0\in \partial D$ is accessible, then  there is only one Martin boundary point associated to $z_0$. Again, no other information about the remaining part of boundary is used.
We will first present proofs for infinity.
For readers' convenience, even though  the structure of proofs is similar, we also provide the proofs for finite boundary points.

The case of inaccessible boundary points  will be discussed in the forthcoming paper \cite{KSVp3}, the main reason being that the treatment of inaccessible points requires  additional assumptions on $j(x,y)$ --
see {\bf E1} and  {\bf E2} in \cite{KSVp3}, and  Theorem 3.1(ii) in \cite{JK}.

Organization of the paper:
In the next section we study the oscillation reduction at an accessible boundary point, first for the infinite point in Proposition \ref{p:oscillation-reduction-I}, and then for a finite boundary point in Proposition \ref{p:oscillation-reduction}. One of the main tools for this, borrowed from \cite{BKK}, is a decomposition of a regular harmonic function into two parts depending on where the process exits the open set. An estimate of one of the parts by the other is derived as a consequence of {\bf F2}, respectively {\bf F1}, cf.~Lemma \ref{l:assumption-2-I} and Lemma \ref{l:assumption-2}. The oscillation reduction result immediately leads to the existence of limits of ratios of non-negative harmonic functions which implies that the Martin kernel is the limit of ratios of Green functions.  This is the key to associating a point on the topological boundary of $D$ with a point on the Martin boundary.
The third section is devoted to the study of the Martin kernel at infinity under the assumption that infinity
 is accessible from $D$ and then of the Martin kernel at a finite
accessible point of an open set $D$.
We first prove that the Martin kernel is harmonic, and then that it is  minimal,
thus showing that a minimal Martin boundary point is associated with an accessible boundary point. In Section 4 we first briefly discuss
examples
satisfying our assumptions and then look at the case of a class of symmetric 
L\'evy processes in detail. 
In the last section we look at minimal thinness at a minimal Martin boundary point of $D$. It is intuitively clear that minimal thinness of a set $F\subset D$ should be a local property depending only on the size of $F$ near the boundary point. This suggests that if $F\subset E$, $E$ open in $D$, and $E$ and $D$ have a common boundary point,  then $F$ should be minimally thin at that boundary point in $E$ if and only if it is minimally thin in $D$. Clearly, the problem is that Martin spaces for $E$ and $D$ are different and one needs some sort of identification of
the underlying boundary points.
This is provided by Theorems \ref{t:main-mb0} and \ref{t:main-mb1}. The second main ingredient in showing local character of minimal thinness is given in Proposition \ref{p:h-minimal} where the Martin kernel with respect to $E$ is given in terms of the Martin kernel with respect to $D$.

Notation: We will use the following conventions in this paper: $c, c_0, c_1, c_2, \dots$ stand for constants
whose values are unimportant and which may change from one
appearance to another. All constants are positive finite numbers.
The labeling of the constants $c_0, c_1, c_2, \dots$ starts anew in
the statement of each result. We will use ``$:=$"
to denote a definition, which is  read as ``is defined to be".
We denote $a \wedge b := \min \{ a, b\}$, $a \vee b := \max \{ a, b\}$.
Notation $f\asymp g$  means that the quotient $f(t)/g(t)$ stays bounded
between two positive numbers on their common domain of definition.  For $x\in \X$ and $r>0$ we denote by $B(x,r)$ be the open ball centered at $x$ with radius $r$ and by $\overline{B}(x,r)$ the closure of $B(x,r)$.


\section{Oscillation reduction  under accessibility assumption}\label{s:osc-red}
It follows easily from the strong Markov property that for all Greenian open sets $U$ and $D$ with $U \subset D$,
$G_D(x,y)= G_U(x,y) + \E_x\left[   G_D(X_{\tau_U}, y); \tau_U<\infty\right]$ for every $(x,y) \in \X \times \X$.

\subsection{Infinity}

In this subsection, we deal with the oscillation reduction at infinity. Throughout this subsection we will  assume that
$R_0=\infty$ and that there exists a point $z_0\in\X$ such that \textbf{C2}$(z_0, R)$ holds, and that $\widehat{X}$ satisfies \textbf{F2}$(z_0, R)$ for some $R>0$.
We will fix $z_0$ and $R$ and use the notation $B_r=B(z_0, r)$.

An immediate consequence of \textbf{F2}$(z_0, R)$ for $\wh{X}$ is the boundary Harnack principle at infinity
in \cite{KSVp1}:
There exists $c>1$ such that for any $r\ge R$, any open set $D\subset \overline{B}_r^c$ and
any non-negative functions $u$ and $v$ on $\X$ that are regular harmonic in $D$ with respect to $\wh{X}$
and vanish on
$\overline{B}_r^c \cap( \overline{D}^c\cup \wh{D}^{\mathrm{reg}})$,
it holds that
\begin{equation}\label{e:bhp-I}
c^{-1}\frac{u(y)}{v(y)} \le \frac{u(x)}{v(x)} \le c \frac{u(y)}{v(y)}\qquad \text{for all } x,y\in D\cap \overline{B}_{8r}^c .
\end{equation}
Note that we can take $c=(C(3/2))^2$. By enlarging $C(3/2)$ in {\bf F2}$(z_0, R)$, without loss of generality we assume the above $c$ is equal to  $C(3/2)$.

For an open set $D$ and $p>q>0$, let
$
D^p=D\cap \overline{B}_p^c$ and $ D^{p,q}=D^q\setminus D^p.
$
For $p>q>1$, $r \ge R$ and non-negative function $f$ on $\X$  define
 \begin{eqnarray*}
 f^{pr,qr}(x)&=&
 \E_{x} \left[f(\wh{X}_{\wh{\tau}_{D^{pr}}}): \wh{X}_{\wh{\tau}_{D^{pr}}} \in D^{pr,qr}\right],
\\
 \wt{f}^{pr,qr}(x)&=&
 \E_{x} \left[f(\wh{X}_{\wh{\tau}_{D^{pr}}}): \wh{X}_{\wh{\tau}_{D^{pr}}} \in (D\setminus D^{qr})\cup \overline{B}_r\right].
 \end{eqnarray*}

\begin{lemma}\label{l:assumption-2-I}
Suppose  that $r\ge R$, $D\subset \overline{B}_r^c$ is an open set,
$f$ is a non-negative function on $\X$ which is  regular harmonic in $D$ with respect to $\wh{X}$,
and vanishes on
$\overline{B}_r^c\cap ( \overline{D}^c\cup \widehat{D}^{\mathrm{reg}})$.
There exists $C_1=C_1(R)>0$ (independent of $D$, $f$ and $r\ge R$) such that
for any  $p/16>q>2$ and any $\epsilon>0$, if
 \begin{equation}\label{e:assumption-2-I}
\int_{ \overline{B}_{qr}} f(y) m(dy) \le \epsilon \int_{ D^{pr/8,qr}} f(y) m(dy),
 \end{equation}
then for every $x\in D^{pr}$,
$ \wt{f}^{pr/8,qr}(x)\le C_1 \epsilon f^{pr/8,qr}(x).$
 \end{lemma}
 \pf
Note that
\begin{eqnarray*}
\wt{f}^{pr/8,qr}(x)&=&
\E_{x} \left[f(\wh{X}_{\wh{\tau}_{D^{pr/8}}}): \wh{X}_{\wh{\tau}_{D^{pr/8}}} \in
\overline{B}_{qr}\right]\\
&=& \int_{ \overline{B}_{qr}} \int_{D^{pr/8}} \wh{G}_{D^{pr/8}} (x,y) \wh{j}(y,z) m(dy) f(z) m(dz).
\end{eqnarray*}
By \textbf{C2}$(z_0, R)$,
$
\wh{j}(y,z) \le c_1 \wh{j}(y,z_0)$,  for all $(y,z) \in  B^c_{2qr} \times   \overline{B}_{qr},
$
where the constant $c_1$ is independent of $p$ and $q$. Thus
\begin{eqnarray*}
 &&\int_{ \overline{B}_{qr}} \int_{D^{pr/8}} \wh{G}_{D^{pr/8}} (x,y) \wh{j}(y,z) m(dy) f(z) m(dz)\\
 & \le & c_1 \int_{ \overline{B}_{qr}} \int_{D^{pr/8}} \wh{G}_{D^{pr/8}} (x,y) \wh{j}(y,z_0) m(dy) f(z) m(dz)\\
 & =& c_1 \wh{P}_{D^{pr/8}}(x,z_0) \int_{ \overline{B}_{qr}}  f(z) m(dz).
\end{eqnarray*}
Now, using \eqref{e:assumption-2-I} and the fact that $f=f^{pr/8,qr}$ on $D^{pr/8,qr}$,  we get that for every $x\in D^{pr}$,
 \begin{eqnarray*}
  \wt{f}^{pr/8,qr}(x) \le c_1 \epsilon\wh{P}_{D^{pr/8}}(x,z_0)
  \int_{ D^{pr/8,qr}} f^{pr/8,qr}(y) m(dy),
  \end{eqnarray*}
   which is less than or equal to
   $$c_1 \epsilon \wh{P}_{D^{pr/8}}(x,z_0)
  \int_{ B_{3pr/8}} f^{pr/8,qr}(y) m(dy).
   $$
Since  ${f}^{pr/8,qr}$ is regular harmonic in $D^{pr/8}$ with respect to $\wh{X}$, and vanishes on
$\overline{B}_{pr/8}^c\cap(\overline{D}^c\cup \wh{D}^{\mathrm{reg}})$,
using  {\bf F2}$(z_0, R)$ (with $a=3/2$), we conclude that for every $x\in D^{pr}$,
 \begin{eqnarray*}
  \wt{f}^{pr/8,qr}(x)  \le
  c_1 \epsilon \wh{P}_{D^{pr/8}}(x,z_0)
  \int_{ B_{3pr/8}} f^{pr/8,qr}(y) m(dy)
  \le  c_1C(3/2) \epsilon f^{pr/8,qr}(x).
 \end{eqnarray*}
 \qed

Again, by enlarging $C(3/2)$ in {\bf F2}, without loss of generality we assume $C_1=(C(3/2))^2$. From now on we let $C=C(3/2)$, so that $C_1=C^2$.

Let $r\ge R$ and $D\subset \overline{B}_r^c$ be an open set. Recall that for any $p>q>0$, $D^p=D\cap \overline{B}_p^c$ and $D^{p,q}=D^q\setminus D^p$. If $f_1$ and $f_2$ are non-negative functions on $\X$, for any $p>1$, we let
 $$
 m^{pr}:=\inf_{D^{pr}}\frac{f_1}{f_2}, \qquad  M^{pr}:=\sup_{D^{pr}}\frac{f_1}{f_2}.
 $$
Note that $f_i=f^{pr,qr}_i+ \wt{f}^{pr,qr}_i$.

 \begin{lemma}\label{l:oscillation-I}
 Let $r \ge R$, $D\subset \overline{B}_r^c$ an open set, and $p/16>q>2$. If $f_1$ and $f_2$ are non-negative functions on $\X$ which are regular harmonic in $D$ with respect
to $\widehat{X}$, and vanish  on
$\overline{B}_r^c\cap (\overline{D}^c\cup \wh{D}^{\mathrm{reg}})$,
then
 \begin{equation}\label{e:oscillation-I}
 (C+1)\left(\sup_{D^{pr}}\frac{f_1^{pr/8,qr}}{f_2^{pr/8,qr}}-\inf_{D^{pr}}\frac{f_1^{pr/8,qr}}{f_2^{pr/8,qr}}\right) \le (C-1)\left(M^{qr}-m^{qr}\right) .
 \end{equation}
  \end{lemma}
 \pf For any $x\in D^{pr/8}$, we define
 \begin{eqnarray*}
 g(x)&:=&f_1^{pr/8,qr} (x)-m^{qr} f_2^{pr/8,qr} (x)\\
  &=&\E_x \left[ (f_1-m^{qr} f_2) (\widehat{X}_{\widehat{\tau}_{D^{  pr/8  }   }}):
  \widehat{X}_{\widehat{\tau}_{D^{  pr/8  }   }} \in D^{pr/8,qr}
\right],
 \end{eqnarray*}
which is regular harmonic in $D^{pr/8}$ with respect to $\widehat{X}$, and vanishes  on
$\overline{B}_{pr/8}^c\cap ( \overline{D}^c\cup \wh{D}^{\mathrm{reg}})$.
Next, it follows from  \eqref{e:bhp-I}  that for any $x_1,x_2\in D^{pr}$ (we assume that $D^{pr}\neq \emptyset$),
$
g(x_1)f_2^{pr/8,qr}(x_2) \le C g(x_2)f_2^{pr/8,qr}(x_1).
$
 Therefore,
 \begin{eqnarray}\label{e:estimate-for-g-I}
 \sup_{D^{pr}}\frac{f_1^{pr/8,qr}}{f_2^{pr/8,qr}}-m^{qr}&=&  \sup_{D^{pr}}\frac{g}{f_2^{pr/8,qr}}
 \le  C \inf_{D^{pr}}\frac{g}{f_2^{pr/8,qr}}\\
& =&C\left(\inf_{D^{pr}}\frac{f_1^{pr/8,qr}}{f_2^{pr/8,qr}}-m^{qr}\right).
\nonumber
 \end{eqnarray}
 We can similarly get that
  \begin{equation}\label{e:estimate-for-h-I}
  M^{qr}-\inf_{D^{pr}}\frac{f_1^{pr/8,qr}}{f_2^{pr/8,qr}} \le C\left(M^{qr}-\sup_{D^{pr}}\frac{f_1^{pr/8,qr}}{f_2^{pr/8,qr}}\right).
  \end{equation}
 Adding \eqref{e:estimate-for-g-I} and \eqref{e:estimate-for-h-I} and rearranging, we arrive at \eqref{e:oscillation-I}. \qed

  \bigskip

 For any positive function $\phi$ on a non-empty open set $U$, let
\begin{align}
\label{e:RO}
 \mathrm{RO}_U \phi=\frac{\sup_U \phi}{\inf_U \phi}.
\end{align}

\begin{lemma}\label{l:oscillation-assumption-2-I}
Let $r \ge R$, $D\subset \overline{B}_r^c$ an open set, $p/16>q>2$, and $\epsilon >0$. Let $f_1$ and $f_2$ be non-negative functions on $\X$ which are regular harmonic in $D$ with respect to $\widehat{X}$ and vanish  on
$\overline{B}_r^c\cap ( \overline{D}^c\cup \wh{D}^{\mathrm{reg}})$.
If
 \begin{equation}\label{e:assumption-2-12-I}
\int_{ \overline{B}_{qr}} f_i(y) m(dy) \le \epsilon
\int_{ D^{pr/8,qr}} f_i(y) m(dy), \quad i=1,2,
 \end{equation}
 then
 \begin{equation}\label{e:oscillation-assumption-2-I}
    \mathrm{RO}_{D^{pr}} \frac{f_1}{f_2} \le (1+C^2\epsilon)^2+(1+C^2\epsilon)\frac{C-1}{C+1}\left(\mathrm{RO}_{D^{qr}}\frac{f_1}{f_2}-1\right).
 \end{equation}
\end{lemma}
\pf 
Applying Lemma \ref{l:assumption-2-I} we get that
\begin{eqnarray*}
M^{pr}&=&\sup_{D^{pr}}\frac{f_1}{f_2}=\sup_{D^{pr}}\frac{f_1^{pr/8,qr}+\wt{f}_1^{pr/8,qr}}{f_2^{pr/8,qr}+\wt{f}_2^{pr/8,qr}}
\le \sup_{D^{pr}}\frac{(1+C^2\epsilon)f_1^{pr/8,qr}}{f_2^{pr/8,qr}},\\
m^{pr}&=&\inf_{D^{pr}}\frac{f_1}{f_2}=\inf_{D^{pr}}\frac{f_1^{pr/8,qr}+\wt{f}_1^{pr/8,qr}}{f_2^{pr/8,qr}+\wt{f}_2^{pr/8,qr}}
\ge \inf_{D^{pr}}\frac{f_1^{pr/8,qr}}{(1+C^2\epsilon)f_2^{pr/8,qr}}.
\end{eqnarray*}
Inserting this in \eqref{e:oscillation-I}, we arrive at
\begin{align*}
&(C+1)\left(\frac{M^{pr}}{1+C^2\epsilon}-(1+C^2\epsilon)m^{pr}\right)\\ 
\le &\ (C+1)\left(\sup_{D^{pr}}\frac{f_1^{pr/8,qr}}{f_2^{pr/8,qr}}-\inf_{D^{pr}}\frac{f_1^{pr/8,qr}}{f_2^{pr/8,qr}}\right)\\
\le &\ (C-1)(M^{qr}-m^{qr}).
\end{align*}
Rearranging and using that $m^{pr}\ge m^{qr}$ we get
$$
\frac{M^{pr}}{m^{pr}}\le (1+C^2\epsilon)^2+(1+C^2\epsilon)\frac{C-1}{C+1}\left(\frac{M^{qr}}{m^{qr}}-1\right),
$$
which implies \eqref{e:oscillation-assumption-2-I}. \qed

In the rest of this subsection, we fix an open set $D$ and a point $x_0\in D$.

\begin{lemma}\label{l:new_estpq}
Suppose that $\infty$ is accessible from  $D$ with respect to $X$.
For any $q\ge 4$, $r\ge 2d(z_0,x_0) \vee R$ and $\epsilon>0$,
 there exists $p=p(\epsilon,q, D, x_0, r)>16q$ such that
$$
\int_{D^{pr,qr}} G_D(x_0, y)m(dy)
> \epsilon \int_{B_{q r}} G_D(x_0, y)m(dy).
$$
 \end{lemma}
\pf
Since
$\infty$ is accessible from  $D$ with respect to $X$,  we have that
\begin{eqnarray*}
\E_{x_0} \tau_D =\int_{D} G_D(x_0, v)m(dv)=\infty.
\end{eqnarray*}
The function $v\mapsto G_D(x_0,v)$ is regular harmonic in $D^r \supset D^{qr/3}$ with respect to $\widehat{X}$ and vanishes on
$\overline{B}^{qr/3}\cap ( \overline{D^{qr/3}}^c\cup \wh{(D^{qr/3})}^{\mathrm{reg}})$.
By using  {\bf F2}$(z_0, R)$ for $\widehat{X}$ (with the open set $D^{qr/3}$, $a=3/2$ and radius $qr/3$)
$$
 \int_{B_{qr}}G_D(x_0,z)m(dz) \le c \inf_{v\in D^{8qr/3}} \frac{ G_D(x_0,v)}{\widehat{P}_{D^{qr/3}} (v,z_0)} < \infty.
$$
Thus
$$
\infty= \int_{D^{qr}}G_D(x_0,z)m(dz)=\lim_{p \to \infty}
\int_{D^{pr,qr}} G_D(x_0, z)m(dz)
$$
and so we can choose
$p=p(\epsilon,q, D, x_0, r)>16q$ large enough so that
$$
\int_{D^{pr,qr}} G_D(x_0, z)m(dz)
> \epsilon \int_{B_{q r}} G_D(x_0, z)m(dz).
$$
\qed

\begin{prop}\label{p:oscillation-reduction-I}
Suppose that $\infty$ is accessible from  $D$ with respect to $X$.
Let $r > 2d(z_0,x_0) \vee R$.
For every $\eta>0$, there exists $s=s(r, D, x_0, \eta)>1$ such that for any two non-negative functions $f_1$, $f_2$ on $\X$ which are regular harmonic in $D^r$ with respect to $\widehat{X}$ and vanish  on
$\overline{B}_r^c\cap ( \overline{D}^c\cup \wh{(D^r)}^{\mathrm{reg}})$,
we have
\begin{equation}\label{e:oscillation-reduction-I}
\mathrm{RO}_{D^{sr}}\frac{f_1}{f_2}\le 1+\eta \, .
\end{equation}

\end{prop}
\pf
Let $\eta>0$ and define
\begin{align}\label{e:phi1}
\phi(t):=1+\frac{\eta}{2}+\frac{C}{C+1}(t-1), \quad t\ge 1.
\end{align}
Then $\phi(t)=t$ for $t=1+\eta(C+1)/2$,
$\phi(t)<t$ if $t>1+\eta(C+1)/2$, and $\phi(t)>t$ if $t<1+\eta(C+1)/2$.
Thus $\lim_{l\to \infty}\phi^l(C)=1+\eta(C+1)/2$, where $\phi^l$ is the $l$-fold composition of $\phi$.
Let $l\in \N$ be such that
\begin{align}\label{e:phi2}
\phi^l(C)<1+\eta(C+1).
\end{align}
Choose $\epsilon=\epsilon(\eta)>0$ small enough so that
\begin{align}\label{e:phi3}
(C\epsilon+1+\epsilon)^2(1+\epsilon)^2 < 1+\eta
\end{align}
and
\begin{align}\label{e:repphi}
(1+C^2\epsilon)^2+(1+C^2\epsilon)\frac{C-1}{C+1}(t-1)<1+\frac{\eta}{2}+\frac{C}{C+1}(t-1)=\phi(t)
\end{align}
for all $t\ge 1$.
Let $k$ be the smallest integer such that $k>C^2\epsilon^{-2}$ and denote $n=lk$.
Let $q_0=8$ and  choose $q_{1}=p(\epsilon,q_0, r, D, x_0)$  as in Lemma \ref{l:new_estpq}.
Inductively, using  Lemma \ref{l:new_estpq}, we choose
$q_{j+1}=p(\epsilon,q_j,r, D, x_0)$ for $j=0,1,\dots, n-1$, and $s=q_n$.
 Then by Lemma \ref{l:new_estpq},
 for $j=0,1,\dots, n-1$, we have
\begin{equation}\label{e:1-or-2-I}
\int_{D^{q_{j+1}r, q_j r}} G_D(x_0, y)m(dy)
> \epsilon \int_{\overline{B}_{q_j r}} G_D(x_0, y)m(dy).
\end{equation}
It follows from {\bf F2}$(z_0, R)$ (applied to the open set $D^{q_j r/3}$ with $a=3/2$ and radius $q_j r/3$) that for every $j=0,1,\dots,n-1$, $i=1, 2$
and $x\in D^{q_{j+1}r, 8q_j r/3}$,
$$
C\frac{f_i(x)}{\int_{B_{q_j r}} f_i(y)m(dy)} \ge \widehat{P}_{D^{q_j r/3}}(x,z_0)
\ge C^{-1}\frac{G_D(x_0, x)}
{ \int_{B_{q_j r}} G_D(x_0, y)m(dy)}.
$$
Hence, by integrating over $D^{q_{j+1}r, 8q_j r/3}$ we get
$$
\frac{\int_{D^{q_{j+1}r, 8q_j r/3}} f_i(x) m(dx)}{\int_{B_{q_j r}} f_i(y)m(dy)}
\ge C^{-1}\frac{\int_{D^{q_{j+1}r, 8q_j r/3}}G_D(x_0, x)m(dx)}
{ \int_{B_{q_j r}} G_D(x_0, y)m(dy)},
\quad i=1,2.
$$
Together with \eqref{e:1-or-2-I} we get that
$$\int_{
D^{q_{j+1}r,q_j r}}f_i(y)m(dy) \ge
\int_{
D^{q_{j+1}r,8q_j r/3}}f_i(y)m(dy) >C^{-2}\epsilon
\int_{B_{q_j r}}f_i(y)m(dy)
$$
for \emph{both} $i=1$ and $i=2$, and all $j=0,1,\dots n-1$. Let $0\le m <l$. By the definition of $k$,
\begin{align*}
&\int_{D^{q_{(m+1)k}r/8, q_{mk}r}}
f_i(y)m(dy)
 \ge  \int_{D^{q_{(m+1)k-1}r, q_{mk}r}}f_i(y)m(dy)\\
 =&\sum_{j=mk }^{(m+1)k-1}
\int_{D^{q_{j+1}r, q_j r}}
f_i(y)m(dy)
\ge \  k C^{-2}\epsilon \int_{B_{q_{mk}r}}f_i(y)m(dy)\\
\ge &\ \epsilon^{-1}\int_{B_{q_{mk}r}}f_i(y)m(dy), \quad i=1,2.
\end{align*}
By using Lemma \ref{l:oscillation-assumption-2-I} with $p=q_{(m+1)k}$ and $q=q_{mk}$ we conclude from
\eqref{e:repphi} that for every integer $m$ such that $0\le m<l$,
\begin{align*}
  \mathrm{RO}_{D^{q_{(m+1)k}r}} \frac{f_1}{f_2} &\le (1+C^2\epsilon)^2+(1+C^2\epsilon)\frac{C-1}{C+1}\left(\mathrm{RO}_{D^{q_{mk}r}}\frac{f_1}{f_2}-1\right)\\
 & <\phi\left(\mathrm{RO}_{D^{q_{mk}r}}\frac{f_1}{f_2}\right).
\end{align*}
By definition of the integer $l$, monotonicity of $\phi$, and the fact that
$\mathrm{RO}_{D^{r/2}}(f_1/f_2)\le C$, it follows that
$$
\mathrm{RO}_{D^{q_{lk}R}}\frac{f_1}{f_2}\le \phi\left(\mathrm{RO}_{D^{q_{(l-1)k}R}}\right)\le \cdots \le \phi^l\left(\mathrm{RO}_{D^{q_{0}R}}\right)\le 1 +\eta(C+1).
$$
This means that
$
\mathrm{RO}_{D^{sr}}\frac{f_1}{f_2}\le 1+\eta(C+1).
$
Since $\eta>0$ is arbitrary, the proof is complete. \qed

\begin{corollary}\label{c:oscillation-reduction-I}
Suppose that $\infty$ is accessible from  $D$ with respect to $X$.
Let $r >  R$
and let $f_1$ and $f_2$ be non-negative functions on $\X$ which are regular harmonic in $D^r$ with respect to $\widehat{X}$ and vanish on
$\overline{B}_r^c\cap (\overline{D}^c\cup \wh{(D^r)}^{\mathrm{reg}})$.
Then the limit
$$
\lim_{D\ni x\to \infty}\frac{f_1(x)}{f_2(x)}
$$
exists and is finite.
\end{corollary}

\pf
Since one can increase $r$ so that $r> 2d(z_0,x_0) \vee R$ without loss of generality, the existence of the limit and its finiteness is a direct consequence of Proposition \ref{p:oscillation-reduction-I}.
\qed

\subsection{Finite boundary point}\label{s:os}

In this subsection, we deal with the oscillation reduction at a boundary point $z_0$ of an open set $D$.
Throughout the subsection, we assume that there exists $R\le R_0$ such that \textbf{C1}$(z_0, R)$ holds,
and that $\widehat{X}$ satisfies \textbf{F1}$(z_0, R)$. We will see that the results and the estimates below  have the same structure as those in the previous subsection,
the difference being that $\wh{P}_D(x,z_0)$ is replaced by $\E_x \wh{\tau}_D$ and $\int_{B_{ar}}f(y)m(dy)$
is replaced by $\wh{\Lambda}_{ar}(f)$ (see below for definition).

As in the previous subsection we begin by recording a simple consequences of \textbf{F1}$(z_0, R)$ for $\wh{X}$,
the boundary Harnack principle in \cite{KSVp1}:
There exists $c>1$ such that  for any $r\in (0, R)$, any open set $D\subset B_r$
and any non-negative functions $u$ and $v$ on $\X$ that are regular harmonic in $D$
with respect to $\wh{X}$ and vanish on
$B_r \cap ( \overline{D}^c\cup \wh{D}^{\mathrm{reg}})$,
it holds that
\begin{equation}\label{e:bhp}
c^{-1}\frac{u(y)}{v(y)} \le \frac{u(x)}{v(x)} \le c \frac{u(y)}{v(y)}\qquad \text{for all } x,y\in D\cap B_{r/8} .
\end{equation}
Note that we can take $c=(C(2/3))^2$. By enlarging $C(2/3)$ in {\bf F1}$(z_0, R)$, without loss of generality we assume the above $c$ is equal to  $C(2/3)$.

Let $D\subset \X$ be an open set.
For $0<p<q$, let $D_p=D\cap B_p$ and $D_{p,q}=D_q\setminus D_p$.
For a function $f$ on $\X$, and $0<p<q$, let
 \begin{equation}\label{e:Lambda}
 \wh{\Lambda}_p(f):=\int_{\overline{B}_p^c}\wh{j}(z_0,y) f(y) m(dy),\qquad \wh{\Lambda}_{p,q}(f):=\int_{D_{p,q}}\wh{j}(z_0,y) f(y) m(dy).
 \end{equation}

Let $r\in (0, R]$. For $0<p<q<1$ and a non-negative function $f$ on $\X$ define
 \begin{eqnarray*}
 f_{pr,qr}(x)&=&
  \E_{x} \left[f(\wh{X}_{\wh{\tau}_{D_{pr}}}): \wh{X}_{\wh{\tau}_{D_{pr}}} \in D_{pr,qr}\right],\\
 \wt{f}_{pr,qr}(x)&=&
 \E_{x} \left[f(\wh{X}_{\wh{\tau}_{D_{pr}}}): \wh{X}_{\wh{\tau}_{D_{pr}}} \in (D\setminus D_{qr})\cup B_r^c\right].
 \end{eqnarray*}

 \begin{lemma}\label{l:assumption-2}
Suppose that $r\le R$, $D$ is an open subset of  $B_r$ and
$f$ is a non-negative function on $\X$
that is regular harmonic in $D$ with respect to $\wh{X}$ and vanishes on
$B_r \cap ( \overline{D}^c\cup \wh{D}^{\mathrm{reg}})$.
There exists $C_2>0$ independent of $f$ and $r\le R$ such that
for any $0<16p<q<1/2$ and any $\epsilon>0$,  if
 \begin{equation}\label{e:assumption-2}
\wh{\Lambda}_{qr}(f)\le \epsilon \wh{\Lambda}_{8pr, qr}(f),
 \end{equation}
 then for every $x\in D_{pr}$,
 $
 \wt{f}_{8pr,qr}(x)\le C_2 \epsilon f_{8pr,qr}(x).
 $
 \end{lemma}
 \pf
 Note that
 \begin{eqnarray*}
 \wt{f}_{8pr,qr}(x)&=&
 \E_{x} \left[f(\wh{X}_{\wh{\tau}_{D_{8pr}}}): \wh{X}_{\wh{\tau}_{D_{pr}}} \in
 {B}_{qr}^c\right]\\
 &=& \int_{{B}_{qr}^c} \int_{D_{8pr}} \wh{G}_{D_{8pr}} (x,y) \wh{j}(y,z) m(dy) f(z) m(dz).
\end{eqnarray*}
By \textbf{C1}$(z_0, R)$,
$
\wh{j}(y,z) \le c_1 \wh{j}(z_0,z)$  for all  $(y,z) \in  \overline{B}_{qr/2} \times   {B}_{qr}^c,
$
where the constant $c$ is independent of $p$ and $q$. Thus
\begin{align*}
 &\int_{{B}_{qr}^c} \int_{D_{8pr}} \wh{G}_{D_{8pr}} (x,y) \wh{j}(y,z) m(dy) f(z) m(dz)\\
 \le& c_1  \int_{D_{8pr}} \wh{G}_{D_{8pr}} (x,y) m(dy)
   \int_{{B}_{qr}^c} \wh{j}(z_0, z) f(z) m(dz)\\
 = & c_1 \E_x \wh{\tau}_{D_{8pr}}    \int_{{B}_{qr}^c}  \wh{j}(z_0, z) f(z) m(dz)
\end{align*}

 Now, using \eqref{e:assumption-2} and the fact that $f=f_{8pr,qr}$ on $D_{8pr,qr}$,
 we get that for every $x\in D_{pr}$,
$$
   \wt{f}_{8pr,qr}(x) \le  c_1 \epsilon (\E_x \wh{\tau}_{D_{8pr}}) \wh{\Lambda}_{8pr,qr}(f_{8pr,qr}),
 $$
 which is less than or equal to
$ c \epsilon (\E_x \wh{\tau}_{D_{8pr}}) \wh{\Lambda}_{8pr/3}(f_{8pr,qr})$.
 Note that  ${f}_{8pr,qr}$ is regular harmonic in $D_{8pr}$ with respect to
 $\wh{X}$ and vanishes on
$B_{8pr}\cap ( \overline{D}^c\cup \wh{D}^{\mathrm{reg}})$.
Thus applying  \textbf{F1}$(z_0, R)$ (with $a=2/3$) to ${f}_{8pr,qr}$ we have that  for every $x\in D_{pr}$,
  \begin{eqnarray*}
   \wt{f}_{8pr,qr}(x) \le c_1 \epsilon (\E_x \wh{\tau}_{D_{8pr}}) \wh{\Lambda}_{8pr/3}(f_{8pr,qr})
   \le c_1C(2/3)\epsilon {f}_{8pr,qr}(x).
 \end{eqnarray*}
  \qed

Again, by enlarging $C(2/3)$ in {\bf F1}, without loss of generality we can assume $C_2=(C(2/3))^2$. From now on we let $C=C(2/3)$, so that $C_2=C^2$.

Let $r\in (0,R]$, $D\subset B_r=B(z_0,r)$ an open set and $z_0\in \partial D$.
Recall that for $0<p<q$,  $D_p=D\cap B_p$ and $D_{p,q}=D_q\setminus D_p$.
If $f_1$ and $f_2$ are non-negative functions on $\X$,  for any $p\in (0,1)$, we let
 $$
 m_{pr}:=\inf_{D_{pr}}\frac{f_1}{f_2}, \qquad  M_{pr}=\sup_{D_{pr}}\frac{f_1}{f_2}.
 $$
Note that $f_i=(f_i)_{pr,qr}+(\wt{f}_i)_{pr,qr}$.

 \begin{lemma}\label{l:oscillation}
 Let $r\le R$, $D\subset B_r$ an open set, and $0<16p<q<1/2$.  If $f_1$ and $f_2$ are non-negative functions on $\X$ which are regular harmonic in $D$ with respect
to $\widehat{X}$ and vanish  on
$B_r\cap ( \overline{D}^c\cup \wh{D}^{\mathrm{reg}})$,
then
 \begin{equation}\label{e:oscillation}
 (C+1)\left(\sup_{D_{pr}}\frac{(f_1)_{8pr,qr}}{(f_2)_{8pr,qr}}
 -\inf_{D_{pr}}\frac{(f_1)_{8pr,qr}}{(f_2)_{8pr,qr}}\right) \le (C-1)\left(M_{qr}-m_{qr}\right) .
 \end{equation}
  \end{lemma}
 \pf For  any $x\in D_{8pr}$, we define
\begin{align*}
 g(x)&:=(f_1)_{8pr,qr} (x)-m_{qr} (f_2)_{8pr,qr} (x)\\
  &=\E_x \left[ (f_1-m_{qr} f_2) (\widehat{X}_{\widehat{\tau}_{D_{8pr }   }}):
  \widehat{X}_{\widehat{\tau}_{D_{ 8pr  }   }} \in D_{8pr,qr} \right],
\end{align*}
 which    is regular harmonic in $D_{8pr}$ with respect to $\widehat{X}$ and vanishes on
$B_{8pr}\cap (\overline{D}^c\cup \wh{D}^{\mathrm{reg}})$.
 Next, it follows from \eqref{e:bhp}  that for any $x_1,x_2\in D_{pr}$
 (we assume that $D_{pr}\neq \emptyset$),
\begin{eqnarray*}
g(x_1)(f_2)_{8pr,qr}(x_2)\le C g(x_2)(f_2)_{8pr,qr}(x_1).
\end{eqnarray*}
 Therefore,
\begin{eqnarray}\label{e:estimate-for-g}
 \sup_{D_{pr}}\frac{(f_1)_{8pr,qr}}{(f_2)_{8pr,qr}}-m_{qr}&=&  \sup_{D_{pr}}\frac{g}{(f_2)_{8pr,qr}}
 \le  C \inf_{D_{pr}}\frac{g}{(f_2)_{8pr,qr}}\\
& =&C\left(\inf_{D_{pr}}\frac{(f_1)_{8pr,qr}}{(f_2)_{8pr,qr}}-m_{qr}\right).\nonumber
 \end{eqnarray}
 We can similarly get that
  \begin{equation}\label{e:estimate-for-h}
  M_{qr}-\inf_{D_{pr}}\frac{(f_1)_{8pr,qr}}{(f_2)_{8pr,qr}} \le C\left(M_{qr}-\sup_{D_{pr}}\frac{(f_1)_{8pr,qr}}{(f_2)_{8pr,qr}}\right).
  \end{equation}
  Adding \eqref{e:estimate-for-g} and \eqref{e:estimate-for-h} and rearranging, we arrive at \eqref{e:oscillation}. \qed

  \bigskip

Recall that $\mathrm{RO}_U \phi$ is defined in   \eqref{e:RO}.
\begin{lemma}\label{l:oscillation-assumption-2}
 Let $r \le R$, $D\subset B_r^c$ an open set, $16p<q<1/2$, and $\epsilon >0$. Let $f_1$ and $f_2$ be non-negative functions on $\X$ which are regular harmonic in $D$ with respect to $\widehat{X}$ and vanish  on
$B_r\cap ( \overline{D}^c\cup \wh{D}^{\mathrm{reg}})$.
 If
$
 \wh{\Lambda}_{qr}(f_i)\le \epsilon \wh{\Lambda}_{8pr,qr}(f_i)$, $i=1,2$,
 then
 \begin{equation}\label{e:oscillation-assumption-2}
    \mathrm{RO}_{D_{pr}} \frac{f_1}{f_2} \le (1+C^2\epsilon)^2+(1+C^2\epsilon)\frac{C-1}{C+1}\left(\mathrm{RO}_{D_{qr}}\frac{f_1}{f_2}-1\right).
 \end{equation}
\end{lemma}
\pf Applying Lemma \ref{l:assumption-2}, we get that
\begin{eqnarray*}
M_{pr}&=&\sup_{D_{pr}}\frac{f_1}{f_2}=\sup_{D_{pr}}\frac{(f_1)_{8pr,qr}
+(\wt{f}_1)_{8pr,qr}}{(f_2)_{8pr,qr}+(\wt{f}_2)_{8pr,qr}} \le \sup_{D_{pr}}\frac{(1+C^2\epsilon)(f_1)_{8pr,qr}}{(f_2)_{8pr,qr}},\\
m_{pr}&=&\inf_{D_{pr}}\frac{f_1}{f_2}=\inf_{D_{pr}}\frac{(f_1)_{8pr,qr}
+(\wt{f}_1)_{8pr,qr}}{(f_2)_{8pr,qr}+(\wt{f}_2)_{8pr,qr}}
\ge \inf_{D_{pr}}\frac{(f_1)_{8pr,qr}}{(1+C^2\epsilon)(f_2)_{8pr,qr}}.
\end{eqnarray*}
By inserting this in \eqref{e:oscillation}, we arrive at
\begin{align*}
&(C+1)\left(\frac{M_{pr}}{1+C^2\epsilon}-(1+C^2\epsilon)m_{pr}\right) \\
\le &\ (C+1)\left(\sup_{D_{pr}}\frac{(f_1)_{8pr,qr}}{(f_2)_{8pr,qr}}-
\inf_{D_{pr}}\frac{(f_1)_{8pr,qr}}{(f_2)_{8pr,qr}}\right)\\
\le & \ (C-1)(M_{qr}-m_{qr}).
\end{align*}
Rearranging and using that $m_{pr}\ge m_{qr}$ we get
$$
\frac{M_{pr}}{m_{pr}}\le (1+C^2\epsilon)^2+(1+C^2\epsilon)\frac{C-1}{C+1}\left(\frac{M_{qr}}{m_{qr}}-1\right),
$$
which implies \eqref{e:oscillation-assumption-2}. \qed

In the remainder of this subsection, we fix an open set $D$ such that $z_0\in \partial D$,  and a point $x_0\in D$.

\begin{lemma}\label{l:new_estpqf}
Suppose that $z_0$ is accessible from  $D$ with respect to $X$.
Assume that $r \le
 R \wedge (\frac12 d(z_0,x_0))
$, $q\le1/4$ and $\epsilon>0$.
Then there exists $p=p(\epsilon,q, D, x_0, r)<q/16$ such that
  $ \wh{\Lambda}_{pr,qr}(G_D(x_0, \cdot))
  > \epsilon \wh{\Lambda}_{qr}(G_D(x_0, \cdot)).
$
 \end{lemma}
\pf Since $z_0$ is accessible from  $D$ with respect to $X$,  we have that
\begin{eqnarray*}
P_D(x_0,z_0)=\int_D G_D(x_0,v)j(v, z_0) m(dv) = \infty.
\end{eqnarray*}
The function $v\mapsto G_D(x_0,v)$ is regular harmonic in $D_r\supset D_{3qr}$ with respect to $\widehat{X}$ and vanishes on
$B_{3qr}\cap ( \overline{D_{3qr}}^c\cup \wh{(D_{3qr})}^{\mathrm{reg}})$.
By using  {\bf F1}$(z_0, R)$ for $\widehat{X}$ (with the open set $D_{3qr}$, $a=2/3$ and radius $3qr$),
  \begin{align*}
\int_{\overline{B}_{qr}^c}G_D(x_0,y)j(y, z_0)m(dv)=\wh{\Lambda}_{qr}(G_D(x_0, \cdot))
 \le
 c  \inf_{v\in D_{3qr/8}} \frac{G_D(x_0,v)}{\E_{v}\widehat{\tau}_{D_{3qr}}} < \infty.
 \end{align*}
Thus
$$
\infty
= \int_{D_{qr}}G_D(x_0,v)j(v, z_0)m(dv)
=\lim_{p\to 0}\int_{D_{pr,qr}} G_D(x_0, v)j(v, z_0)m(dy)
$$
and so we can choose $p=p(\epsilon,q, D, x_0, r)<q/16$ small  so that
$$
\int_{D_{pr,qr}} G_D(x_0, y)j(y, z_0)m(dy) > \epsilon \int_{\overline{B}_{qr}^c} G_D(x_0, y)j(y, z_0)m(dy).
$$
\qed

\begin{prop}\label{p:oscillation-reduction}
 Suppose that $z_0$ is accessible from  $D$ with respect to $X$.
 Assume that $r  \le
R \wedge (\frac12 d(z_0,x_0))$.
 For every $\eta>0$, there exists $s=s(r, D, x_0, \eta) \in (0,1)$ such that for any two non-negative functions $f_1$, $f_2$ on $\X$ which are regular harmonic in $D_r$ with respect to $\widehat{X}$ and vanish  on
$B_r\cap ( \overline{D}^c\cup \wh{(D_r)}^{\mathrm{reg}})$,
 we have
\begin{equation}\label{e:oscillation-reduction}
\mathrm{RO}_{D_{sr}}\frac{f_1}{f_2}\le 1+\eta .
\end{equation}
\end{prop}
\pf
Let $\eta>0$ and define $\phi$ as in \eqref{e:phi1}
and let $\phi^l$ be the $l$-fold composition of $\phi$.
Let $l\in \N$ be such that \eqref{e:phi2} holds.
Choose $\epsilon>0$ small enough so that
\eqref{e:phi3}
and \eqref{e:repphi} holds.
Let $k$ be the smallest integer such that $k>C^2\epsilon^{-2}$ and denote $n=lk$.
Let $q_0=1/8$ and choose $q_{1}=p(\epsilon,q_0, r,  D, x_0)$  as in Lemma \ref{l:new_estpqf}.
Inductively, using  Lemma \ref{l:new_estpqf}, we choose
$q_{j+1}=p(\epsilon,q_j,r, D, x_0)$ for $j=0,1,\dots, n-1$, and $s=q_n$.
Then it follows from Lemma \ref{l:new_estpqf} that
 for $j=0,1,\dots, n-1$, we have
\begin{equation}\label{e:1-or-2}
\wh{\Lambda}_{q_{j+1}r,q_j r}(G_D(x_0, \cdot)) > \epsilon \wh{\Lambda}_{q_j r}(G_D(x_0, \cdot)).
\end{equation}
It follows from {\bf F1}$(z_0, R)$ (applied to the open set $D_{3q_j r}$ with $a=2/3$ and radius $3q_j r$)  that for every $j=0,1,\dots,n-1$,
$$
C\frac{f_i(x)}{\wh{\Lambda}_{q_j r}(f_i)} \ge \E_x \widehat{\tau}_{D_{3q_j r}} \ge
C^{-1}\frac{G_D(x_0, x)}{\wh{\Lambda}_{q_j r}(G_D(x_0, \cdot))},
\quad x\in D_{q_{j+1}r, 3q_j r/8}.
$$
Hence, by integrating over $D_{q_{j+1}r, 3q_j r/8}$ we get
$$
\frac{\wh{\Lambda}_{q_{j+1}r, 3q_j r/8}(f_i)}{\wh{\Lambda}_{q_j r}(f_i)}
\ge C^{-1} \frac{\wh{\Lambda}_{q_{j+1}r, 3q_j r/8}(G_D(x_0, \cdot))}{\wh{\Lambda}_{q_j r}(G_D(x_0, \cdot))} ,
\quad i=1,2.
$$
Together with \eqref{e:1-or-2} it follows that
$
\wh{\Lambda}_{q_{j+1}r, q_j r}(f_i) \ge
\wh{\Lambda}_{q_{j+1}r, 3q_j r/8}(f_i)
>C^{-2}\epsilon \wh{\Lambda}_{q_j r}(f_i)
$
for \emph{both} $i=1$ and $i=2$, and all $j=0,1,\dots n-1$. Let $0\le m <l$; then
\begin{eqnarray*}
\Lambda_{8q_{(m+1)k}r, q_{mk}r}(f_i)& \ge& \wh{\Lambda}_{q_{(m+1)k-1}r, q_{mk}r}(f_i)
=\sum_{j=mk }^{(m+1)k-1}\wh{\Lambda}_{q_{j+1}r, q_j r}(f_i)\\
&\ge & k C^{-2}\epsilon \wh{\Lambda}_{q_{mk}r}(f_i)\ge \epsilon^{-1}\wh{\Lambda}_{q_{mk}r}(f_i), \quad i=1,2.
\end{eqnarray*}
By using Lemma \ref{l:oscillation-assumption-2} with $p=q_{(m+1)k}$ and $q=q_{mk}$
we conclude that for every integer $m$ such that $0\le m<l$,
\begin{align*}
  \mathrm{RO}_{D_{q_{(m+1)k}r}} \frac{f_1}{f_2} &\le (1+C^2\epsilon)^2+(1+C^2\epsilon)\frac{C-1}{C+1}\left(\mathrm{RO}_{D_{q_{mk}r}}
  \frac{f_1}{f_2}-1\right)\\
  &<\phi\left(\mathrm{RO}_{D_{q_{mk}r}}\frac{f_1}{f_2}\right).
\end{align*}
The remainder of the proof is the same as the corresponding part of the proof of Proposition \ref{p:oscillation-reduction-I}.
 \qed

\bigskip

\begin{corollary}\label{c:reductionac}
 Suppose that $z_0$ is accessible from  $D$ with respect to $X$.
Let $r  \le
R
$ and let $f_1$ and $f_2$ be non-negative functions on $\X$ which are regular harmonic in $D_r$ with respect to $\widehat{X}$ and vanish  on
$B_r\cap (\overline{D}^c\cup \wh{(D_r)}^{\mathrm{reg}})$.
 Then the limit
$$
\lim_{D\ni x\to z_0}\frac{f_1(x)}{f_2(x)}
$$
exists and is finite.
\end{corollary}

\pf
Since one can decrease $r$ so that $r \le  R \wedge (\frac12 d(z_0,x_0))$ without loss of generality, the existence of the limit is a direct consequence of Proposition \ref{p:oscillation-reduction}. \qed

\section{Martin boundary for accessible points }\label{s:mb-infty}

 Recall that $D$ is a Greenian open subset of $\X$ and that $X^D$ is the process $X$ killed upon exiting $D$. In order to apply the theory of Martin boundary developed in \cite{KW}, we have to check that their Hypothesis (B), 
see p.498 in \cite{KW}, holds in our setting. 
 Since $X^D$ is strongly Feller, it follows by the dominated convergence theorem that the $\alpha$-resolvent operator $G_D^{\alpha}f(x)=\int_0^{\infty} e^{-\alpha t}P_t^D f(x)\, dt$, $\alpha>0$, is also strongly Feller. Here $(P_t^D)_{t\ge 0}$ denotes the semigroup of $X^D$. In particular, $G_D^{\alpha}f$ is continuous for every bounded non-negative measurable $f$ on $D$. It follows that $G_D^{\alpha}f$ is lower semi-continuous for every non-negative $f$ on $D$ and every $\alpha>0$. Since $G_D f=\, \uparrow \lim_{\alpha\to 0}G_D^{\alpha}f$, we see that $G_D f$ is also lower semi-continuous for every  non-negative $f$. 
Hence, conditions (11) and (12) on p.126 of \cite{CW} are satisfied. 
It follows from Theorem 2 on  p.268 of \cite{CW} that  $G_D\ind_K$ 
 is bounded for every compact set $K\subset D$. Let $f:D\to [0,\infty)$ be bounded measurable and vanish outside of a compact set $K\subset D$. Then 
 $0\le f\le \|f\|_\infty \ind_K$. Thus $G_Df \le G_D(\|f\|_\infty \ind_K)\le \|f\|_\infty G_D\ind_K$ is bounded. Since $X^D$ is strongly Feller, it follows that $P^D_t G_D f$ is continuous. Further,
$$
G_D f-P^D_t G_D f=\int_0^t P^D_s f ds \le 
\|f\|_\infty \int_0^t P^D_s 1 ds \le \|f\|_\infty t \, .
$$
The right-hand side converges to 0 uniformly in $x\in D$. Hence
$$
G_D f=\lim_{t\to 0} P^D_t G_D f
$$
uniformly in $D$. Thus $G_D f$ is a uniform limit of continuous functions, hence continuous. Finally, if $f\in C_c(D)$ (continuous functions on $D$ with compact support), it is clear that $\alpha G_D^{\alpha}f(x)=\E_x \int_0^{\infty}e^{-t}f(X_{t/\alpha})\, dt \to f(x)$ boundedly on compacts as $\alpha\to \infty$. 
Since the same conclusions are valid for $\widehat{X}$,
we have checked that Hypothesis (B) from \cite{KW} holds true.

Fix $x_0\in D$ and define
$$
M_D(x, y):=\frac{G_D(x, y)}{G_D(x_0, y)}, \qquad x, y\in D,~y\neq x_0.
$$
By Theorem 3 in \cite{KW}, $D$  has
a Martin boundary $\partial_M D$ with respect to $X^D$ satisfying the following properties:
\begin{description}
\item{(M1)} $D\cup \partial_M D$ is
a compact metric space (with the metric denoted by $d_M$);
\item{(M2)} $D$ is open and dense in $D\cup \partial_M D$,  and its relative topology coincides with its original topology;
\item{(M3)}  $M_D(x ,\, \cdot\,)$ can be uniquely extended  to $\partial_M D$ in such a way that
\begin{description}
\item{(a)}
$ M_D(x, y) $ converges to $M_D(x, w)$ as $y\to w \in \partial_M D$ in the Martin topology;
\item{(b)} for each $ w \in D\cup \partial_M D$, the function $x \mapsto M_D(x, w)$  is excessive with respect to $X^D$;
\item{(c)} the function $(x,w) \mapsto M_D(x, w)$ is jointly continuous on $D\times ((D\setminus\{x_0\})\cup \partial_M D)$ in the Martin topology and
\item{(d)} $M_D(\cdot,w_1)\not=M_D(\cdot, w_2)$ if $w_1 \not= w_2$ and $w_1, w_2 \in \partial_M D$.
\end{description}
\end{description}

Recall that a positive harmonic function $f$ for  $X^{D}$ is minimal if, whenever $g$ is a positive harmonic function for $X^{D}$ with $g\le f$ on $D$, one must have $f=cg$ for some constant $c$. If $M_D(\cdot, z)$, $z\in \partial_M D$, is a minimal harmonic function, the point $z$ is called a minimal Martin boundary point. The set of all minimal Martin boundary points is denoted by $\partial_m D$. 
Then the following Martin representation is valid, see Theorem 4 in \cite{KW}:
For every non-negative function $h$ harmonic with respect to $X^D$,
there is a unique finite measure $\mu$ on $\partial_M D$ concentrated on $\partial_m D$ such that
\begin{equation}\label{e:martin-representation}
h(x)=\int_{\partial_M D}M_D(x,z)\, \mu(dz)=\int_{\partial_m D}M_D(x,z)\, \mu(dz)\, ,\quad x\in D\, .
\end{equation}

Recall that a point $w\in \partial_M D$ is a finite Martin boundary point if there exists a bounded sequence
$(y_n)_{n\ge 1}\subset D$ converging to $w$ in the Martin topology. The finite part of the Martin boundary will be denoted by $\partial_M^f D$.
Recall that a point $w$ on the Martin boundary $\partial_MD$ of $D$ is said to be associated
with $z_0\in \partial D$ if there is a sequence $(y_n)_{n\ge 1}\subset D$ converging to $w$
in the Martin topology and to $z_0$ in the topology of $\X$. The set of Martin
boundary points associated with $z_0$ is denoted by $\partial_M^{z_0} D$.

The proof of part (b) of the following result is a direct extension of that of
Lemma 4.18 in \cite{KSV14} and part (a) is even simpler.
So we omit the proof.

\begin{lemma}\label{l:irregular}
(a) Let $D$ be a bounded open set and
suppose that $u$ is a bounded non-negative harmonic function for $X^D$. If there exists a polar set $N\subset \partial D$ such that for any $z\in \partial D\setminus N$
\begin{equation}\label{e:not-polar}
\lim_{D\ni x\to z} u(x)=0\, ,
\end{equation}
then $u$ is identically equal to zero.

\noindent
(b) Let $D$ be an unbounded open set and
suppose that $u$ is a bounded non-negative harmonic function for $X^D$. If there exists a polar set $N\subset \partial D$ such that for any $z\in \partial D\setminus N$ \eqref{e:not-polar} holds true and additionally,
$$
\lim_{D\ni x\to \infty} u(x)=0\, ,
$$
then $u$ is identically equal to zero.
\end{lemma}

\subsection{Martin boundary at infinity}\label{ss:mb-infty}

In this subsection we assume that $R_0=\infty$, that there exists a point $z_0\in\X$
such that \textbf{C2}$(z_0, R)$ holds,
and that $\widehat{X}$ satisfies \textbf{F2}$(z_0, R)$ for some $R>0$.
We will fix $z_0$ and $R$ and use the notation $B_r=B(z_0, r)$.
Let $D$ be an unbounded open subset of $\X$ such that $\infty$ is accessible from
$D$ with respect to $X$.
We will deal with the Martin boundary of $D$ at infinity.
Recall that $x_0$ is a fixed point in $D$.

\begin{lemma}\label{l:limit-exists-i}
For every $x\in D$ the limit
\begin{align}\label{e:martin-kernel}
M_D(x,\infty):=\lim_{D\ni v\to \infty}\frac{G_D(x,v)}{G_D(x_0,v)}
\end{align}
exists and is finite.
\end{lemma}
\pf Fix $x\in D$, and let 
$r\ge  2 \min\{d(z_0,x),d(z_0,x_0),R\}$.
As before, let $D^r=D\cap \overline{B}_r^c$. The functions $G_D(x,\cdot)$ and $G_D(x_0,\cdot)$ are regular harmonic in $D^r$
with respect to $\widehat{X}$ and vanish in $\overline{B}_r^c \cap (\overline{D}^c\cup \wh{(D^r)}^{\mathrm{reg}})$,
hence
by Corollary \ref{c:oscillation-reduction-I} we deduce that the limit
$$
M_D(x,\infty):=\lim_{D\ni v\to \infty}\frac{G_D(x,v)}{G_D(x_0,v)}
$$
exists and is finite.
\qed

\noindent
\textbf{Proof of Theorem \ref{t:main-mb1}(a)}:
We first note that $\partial_M^{\infty}D$ is not empty.
Indeed, let $(y_n)_{n\ge 1}\subset D$ converge to ${\infty}$ in the topology of $\X$.
Since $D\cup \partial_M D$ is a compact metric space with the Martin metric $d_M$,
there exist a subsequence $(y_{n_k})_{k\ge 1}$ and $w\in D\cup \partial_M D$ such
that $\lim_{k\to \infty}d_M(y_{n_k},w)=0$. Clearly, $w\notin D$ (since relative
topologies on $D$ are equivalent). Thus we have found 
an unbounded sequence
$(y_{n_k})_{k\ge 1}\subset D$ which converges to $w\in \partial_M D$ in
the Martin topology and to ${\infty}$ in the topology of $\X$.

Let $w\in \partial_M^{\infty}D$ and let $M_D(\cdot, w)$ be the corresponding Martin kernel.
If  $(y_n)_{n\ge 1}$ is a sequence in $D$ converging to $w$
in the Martin topology and to ${\infty}$ in the topology of $\X$,
then, by (M3)(a), $M_D(x,y_n)$ converges to $M_D(x,w)$. On the other hand,
since $y_n$ converges ${\infty}$ in the topology of $\X$,
by Lemma \ref{l:limit-exists-i},
$
\lim_{n\to \infty}M_D(x,y_n)=M_D(x, {\infty}).
$
Hence, for each $w\in \partial_M^{\infty} D$
it holds that $M_D(\cdot, w)=M_D(\cdot, {\infty})$.
Since, by (M3)(d),  for two different Martin boundary points $w^{(1)}$
and $w^{(2)}$ it always holds that $M_D(\cdot, w^{(1)})\neq M_D(\cdot, w^{(2)})$,
we conclude that $\partial_M^{\infty}D$ consists of exactly one point.
\qed

\noindent
\textbf{Proof of Theorem \ref{t:main-mb1}(b)}:
We claim that for every $r>4\max(d(z_0,x_0), R)$ and $U:=D\cap B_r$ it holds that
\begin{equation}\label{e:M-mean-value-I}
M_D(x,\infty)=\E_x\left[M_D(X_{\tau_U},\infty)\right],\qquad x\in U.
\end{equation}
For any $z\in D^{2r}$, since $G_D(\cdot, z)$ is regular harmonic in $U$, we have
$$
\frac{G_D(x,z)}{G_D(x_0,z)}=\E_x\left[\frac{G_D(X_{\tau_U},z)}{G_D(x_0,z)}\right], \qquad x\in U.
$$
Hence, in view of Lemma \ref{l:limit-exists-i}, in order to prove \eqref{e:M-mean-value-I} it suffices to show that,
for any fixed $x\in U$, there exists $s>16r$ such that the family
$$
\left\{\frac{G_D(X_{\tau_U},z)}{G_D(x_0,z)}: z \in D^{4s}\right\}
$$
is uniformly integrable with respect to the distribution of $X_{\tau_U}$ under $\P_x$.

In the remainder of this proof, we fix an $x\in U$.
Let $s>8r$. Then for any Borel set $E  \subset D^r$,
\begin{eqnarray*}
\lefteqn{\E_x\left[\frac{G_D(X_{\tau_U},z)}{G_D(x_0,z)}, X_{\tau_U}\in E\right]}\\
 &\le & \E_x\left[\frac{G_D(X_{\tau_U},z)}{G_D(x_0,z)}, X_{\tau_U}\in (D^r\setminus D^{s/3})\cap E\right]
+\E_x\left[\frac{G_D(X_{\tau_U},z)}{G_D(x_0,z)}, X_{\tau_U}\in D^{s/3}\right]\\ &=&:I+II\, .
\end{eqnarray*}

We first show that $II$ is small for large $s$.
Let $w\in U$ and $d(z_0,y)>4r/3$.
By \textbf{C2}$(z_0, R)$
we have that $j(w,y)\le c_1 j(z_0,y)$ with $c_1=c_1(z_0,4/3)$.
It follows that
\begin{eqnarray*}
P_U(x,y)&=&\int_U G_U(x,w)j(w,y)m(dw) \le c_1 (\E_x \tau_U) j(z_0,y)\\
&\le & c_1 (\E_x \tau_{B_r})  j(z_0,y)=c_2 j(z_0,y),
\end{eqnarray*}
where $c_2=c_2(z_0, x, r)$. Hence,
\begin{eqnarray}\label{e:mkh-1-I}
\E_x[G_D(X_{\tau_U},z), X_{\tau_U}\in D^{s/3}]&=&\int_{D^{s/3}}G_D(y,z) P_U(x,y) m(dy)
 \nonumber \\
&\le& c_2  \int_{D^{s/3}}G_D(y,z) j(z_0,y) m(dy) .
\end{eqnarray}
Next, for $z\in D^{4s}$,
\begin{eqnarray*}
G_D(x_0,z)&\ge &\int_{D^{s/3}}G_D(y,z)P_{D\setminus \overline{D}^{s/4}}(x_0,y)m(dy)\\
&=&\int_{D^{s/3}}\int_{D\setminus \overline{D}^{s/4}}G_D(y,z)
G_{D\setminus\overline{D}^{s/4}}(x_0,u)j(u,y)\, m(du)\, m(dy) .
\end{eqnarray*}
Let $y\in D^{s/3}$ and $u\in D\setminus \overline{D}^{s/4}$.
By \textbf{C2}$(z_0, R)$
we have that $j(z_0,y)\le c_4 j(u,y)$ with $c_3=c_3(z_0)$.
Continuing the above display, we get that
\begin{equation}\label{e:mkh-2-I}
G_D(x_0,z)
\ge c_3^{-1}\left(\int_{D^{s/3}}G_D(y,z) j(z_0, y) m(dy) \right)
\left(\int_{D\setminus \overline{D}^{s/4}}G_{D\setminus \overline{D}^{s/4}}(x_0,u) m(du) \right) .
\end{equation}
By combining \eqref{e:mkh-1-I} and \eqref{e:mkh-2-I}
we get that for all $z\in D^{4s}$,
$$
II=\int_{D^{s/3}}\frac{G_D(y,z)}{G_D(x_0,z)}P_U(x,y)m(dy) \le c_2 c_3
\left(\int_{D\setminus\overline{D}^{s/4}}G_{D\setminus \overline{D}^{s/4}}(x_0,u)m(du) \right)^{-1}.
$$
Since
$$
\lim_{s\to \infty}\int_{D\setminus \overline{D}^{s/4}}G_{D\setminus \overline{D}^{s/4}}(x_0,u)m(du)
=\int_D G_D(x_0,u) m(du) = \E_{x_0} \tau_D =\infty ,
$$
we see that for any $\epsilon >0$ we can find $s>16r$ such that
$$
\left(\int_{D\setminus \overline{D}^{s/4}}G_{D\setminus \overline{D}^{s/4}}(x_0,u) m(du)\right)^{-1}
< \frac{\epsilon}{2 c_2 c_3}.
$$
Thus $II<\epsilon/2$ for all $z  \in D^{4s}$.

We now fix an $s>16r$ as above and estimate $I$ for all $z  \in D^{4s}$.
If $y\in D^r\setminus D^{s/3}$, then  both $G_D(y,\cdot)$ and
$G_D(x_0,\cdot)$ are regular harmonic with respect to $\wh{X}$ in $D^{s/2}$ and vanish on
$B_{s/2}\cap (\overline{D}^c\cup \wh{D}^{\mathrm{reg}})$.
Choose $z_1\in D^s$.
By the boundary Harnack principle
\eqref{e:bhp-I}, we have that
$$
\frac{G_D(y,z)}{G_D(x_0,z)}\le c_4 \frac{G_D(y,z_1)}{G_D(x_0,z_1)},\qquad z\in D^{4s}.
$$
Since $z_1\in D^s$ it follows from \eqref{e:nGassup1-infty} that
$c_5:=\sup_{y\in D^r\setminus D^{s/3}}G_D(y,z_1)<\infty .$
Hence,
\begin{align}\label{e:estimate-I} 
I & \le \ c_4\E_x\left[ \frac{G_D(X_{\tau_U},z_1)}{G_D(x_0,z_1)}, X_{\tau_U}\in (D^r\setminus D^{s/3})\cap E \right]  \nonumber \\
&\le \ \frac{c_4 c_5}{G_D(x_0,z_1)}\P_x\big( X_{\tau_U}\in (D^r\setminus D^{s/3})\cap E\big)
\le c_6\P_x(X_{\tau_U}\in E)\, ,
\end{align}
where $c_6=c_4 c_5/G_D(x_0,z_1)$.
Thus, given $\epsilon >0$, for any set $E\subset D^r$ with $\int_E P_U(x,y)m(dy) <\epsilon/(2c_7)$,
we have $I<\epsilon/2$ for all $z  \in D^{4s}.$

Therefore we have proved the claimed uniform integrability
for the $s$ chosen above, and consequently \eqref{e:M-mean-value-I}.

Now let $U_1\subset D$ be any bounded open set such that $\overline{U}_1\subset D$.
Then there is $r>4R$ such that $U_1\subset D\cap B_r=:U$. Then by
\eqref{e:M-mean-value-I} and the strong Markov property we have that
\begin{equation}\label{e:M-mean-value-I-2}
M_D(x,\infty)=\E_x\left[M_D(X_{\tau_{U_1}},\infty)\right], \qquad x\in U_1,
\end{equation}
which finishes the proof. \qed

Because of Theorem \ref{t:main-mb1}(a), we will also use $\infty$ to denote the Martin boundary point
$\partial_M^{\infty}D$ associated with $\infty$.
Note that it follows from the proof of Theorem \ref{t:main-mb1}(a) that
if $(y_n)_{n\ge 1}$ converges to $\infty$ in the topology of $\X$,
then it also converges to $\infty$ in the Martin topology.

For any $\epsilon >0$, define
\begin{equation}\label{e:definition-U_K-infty}
K^{\infty}_{\epsilon}:=\left\{w\in \partial_M^f D: 
d_M(w, \infty) \ge \epsilon\right\}.
\end{equation}
By the definition of the finite part of the Martin boundary,
for each $w\in K^{\infty}_{\epsilon}$ there exists a bounded sequence $(y_n^w)_{n\ge 1}\subset D$ such that
$\lim_{n\to \infty} d_M(y_n^w, w)=0$. Without loss of generality we may assume that
$d_M(y_n^w, w)<\frac{\epsilon}{2}$ for all $n\ge 1$.

\begin{lemma}\label{l:boundedness-U_K}
There exists $c=c(\epsilon)>0$ such that $d(y_n^w,z_0)\le c$ for all $w\in K_{\epsilon}^{\infty}$ and all $n\ge 1$.
\end{lemma}
\pf
We first claim that for any sequence $(y_n)_{n\ge 1}$ in $D$ satisfying $d(y_n, z_0)\to \infty$, we have
$\lim_{n\to \infty}d_M(y_n, {\infty})=0$, i.e., $(y_n)_{n\ge 1}$ converges to ${\infty}$
in the Martin topology.
Indeed, since $D\cup \partial_M D$ is a compact metric space, $(y_n)$ has  a convergent
subsequence $(y_{n_k})$. Let $w=\lim_{k\to \infty}y_{n_k}$ (in the Martin topology).
Then $\lim_{k\to \infty}M_D(\cdot, y_{n_k})=M_D(\cdot, w)$. On the other hand, from
Lemma \ref{l:limit-exists-i}
and Theorem \ref{t:main-mb1}(a) we see that $\lim_{k\to \infty}M_D(\cdot, y_{n_k})=M_D(\cdot, \infty)$.
Therefore, $M_D(\cdot,w)= M_D(\cdot, {\infty})$, which  implies
that $w={\infty}$ by (M3)(d). Since this argument also holds for
any subsequence of $(y_n)_{n\ge 1}$, we conclude that $y_n\to {\infty}$ in the Martin topology.

Now suppose the lemma is not true. Then $\{y_n^w:\, w\in K_{\epsilon}^{\infty}, n\in \N\}$
contains a sequence $(y_{n_k}^{w_k})_{k\ge 1}$ such that $\lim_{k\to \infty}d(y_{n_k}^{w_k}, z_0)= \infty$.
By the paragraph above, we have
that $\lim_{k\to \infty}d_M(y_{n_k}^{w_k},{\infty})= 0$. On the other hand,
$
d_M(y_{n_k}^{w_k},{\infty})\ge d_M(w_k,{\infty})-d_M(y_{n_k}^{w_k}, w_k)
\ge {\epsilon}/{2}\, .
$
This contradiction proves the claim. \qed

\noindent
\textbf{Proof of Theorem \ref{t:main-mb1}(c)}:
Let $h$ be a positive harmonic function for  $X^{D}$
such that $h\le M_D(\cdot, {\infty})$. By
the Martin representation \eqref{e:martin-representation},
there is a finite measure $\mu$ on 
$\partial_M D$ (concentrated on $\partial_m D$) such that
$$
h(x)=\int_{\partial_M D}M_D(x,w)\, \mu(dw)=\int_{\partial_M D  \setminus \{\infty\}}M_D(x,w)\, \mu(dw)+ M_D(x,
\infty)\mu(\{\infty\})\, .
$$
In particular, $h(x_0)=\mu(\partial_M D)\le  M_D(x_0,  \infty)=1$ 
(because of the normalization at $x_0$). Hence, $\mu$ is a sub-probability measure.

For $\epsilon >0$, let $K^{\infty}_{\epsilon}$
be the closed subset of $\partial^f_M D$ defined in \eqref{e:definition-U_K-infty}.
Define
\begin{equation}\label{d:definition-u-I}
    u(x):=\int_{ K^{\infty}_{\epsilon} }M_D(x,w)\, \mu(dw).
 \end{equation}
Then $u$ is a positive harmonic function with respect to  $X^{D}$
and bounded above by
\begin{align}
u(x) \le  h(x)-\mu(\{\infty\})M_D(x,
{\infty})\le \big(1-\mu(\{\infty\})\big)M_D(x,
{\infty})\, . \label{e:newm11}
\end{align}
We claim that $\lim_{x\to \infty} u(x)=0$. Let $p=c(\epsilon) \vee R$, where $c(\epsilon)$ is the constant from Lemma \ref{l:boundedness-U_K}.
Hence, for $w\in K^{\infty}_{\epsilon}$ 
and $(y_n^w)_{n\ge 1}$ a sequence such that
$\lim_{n\to \infty}d_M(y_n^w,w)=0$, it holds that $d(y_n^w, {z_0})\le p$.
Fix a point $x_1\in D^{8p}$ and choose an arbitrary point $y_0\in D_p$. Then for any $x\in D^{8p}$ and any $y\in D_p$ we have that
\begin{eqnarray*}
\frac{G_D(x,y)}{G_D(x_0,y)}=\frac{G_D(x,y)}{G_D(x_1,y)}\, \frac{G_D(x_1,y)}{G_D(x_0,y)}
\le c_1 \frac{G_D(x,y_0)}{G_D(x_1,y_0)}\, \frac{G_D(x_1,y)}{G_D(x_0,y)},
\end{eqnarray*}
where the inequality follows from 
the dual version of \eqref{e:bhp-I}
since $X$ satisfies \textbf{F2}$(z_0, R)$.
Therefore for each $w\in K_{\epsilon}^{\infty}$  we have
\begin{eqnarray*}
M_D(x,w)&=&\lim_{n\to \infty}\frac{G_D(x,y_n^w)}{G_D(x_0,y_n^w)}\le c_1 \frac{G_D(x,y_0)}{G_D(x_1,y_0)}\, \lim_{n\to \infty} \frac{G_D(x_1,y_n^w)}{G_D(x_0,y_n^w)}\\
&=& c_1\frac{G_D(x,y_0)}{G_D(x_1,y_0)}\ M_D(x_1, w) \le c_1\frac{G_D(x,y_0)}{G_D(x_1,y_0)}\ \sup_{w\in K_{\epsilon}} M_D(x_1, w)\\
&=& c_2 G_D(x,y_0)
\end{eqnarray*}
by the continuity of the Martin kernel (M3)(c).
This inequality together with the definition of $u$ shows that $u(x)\le c_2 G_D(x, y_0)$. Now using \eqref{e:nGassup2i}, we can conclude that $\lim_{x\to \infty, x\in D}u(x)=0$ uniformly for $w\in K_{\epsilon}^{\infty}$.

Choose $r>16p$. For any $x\in D_{r/2}$ and 
$y, y_1\in D^{8r}$,
by \eqref{e:bhp-I} applied to $G_D(x, \cdot)$ and $G_D(x_0,\cdot)$, we have
$$
\frac{G_D(x, y)}{G_D(x_0, y)}\le c_3\frac{G_D(x, y_1)}{G_D(x_0, y_1)}.
$$
Letting $D\ni y\to \infty$, by Lemma \ref{l:limit-exists-i} we get
\begin{equation}\label{e:martinkernelub-I}
M_D(x, \infty)\le c_3\frac{G_D(x, y_1)}{G_D(x_0, y_1)}=c_4 G_D(x, y_1)\, ,\quad x\in D_{r/2}.
\end{equation}
Recall that by \eqref{e:zreg} $\lim_{D\ni x\to z}G_D(x,y)=0$ for every $z\in \partial D$ which is regular for $D^c$
with respect to $X$.
Since $r>16p$ can be arbitrarily large, we see from \eqref{e:martinkernelub-I} and
\eqref{e:newm11} that  $\lim_{D\ni x, x\to z}u(x)=0$ for every $z\in \partial D$
which is regular for $D^c$ with respect to $X$.

Fix $r>16p$ and $y_1\in D^{8r}$.
 It follows from \eqref{e:nGassup1-infty}  that for all
$x\in D_{r/2}$,
\begin{equation}\label{e:martinkernelub3-i}
G_D(x, y_1)\le c_5\, .
\end{equation}
From  \eqref{e:newm11}--\eqref{e:martinkernelub3-i}
we conclude that $u$ is
bounded in $x\in D_{r/2}$.
Similarly by \eqref{e:nGassup1-infty}, for every $x\in D^{8p}$
we have that  $G_D(x,y_0)\le c_6$
(recall $y_0\in D_p$). Since 
$M_D(x,w)\le c_2 G_D(x,y_0)$
 for each $x\in D^{8p}$ and each $w\in K_{\epsilon}^{\infty}$,
by using \eqref{d:definition-u-I} we see that $u$ is bounded on 
$D^{8p}$. 
Thus $u$ is bounded on $D$. 

Now it follows from Lemma \ref{l:irregular}(b) that $u\equiv 0$ in $D$.
This means that $\nu=\mu_{|
K^{\infty}_{\epsilon}}=0$. Since
$\epsilon >0$ was arbitrary and 
$\partial_M D  \setminus \{\infty\}=\cup_{\epsilon >0} K^{\infty}_{\epsilon}$,
we see that $\mu_{|\partial_M D \setminus \{\infty\}}=0$. 
Hence $h=\mu(\{\infty\})M(\cdot,
\infty)$ showing that $M(\cdot, \infty)$ is minimal. Therefore we have proved the theorem.
\qed

\subsection{Finite part  of  Martin boundary}\label{ss:mb-finite}

In this subsection, we deal with the oscillation reduction at a boundary point $z_0$ of an open set $D$.
We will fix $D$ and $z_0$ in this subsection, and use the notation $B_r=B(z_0, r)$.
In this subsection, we will always assume that there exists
$R\le R_0$ such that \textbf{C1}$(z_0, R)$ holds,
and that $\widehat{X}$ satisfies \textbf{F1}$(z_0, R)$.
We also assume that $z_0$
is accessible from $D$ with respect to $X$.

\begin{lemma}\label{l:limit-exists}
For every $x\in D$, the limit
$$
M_D(x,z_0):=\lim_{D\ni v\to z_0}\frac{G_D(x,v)}{G_D(x_0,v)}$$
exists and is finite.
\end{lemma}

\pf
Fix $x\in D$, and let $r\le \frac{1}{2}\min\{d(z_0,x),d(z_0,x_0),R\}$.
As before, let $D_r=D\cap B_r$. The functions $G_D(x,\cdot)$ and $G_D(x_0,\cdot)$ are regular harmonic in $D_r$ with respect to $\widehat{X}$ and vanish in
$B_r\cap ( \overline{D}^c\cup \wh{(D_r)}^{\mathrm{reg}})$,
hence
by Corollary \ref{c:reductionac} we deduce that the limit
$$
M_D(x,z_0):=\lim_{D\ni v\to z_0}\frac{G_D(x,v)}{G_D(x_0,v)}
$$
exists and is finite.
\qed

\noindent
\textbf{Proof of Theorem \ref{t:main-mb0}(a)}:
We first note that $\partial_M^{z_0}D$ is not empty.
Indeed, let $(y_n)_{n\ge 1}\subset D$ converge to ${z_0}$ in the topology of $\X$.
Since $D\cup \partial_M D$ is a compact metric space with the Martin metric $d_M$,
there exist a subsequence $(y_{n_k})_{k\ge 1}$ and $w\in D\cup \partial_M D$ such
that $\lim_{k\to \infty}d_M(y_{n_k},w)=0$. Clearly, $w\notin D$ (since relative
topologies on $D$ are equivalent). Thus we have found a sequence
$(y_{n_k})_{k\ge 1}\subset D$ which converges to $w\in \partial_M D$ in
the Martin topology and to ${z_0}$ in the topology of $\X$.

Let $w\in \partial_M^{z_0}D$ and let $M_D(\cdot, w)$ be the corresponding Martin kernel.
If  $(y_n)_{n\ge 1}$ is a sequence in $D$ converging to $w$
in the Martin topology and to ${z_0}$ in the topology of $\X$,
then, by (M3)(a), $M_D(x,y_n)$ converge to $M_D(x,w)$. On the other hand, $d(y_n, {z_0})\to 0$, thus by Lemma \ref{l:limit-exists},
$
\lim_{n\to \infty}M_D(x,y_n)=M_D(x,{z_0}).
$
Hence, for each $w\in \partial_M^{z_0} D$ it holds that $M_D(\cdot, w)=M_D(\cdot, {z_0})$.
Since, by (M3)(d),  for two different Martin boundary points $w^{(1)}$
and $w^{(2)}$ it always holds that $M_D(\cdot, w^{(1)})\neq M_D(\cdot, w^{(2)})$,
we conclude that $\partial_M^{z_0}D$ consists of exactly one point.
\qed

\noindent
\textbf{Proof of Theorem \ref{t:main-mb0}(b)}:
We claim that for every $r\le \frac14 \min\{d(z_0,x_0),R\}$ and $U:=D\setminus \overline{B}_r$ it holds that
\begin{equation}\label{e:M-mean-value}
M_D(x,z_0))=\E_x\left[M_D(X_{\tau_U},z_0)\right]\, ,\qquad x\in U.
\end{equation}
For any $z\in D_{r/2}$, since $G_D(\cdot, z)$ is regular harmonic in $U$, we have
$$
\frac{G_D(x,z)}{G_D(x_0,z)}=\E_x\left[\frac{G_D(X_{\tau_U},z)}{G_D(x_0,z)}\right], \quad x\in U.
$$
Hence, in view of Lemma \ref{l:limit-exists}, in order to prove \eqref{e:M-mean-value} it suffices to show that,
for any fixed $x\in U$, there exists $s<r/(16)$ such that the family
$$
\left\{\frac{G_D(X_{\tau_U},z)}{G_D(x_0,z)}: z \in D_{s/4}\right\}
$$
is uniformly integrable with respect to the distribution of $X_{\tau_U}$ under $\P_x$.

In the remainder of this proof we fix an $x\in U$.
Let $0<s<r/8$. Then for any Borel set $E \subset D_r$,
\begin{eqnarray*}
\E_x\left[\frac{G_D(X_{\tau_U},z)}{G_D(x_0,z)}, X_{\tau_U}\in E\right] &\le & \E_x\left[\frac{G_D(X_{\tau_U},z)}{G_D(x_0,z)}, X_{\tau_U}\in (D_r\setminus D_{3s})\cap E\right]\\
& & +\E_x\left[\frac{G_D(X_{\tau_U},z)}{G_D(x_0,z)}, X_{\tau_U}\in D_{3s}\right]\\
& =&:I+II\, .
\end{eqnarray*}

We first show that $II$ is small for small $s$.
We claim that $P_U(x,\cdot)$ is bounded on $D_{3r/4}$.
Indeed, let $y\in D_{3r/4}$. If $w\in U$, then by
\textbf{C1}$(z_0, R)$,
we have that $j(w,y)\le  c_1 j(w, z_0)$ where $c_1=c_1(z_0, 4/3)$. Hence,
\begin{align*}
P_U(x,y)&=\int_U G_U(x,w)j(w,y)m(dw)\le c_1 \int_U G_U(x,w) j
(w,z_0)m(dw)\\
&=c_1 P_U(x,z_0)=: c_2.
\end{align*}
This implies that
\begin{eqnarray}\label{e:mkh-1}
\E_x[G_D(X_{\tau_U},z)), X_{\tau_U}\in D_{3s}] &=&\int_{D_{3s}}G_D(y,z)P_U(x,y)m(dy) \nonumber \\
&\le& c_2 \int_{D_{3s}}G_D(y,z)m(dy) .
\end{eqnarray}
Next, for $z\in D_{s/4}$,
\begin{eqnarray*}
G_D(x_0,z)&\ge &\int_{D_{3s}}G_D(y,z)P_{D\setminus \overline{D}_{4s}}(x_0,y)m(dy)\\
&=&\int_{D_{3s}}\int_{D\setminus \overline{D}_{4s}}G_D(y,z)G_{D\setminus
\overline{D}_{4s}}(x_0,u)j(u,y)\, m(du)\, m(dy) .
\end{eqnarray*}
Let $y\in D_{3s}$ and $u\in D\setminus \overline{D}_{4s}$.
By \textbf{C1}$(z_0, R)$,
we have that $j(u,z_0)\le c_3 j(u,y)$ with $c_3=c_3(z_0)$.
Continuing the above display, we get that
\begin{equation}\label{e:mkh-2}
G_D(x_0,z)\ge c_3^{-1}\left(\int_{D_{3s}}G_D(y,z) m(dy) \right)
\left(\int_{D\setminus\overline{D}_{4s}}G_{D\setminus
\overline{D}_{4s}}(x_0,u)j(u,z_0)m(du) \right) .
\end{equation}
Combining \eqref{e:mkh-1} and \eqref{e:mkh-2} we arrive at
\begin{eqnarray*}
II&=&\int_{D_{3s}}\frac{G_D(y,z)}{G_D(x_0,z)}P_U(x,y)m(dy)\\
& \le& c_2 c_3
\left(\int_{D\setminus \overline{D}_{4s}}G_{D\setminus
\overline{D}_{4s}}(x_0,u)j(u,z_0)m(du) \right)^{-1}.
\end{eqnarray*}
Since $z_0$ is accessible from $D$,
\begin{align*}
&\lim_{s\to 0}\int_{D\setminus \overline{D}_{4s}}G_{D\setminus \overline{D}_{4s}}(x_0,u)j(u,z_0)m(du)\\
&=\int_D G_D(x_0,u)j(u,z_0) m(du) =P_D(x_0,z_0)=\infty.
\end{align*}
Therefore, for any $\epsilon >0$ one can find $s>0$ such that
$$
\left(\int_{D\setminus \overline{D}_{4s}}G_{D\setminus \overline{D}_{4s}}(x_0,u)j(u,z_0)m(du)
\right)^{-1}< \frac{\epsilon}{2 c_2 c_3}.
$$
Thus $II<\epsilon/2$ for all $z \in D_{s/4}$.

We fix $s<r/16$ as above and estimate $I$ for all $z  \in D_{s/4}$.
Choose $z_1\in D_s$. If $y\in D_r\setminus D_{3s}$,
then  both $G_D(y,\cdot)$ and $G_D(x_0,\cdot)$ are regular harmonic in $D_{2s}$ with
respect to $\widehat{X}$ and vanish on
$B_{2s}\cap ( \overline{D}^c\cup \wh{(D_{2s})}^{\mathrm{reg}})$.
Hence, by the boundary Harnack principle \eqref{e:bhp}, we have that
$$
\frac{G_D(y,z)}{G_D(x_0,z)}\le c_4 \frac{G_D(y,z_1)}{G_D(x_0,z_1)},\qquad z\in D_{s/4}.
$$
Since $z_1\in D_s$ it follows from \eqref{e:nGassup1} that
$c_5:=\sup_{y\in D_r\setminus D_{3s}}G_D(y,z_1)<\infty .$
Hence,
\begin{eqnarray}\label{e:estimate}
I & \le & c_4\E_x\left[ \frac{G_D(X_{\tau_U},z_1)}{G_D(x_0,z_1)}, X_{\tau_U}\in (D_r\setminus D_{3s})\cap E \right]  \nonumber \\
&\le &\frac{c_4 c_5}{G_D(x_0,z_1)}\P_x\big( X_{\tau_U}\in (D_r\setminus D_{3s})\cap E\big)
\le c_6\P_x(X_{\tau_U}\in E), 
\end{eqnarray}
where $c_6=c_4c_5/G_D(x_0,z_1)$.
Thus, given $\epsilon >0$, for any set $E\subset D_r$ with $\P_x(X_{\tau_U}\in E) < \epsilon/(2c_6)$
we have that $I<\epsilon/2$ for all $z  \in D_{s/4}$.

Therefore we have proved the claimed uniform integrability for the $s$ chosen above, and consequently \eqref{e:M-mean-value}.

Now let $U_1\subset D$ be any open set such that $z_0$ is not in $\overline{U}_1$.
Then there is $r>0$ such that $U_1\subset D\setminus \overline{B}_r=:U$.
Then by \eqref{e:M-mean-value} and the strong Markov property  we get that
\begin{equation}\label{e:M-mean-value-2}
M_D(x,z_0)=\E_x\left[M_D(X_{\tau_{U_1}},z_0)\right], \qquad x\in U_1,
\end{equation}
which finishes the proof. \qed

Because of Theorem \ref{t:main-mb0}(a), we will also use ${z_0}$ to denote the Martin boundary point
$\partial_M^{z_0}D$ associated with $z_0\in \partial D$.
Note that it follows from the proof of Theorem \ref{t:main-mb0}(a) that
if $(y_n)_{n\ge 1}$ converges to ${z_0}$ in the topology of $\X$,
then it also converges to ${z_0}$ in the Martin topology.

For any $\epsilon >0$, define
\begin{equation}\label{e:definition-U_K}
K^{z_0}_{\epsilon}:=\left\{w\in \partial_M^f D: 
d_M(w,z_0) \ge \epsilon\right\}.
\end{equation}
By the definition of the finite part of the Martin boundary,
for each $w\in K^{z_0}_{\epsilon}$ there exists a bounded sequence $(y_n^w)_{n\ge 1}\subset D$ such that
$\lim_{n\to \infty} d_M(y_n^w, w)=0$. Without loss of generality we may assume that
$d_M(y_n^w, w)<\frac{\epsilon}{2}$ for all $n\ge 1$.

\begin{lemma}\label{l:boundedness-U-K}
There exists $c=c(\epsilon)>0$ such that $d(y_n^w, {z_0})\ge c$ for all $w\in K_{\epsilon}^{z_0}$ and all $n\ge 1$.
\end{lemma}
\pf
Suppose the lemma is not true. Then $\{y_n^w:\, w\in K_{\epsilon}^{z_0}, n\in \N\}$
contains a subsequence $(y_{n_k}^{w_k})_{k\ge 1}$ such that $\lim_{k\to \infty}d(y_{n_k}^{w_k}, {z_0})= 0$.
We also have
$\lim_{k\to \infty}d_M(y_{n_k}^{w_k}, {z_0})= 0$. On the other hand,
$$
d_M(y_{n_k}^{w_k}, {z_0})\ge d_M(w_k, {z_0})-d_M(y_{n_k}^{w_k}, w_k)\ge {\epsilon}/{2}.
$$
This contradiction proves the claim. \qed

\medskip
\noindent
\textbf{Proof of Theorem \ref{t:main-mb0}(c):}
Let $h$ be a positive harmonic function for  $X^{D}$ such that $h\le M_D(\cdot, {z_0})$. By the Martin representation  \eqref{e:martin-representation},
there is a finite measure $\mu$ on $\partial_M D$ (concentrated on $\partial_m D$)  such that
$$
    h(x)=\int_{\partial_M D}M_D(x,w)\, \mu(dw)=\int_{\partial_M D\setminus\{{z_0}\}}M_D(x,w)\, \mu(dw)+ M_D(x,{z_0})\mu(\{{z_0}\})\, .
$$
In particular, $\mu(\partial_M D)=h(x_0)\le M_D(x_0, {z_0})=1$ (because of the normalization at $x_0$). Hence, $\mu$ is a sub-probability measure.

For $\epsilon >0$, let $K^{z_0}_{\epsilon}$ be the closed subset of $\partial_M D$ defined in \eqref{e:definition-U_K}.
Define
\begin{equation}\label{d:definition-u}
    u(x):=\int_{ K^{z_0}_{\epsilon} }M_D(x,w)\, \mu(dw).
\end{equation}
Then $u$ is a positive harmonic function with respect to  $X^{D}$
satisfying
\begin{align}\label{e:newm1}
u(x)\le h(x)-\mu(\{{z_0}\})M_D(x,
{z_0})\le \big(1-\mu(\{{z_0}\})\big)M_D(x, {z_0})\, .
\end{align}

Let $p=c(\epsilon) \wedge R$, where $c(\epsilon)$ is the constant from Lemma \ref{l:boundedness-U-K}.
Hence, for $w\in K^{z_0}_{\epsilon}$ and $(y_n^w)_{n\ge 1}$ a sequence such that
$\lim_{n\to \infty}d_M(y_n^w,w)=0$, it holds that 
$d(y_n^w, {z_0})\ge p$.
Fix $x_1\in D_{p/8}$ and choose arbitrary $y_0\in D^p$.
For any $x\in D_{p/8}$ and any $y\in D^p$ we have that
$$
\frac{G_D(x,y)}{G_D(x_0,y)}=\frac{G_D(x,y)}{G_D(x_1,y)}\, \frac{G_D(x_1,y)}{G_D(x_0,y)}\le c_1 \frac{G_D(x,y_0)}{G_D(x_1,y_0)}\, \frac{G_D(x_1,y)}{G_D(x_0,y)}\, .
$$
Here the inequality follows from the dual version of \eqref{e:bhp} applied to functions $G_D(\cdot, y)$
and $G_D(\cdot, y_0)$ which are regular harmonic in $D_p$ with respect to $X$ and vanish in
$B({z_0},c)\cap ( \overline{D}^c\cup D^{\mathrm{reg}})$.
Now fix $w\in K^{z_0}_{\epsilon}$ and apply the
above inequality to $y_n^w$ to get
\begin{align*}
&M_D(x,w)=\lim_{n\to \infty}\frac{G_D(x,y_n^w)}{G_D(x_0, y_n^w)}\le c_1 \frac{G_D(x,y_0)}{G_D(x_1,y_0)}\, \lim_{n\to
\infty}\frac{G_D(x_1, y_n^w)}{G_D(x_0, y_n^w)}\\
&=c_1  \frac{G_D(x,y_0)}{G_D(x_1,y_0)}\, M_D(x_1,w)\le c_1  \frac{G_D(x,y_0)}{G_D(x_1,y_0)} \sup_{w\in K_{\epsilon}^{z_0}}M_D(x_1,w)\\
&\le c_2  \frac{G_D(x,y_0)}{G_D(x_1,y_0)} = c_3  G_D(x,y_0)\, .
\end{align*}
In the  last inequality we used property (M3)(c) of the Martin kernel. Thus,
\begin{equation}\label{e:bound-on-K-epsilon}
M_D(x,w)\le c_3 G_D(x,y_0)\, , \qquad x\in D_{p/8}, w\in K^{z_0}_{\epsilon}\, .
\end{equation}

Choose $r<p/16$. For any $x\in D^{2r}$ and $y, y_1\in D_{r/8}$,
by \eqref{e:bhp} applied to $G_D(x, \cdot)$ and $G_D(x_0,\cdot)$, we have
$$
\frac{G_D(x, y)}{G_D(x_0, y)}\le c_4\frac{G_D(x, y_1)}{G_D(x_0, y_1)}.
$$
Letting $D\ni y\to {z_0}$, we get
\begin{equation}\label{e:martinkernelub}
M_D(x, {z_0})\le c_4\frac{G_D(x, y_1)}{G_D(x_0, y_1)}=c_5 G_D(x, y_1)\, ,\quad x\in D^{2r}.
\end{equation}
Recall that by \eqref{e:zreg} $\lim_{D\ni x\to z}G_D(x,y)=0$ for every $z\in \partial D$ which is regular for $D^c$
with respect to $X$.
Since $r<p/16$ can be arbitrarily small, we see from \eqref{e:martinkernelub} and
\eqref{e:newm1} that  $\lim_{D\ni x, x\to z}u(x)=0$ for every $z\in \partial D$, $z\neq {z_0}$,
which is regular for $D^c$ with respect to $X$.

Assume $D$ is bounded. Fix $r<p/16$. It follows from \eqref{e:nGassup1} that for all
$x\in D^{2r}$,
\begin{equation}\label{e:martinkernelub3}
G_D(x, y_1)\le c_7\, .
\end{equation}
From \eqref{e:martinkernelub} and \eqref{e:newm1}  we conclude that $u$ is
bounded in $x\in D^{2r}$.
Similarly by \eqref{e:nGassup1}, for every $x\in D_{p/8}$
we have that  $G_D(x,y_0)\le c_9$
(recall $y_0\in D^p$). Hence by \eqref{e:bound-on-K-epsilon}
and  \eqref{d:definition-u} we see that $u$ is bounded on $D_{p/8}$.
Thus $u$ is bounded on $D$. Now it follows from Lemma \ref{l:irregular}(a) that $u\equiv 0$ in $D$.

If $D$ is unbounded, we argue as follows.
It follows from \eqref{e:bound-on-K-epsilon} and the assumption \eqref{e:nGassup2} that
$\lim_{D\ni x\to \infty} M_D(x,{z_0})=0$. Hence by \eqref{e:newm1}
$\lim_{D\ni x\to \infty}u(x)=0$. Thus, there exists
$\overline{r}\ge 2$ such that $u(x)\le 1$
for all $x\in D^{\overline{r}})$. Fix $r<p/16\wedge 1$ and let $x\in D\cap (B({z_0},\overline{r})\setminus B({z_0},2r))$. By \eqref{e:martinkernelub} and \eqref{e:nGassup1},
$$
M_D(x,{z_0})\le c_5 G_D(x, y_1)\le c_{11}\, .
$$
It follows that $u$ is bounded in $D\cap (B({z_0},\overline{r})\setminus B({z_0},2r))$.
The proof that $u$ is bounded on $D\cap B({z_0},p/16)$ is the same as in the
case of a bounded $D$. Hence, $u$ is bounded, and again we conclude from
Lemma \ref{l:irregular} (b) that $u\equiv 0$ in $D$.

We see from \eqref{d:definition-u} that $\nu=\mu_{| K_{\epsilon}}=0$. Since
$\epsilon >0$ was arbitrary and $\partial_M D\setminus\{{z_0}\}=\cup_{\epsilon >0} K^{z_0}_{\epsilon}$,
it follows that $\mu_{|\partial_M D\setminus\{{z_0}\}}=0$. Therefore $h=\mu(\{{z_0}\})
M_D(\cdot, {z_0})$ showing that $
M_D(\cdot, {z_0})$ is minimal.
\qed

\noindent
\textbf{Proof of Corollary \ref{c:main-mb0}} (a)
We first note that since $D$ is bounded, all Martin boundary points are finite, hence $\partial_M^f D=\partial_M D$. 
Let $\Xi:\partial D\to \partial^f_MD$ so that $\Xi(z)$ is the unique Martin boundary
point associated with $z\in \partial D$.
Since every finite Martin boundary point is associated with some $z\in \partial D$, we see that $\Xi$ is onto.
We show now that $\Xi$ is 1-1.
If not, there are $z, z'\in \partial D$, $z\neq z'$, such that $\Xi(z)=\Xi(z')=w$.
Then $M_D(\cdot, z)=M_D(\cdot, w)= M_D(\cdot, z')$. Choose $r>0$ small enough and
satisfying  $r<d(z, z')/4$. By \eqref{e:martinkernelub} and \eqref{e:martinkernelub3}
we see that there exists a constant $c_1=c_1(z)$ such that $M_D(x, z)\le c_1$ for
all $x\in D\setminus B(z,2r)$. Similarly, there exists $c_2=c_2(z')$  such that
$M_D(x, z')\le c_2$ for all $x\in D\setminus B(z',2r)$. Since $B(z, 2r)$ and
$B(z',2r)$ are disjoint, we conclude that $M_D(\cdot, z)=M_D(\cdot, z')$ is
bounded  on $D$ by $c_1\vee c_2$. Again by \eqref{e:martinkernelub},
$\lim_{D\ni x\to \zeta}M_D(x,z)=0$ for all regular $\zeta\in \partial D$.
In case of unbounded $D$, we showed in the proof of Theorem \ref{t:main-mb0}(b)
that $\lim_{x\to \infty}M_D(x,z)=0$. Hence by Lemma \ref{l:irregular}
we conclude that $M_D(\cdot, z)\equiv 0$. This is a contradiction with $M_D(x_0, z)=1$.

The statement about the minimal Martin boundary follows from Theorem 
\ref{t:main-mb0}(c).

\noindent (b)
We will show that  $\Xi:\partial D\to \partial^f_MD$ is actually a homeomorphism.
Let $z_0\in \partial D$ and $x\in D$. Choose $r<\frac12 \min\{R,\mathrm{dist}(x,z_0),\mathrm{dist}(x_0,z_0)\}$
so that $x\in D\setminus B(z_0, 2r)$.
It follows from Lemma \ref{l:limit-exists} that for any $s<1$ and $y\in D_{sr}$,
\begin{equation}\label{e:martin-estimate-boundary}
\Big|\frac{G_D(x, y)}{G_D(x_0, y)}-M_D(x, z_0)\Big|\le M_D(x, z_0)\left(\mbox{RO}_{D_{sr}}
\frac{G_D(x, \cdot)}{G_D(x_0, \cdot)}-1\right).
\end{equation}
Let $s<1$ and  $z'\in \partial D\cap B(z_0, sr/2)$.
It follows from Lemma \ref{l:limit-exists} that
there exists $M_D(x, z')=\lim_{D\ni y\to z'}M_D(x, y)$.
Letting $y\to z'$ in \eqref{e:martin-estimate-boundary} we get that
$$
|M_D(x, z')-M_D(x, z_0)|\le M_D(x, z_0)\left(\mbox{RO}_{D_{sr}}
\frac{G_D(x, \cdot)}{G_D(x_0, \cdot)}-1\right).
$$
Together with Proposition \ref{p:oscillation-reduction} we get
that if $(z_n)_{n\ge 1}$ is a sequence of points in $\partial D$ converging to
$z_0\in \partial D$, then $M_D(\cdot,z_0)=\lim_{n\to \infty}M_D(\cdot, z_n)$.

In order to show that $\Xi$ is continuous we proceed as follows.
Let $z_n\to z_0$ in $\partial D$.
Since $\partial_M D$ is compact, $(\Xi(z_n))_{n\ge 1}$ has
a subsequence $(\Xi(z_{n_k}))_{k\ge 1}$ converging in the Martin topology to some $w\in \partial_M D$.
By property (M3), $M_D(\cdot, \Xi(z_{n_k}))\to M_D(\cdot, w)$. On the
other hand, by the first part of the proof,
$M_D(\cdot, \Xi(z_{n_k}))=M_D(\cdot, z_{n_k})\to M_D(\cdot, z_0)$, implying that $w=\Xi(z_0)$.
This shows in fact that $(\Xi(z_n))_{n\ge 1}$ is convergent with the limit $\Xi(z_0)$.
Using the fact that $\partial D$ is compact, the proof of the continuity of the inverse is similar.

\noindent
(c) The Martin representation for non-negative harmonic functions is now a consequence of the general result 
from \cite{KW}, cf.~\eqref{e:martin-representation}.
\qed

\noindent
\textbf{Proof of Corollary \ref{c:finite-not-infinite}} 
(a) Assume that $w\in \partial_M^{\infty} D\cap \partial_M^f D$. Then there exist an unbounded sequence $(y_n)_{n\ge 1}\subset D$ and a bounded sequence $(z_n)_{n\ge 1}\subset D$ both converging to $w$ in the Martin topology. Since there is a subsequence $(y_{n_k})_{k\ge 1}$ such that $y_{n_k}\to \infty$, we have that $w=\infty$, i.e., $M_D(\cdot, w)=M_D(\cdot, \infty)$. Similarly, there is a subsequence $(z_{n_k})_{k\ge 1}$ and $z\in \partial D$ such that $z_{n_k}\to z$, hence $M_D(\cdot, w)=M(\cdot, z)$. This implies that $M_D(\cdot, \infty)=M_D(\cdot,z)$. We are going to show now that this is impossible. The proof of this fact is similar to the proof of Corollary \ref{c:main-mb0}(a). 

As in the proof of Theorem \ref{t:main-mb0}(c), choose $r$ small enough so that $M_D(x,z)\le c_1$ for all $x\in D\setminus B(z,2r)$, cf.~\eqref{e:martinkernelub} and \eqref{e:martinkernelub3}. Let $z_0\in\X$ be the point in the statement of Theorem \ref{t:main-mb1}. 
As in the proof of Theorem \ref{t:main-mb1}(c), choose $r'$ large enough satisfying
$r'>2(d(z, z_0)+4r)$ so that $M_D(x,\infty)\le c_2$ for all $x\in D\cap B(z_0, r'/2)$, cf.~\eqref{e:martinkernelub-I} and \eqref{e:martinkernelub3-i}.
Since $(D\setminus B(z,2r))\cup B(z_0, r'/2)=D$, we conclude that $M_D(\cdot, \infty)=M_D(\cdot, z)$ is bounded on $D$ by $c_1\vee c_2$. In the same way as in the proof of Corollary \ref{c:main-mb0}(a) we conclude that $M_D(\cdot, z)\equiv 0$ which is a contradiction.

\noindent (b) In the proof of Corollary \ref{c:main-mb0}(a)  we defined the mapping  $\Xi:\partial D\to \partial^f_MD$ and showed that it is 1-1 and onto. By inspecting the proof of Corollary \ref{c:main-mb0}(b), we can see that it carries over to the case when $D$ is unbounded. Hence, $\Xi$ is a homeomorphism from $\partial D$ to $\partial_M^f D$. Let $\partial D\cup \{\partial_{\infty}\}$ be the one-point compactification of $\partial D$. Extend $\Xi$ to this compactification by defining $\Xi(\partial_{\infty})=\infty\in \partial_M^{\infty}D$. By part (a), $\Xi$ is 1-1 and onto. 

Let $x\in D$. Choose 
$r>2 \max\{R,\mathrm{dist}(x,z_0),\mathrm{dist}(x_0,z_0)\}$
so that $x\in D\cap  B(z_0, r/2)$.
It follows from Lemma \ref{l:limit-exists-i} that for any $s>1$
\begin{equation}\label{e:martin-estimate-boundary-i}
\Big|\frac{G_D(x, y)}{G_D(x_0, y)}-M_D(x, \infty)\Big|\le M_D(x,\infty)\left(\mbox{RO}_{D^{sr}}
\frac{G_D(x, \cdot)}{G_D(x_0, \cdot)}-1\right), \qquad y\in D^{sr}.
\end{equation}
Let $s>1$ and  
$z'\in \partial D\cap B(z_0, 2sr)^c$.
It follows from Lemma \ref{l:limit-exists} that
there exists $M_D(x, z')=\lim_{D\ni y\to z'}M_D(x, y)$.
Letting $y\to z'$ in \eqref{e:martin-estimate-boundary-i} we get that
$$
|M_D(x, z')-M_D(x, \infty)|\le M_D(x, \infty)\left(\mbox{RO}_{D^{sr}}
\frac{G_D(x, \cdot)}{G_D(x_0, \cdot)}-1\right).
$$
Together with Proposition \ref{p:oscillation-reduction-I} we get
that if $(z_n)_{n\ge 1}$ is a sequence of points in $\partial D$ converging to
$\infty$, then $M_D(\cdot,\infty)=\lim_{n\to \infty}M_D(\cdot, z_n)$.
\ref{c:main-mb0}(b).
The rest of the proof of (b) and the proof of (c) is exactly the same as the proof of Corollary  
\ref{c:main-mb0}(b) and (c), respectively.

\section{Examples}\label{s:exam}
Several classes of Feller processes satisfying the assumptions of \cite{KSVp1} were studied in that paper. These examples include some symmetric and isotropic L\'evy processes in $\R^d$, strictly stable (not necessarily symmetric) processes in $\R^d$, processes obtained by subordinating a Feller diffusion on unbounded Ahlfors regular $n$-spaces, and space non-homogeneous processes on $\R^d$ whose Dirichlet form is comparable to 
the Dirichlet forms of certain subordinate Brownian motions. 
Since the conditions of the present paper are implied by the conditions of \cite{KSVp1}, we refer the readers to that paper for details. Here we will focus on certain symmetric and isotropic L\'evy processes where we can say more regarding accessible boundary points, and a class of subordinate Brownian motions not covered by \cite{JK}.

\subsection{Symmetric and isotropic L\'evy processes}
Let $X=(X_t, \P_x)$ be a purely discontinuous  symmetric L\'evy process in $\R^d$ with L\'evy exponent $\Psi(\xi)$ so that
$$
\E_x\left[e^{i\xi\cdot(X_t-X_0)}\right]=e^{-t\Psi(\xi)}, \qquad t>0, x\in \R^d, \xi\in\R^d.
$$
Thus the state space $\X=\R^d$, the measure $m$ is the $d$-dimensional Lebesgue measure and the localization radius $R_0=\infty$.
Assume that $r\mapsto j_0(r)$ is a strictly positive and nonincreasing function on $(0, \infty)$ satisfying
\begin{equation}\label{e:fuku1.1}
j_0(r)\le cj_0(r+1), \qquad r>1\, ,
\end{equation}
for some $c>1$ and that the L\'evy measure of $X$ has a density $j$ such that
\begin{equation}\label{e:fuku1.2}
\gamma^{-1}j_0(|y|)\le j(y) \le \gamma j_0(|y|), \qquad y\in \R^d
\end{equation}
for some $\gamma\ge 1$.
Since $\int_0^\infty j_0(r) (1\wedge r^2) r^{d-1}dr < \infty$ by \eqref{e:fuku1.2},
the function $x \to j_0(|x|)$ is the L\'evy density of  an isotropic unimodal L\'evy process whose characteristic exponent is
\begin{equation}\label{e:fuku1.3}
\Psi_0(|\xi|)= \int_{\R^d}(1-\cos(\xi\cdot y))j_0(|y|)dy.
\end{equation}
The L\'evy exponent $\Psi$ can be written as
$$
\Psi(\xi)= \int_{\R^d}(1-\cos(\xi\cdot y))j(y)dy
$$
and, clearly by \eqref{e:fuku1.2}, it satisfies
\begin{equation}\label{e:fuku1.4}
\gamma^{-1} \Psi_0(|\xi|)\le \Psi(\xi) \le \gamma \Psi_0(|\xi|),
\quad \mbox{for all } \xi\in \R^d\, .
\end{equation}

Under the above assumptions, the process $X$ satisfies Assumptions {\bf (A)} and {\bf (C)} (with $j(y,z)=\wh{j}(y,z)=j(z-y)$),
It also satisfies the assumption  {\bf (B)}, {\bf B1-a}$(0,R)$, {\bf B1-b}$(0,R)$,
{\bf B1-c}$(0,R)$, {\bf B2-a}$(0,R)$ of \cite{KSVp1} (for some $R>0$); see \cite{KSVp1} for more details.

Assume further that $\Psi_0$ satisfies the following scaling condition at infinity:

\medskip
\noindent
{\bf H1}:
There exist constants $0<\delta_1\le \delta_2 <1$ and $a_1, a_2>0$  such that
\begin{equation}\label{e:fuku1.6}
a_1\left(\frac{t}{s}\right)^{2\delta_1} \le \frac{\Psi_0( t)}{\Psi_0( s)} \le a_2 \left(\frac{t}{s}\right)^{2\delta_2} , \quad t \ge s \ge 1\, .
\end{equation}
Then by (15) and Corollary 22 in  \cite{BGR14}, for every $R>0$, 
there exists $c=c(R)>1$ such that
\begin{equation}\label{e:fuku1.7}
c^{-1}\frac{\Psi_0(r^{-1})}{r^d} \le j(r)\le c \frac{\Psi_0(r^{-1})}{r^d}
\quad \hbox{for } r\in (0, R]\, .
\end{equation}

Let $\Phi(r)=(\Psi_0(r^{-1}))^{-1}$.
Using \eqref{e:fuku1.1} and \eqref{e:fuku1.7}, one can easily see that there exists $R>0$ such that
Assumption \textbf{C1}$(0, R)$ is satisfied.
It is shown in Example 5.1 in \cite{KSVp1} that $X$ 
also satisfies assumptions {\bf C1}$(0,R)$ and {\bf D1}$(0,R)$  
of that paper (for some $R>0$).
Consequently, Theorem 4.1 of \cite{KSVp1} is valid which is precisely the assumption \textbf{F1}$(0, R)$. Further, it follows from 
Lemma 2.7 in \cite{KSV14b}  that \eqref{e:nGassup1} is also satisfied.
Using \textbf{F1}$(0, R)$ and the fact that open balls
are Greenian, we can apply Proposition 6.5 in \cite{KSVp1}.
Thus for any Greenian open set
  $D$, $\lim_{x\to z}G_D(x,y)=0$ for every regular point $z\in \partial D$, so \eqref{e:zreg} holds.
 In case of an unbounded $D$ we assume that $X$ is transient. Then 
 $\lim_{x\to \infty}G_D(x,y)=0$ by Lemma 2.10 in \cite{KSV14b}.
 We conclude that Theorem \ref{t:main-mb0} and Corollary 
 \ref{c:main-mb0} apply.

Instead of {\bf H1}, assume that $\Psi_0$ satisfies the following scaling condition at zero:

\medskip
\noindent
{\bf H2}:
There exist constants $0<\delta_3\le \delta_4 <1$ and $a_3, a_4>0$  such that
\begin{equation}\label{e:fuku1.6b}
a_3\left(\frac{t}{s}\right)^{2\delta_3} \le \frac{\Psi_0( t)}{\Psi_0( s)} \le a_4 \left(\frac{t}{s}\right)^{2\delta_4} , \quad s \le t\le 1\, .
\end{equation}
It is shown in Example 5.1 of \cite{KSVp1}
that for every $R>0$ there exists $c=c(R)>1$ such that
\begin{equation}\label{e:fuku1.7gl}
c^{-1}\frac{\Psi_0(r^{-1})}{r^d} \le j_0(r)\le c \frac{\Psi_0(r^{-1})}{r^d}
\quad \hbox{for } r\in [R, \infty).
\end{equation}
Together with \eqref{e:fuku1.2} this implies that there is $R>0$ such that {\bf C2}$(0,R)$ is true.

Let $d\ge 3$.
Then $X$ is transient and let $G(x)=G(x,0)$ be its Green function.
Then by Lemma 2.10 of \cite{KSV14b},  \eqref{e:nGassup1-infty} holds. 
Assume that there exists a
non-increasing
 function $r\mapsto G_0(r)$ and a constant $c\ge 1$ such that
\begin{align}
\label{e:Gcomp}
c^{-1}G_0(|x|)\le G(x) \le c G_0(|x|)\, ,\qquad x\in \R^d\, .
\end{align}
It is shown in Example 5.1 of \cite{KSVp1} that $X$ also satisfies assumptions
{\bf B2-b}$(0,R)$, {\bf C2}$(0,R)$ and {\bf D2}$(0,R)$  of that paper (for some $R>0$). Consequently, Theorem 2.1 of \cite{KSVp1} is 
valid which is precisely the assumption \textbf{F2}$(0, R)$.
If we assume that the Green function of $X$ is continuous then using 
the upper bound $G(x)\le c {|x|^{-1}\Psi_0(|x|^{-1})^{-1}}$ in  
(5.16) in \cite{KSVp1} and the strong Markov property,
 the Green function of $X^D$ is continuous  for all open set $D$. Thus
by Proposition 6.2 of \cite{KSVp1},  \eqref{e:zreg} holds.
Further, it follows from
(5.16) in \cite{KSVp1}
that, if $d \ge 3$, \eqref{e:nGassup2i} is also satisfied.
We conclude that Theorem \ref{t:main-mb1} applies for $d \ge 3$ 
under the assumption that $G$ is continuous and satisfies \eqref{e:Gcomp}.
In fact, it is also  shown in Example 5.1 of \cite{KSVp1} that, if $X$
 is a subordinate Brownian motion
whose
Laplace exponent $\phi$ is a complete Bernstein function
and that  $\xi \to \phi(|\xi|^2)$ satisfies Assumption {\bf H2}, then \eqref{e:nGassup2i} is satisfied for $d >2 \delta_4$.
Since $G(x)=g(|x|)$ is continuous and $r\mapsto g(r)$ is decreasing,  in this case 
Theorem \ref{t:main-mb1} applies for $d >2 \delta_4$.

\smallskip
In the next proposition we give a criterion for the accessibility of infinity and a finite boundary point. Let $B_r=B(0,r)$ and for an open set $D$, $D_r=D\cap B_r$ and $D^r=D\cap \overline{B}_r^c$.
\begin{prop}\label{p:criterion-accessibility}
(a) Let $D\subset \R^d$ be a Greenian open set such that $0\in \partial D$ and assume that {\bf H1} holds. Then $0$ is inaccessible from $D$ with respect to $X$ if and only if
\begin{equation}\label{e:criterion-accessibility-finite}
\int_{D_1} (\E_y \tau_{D_1}) j(y)\, dy <\infty\, .
\end{equation}

\noindent
(b) Let $D\subset  \R^d$ be a Greenian open set and assume that {\bf H2} holds. Then $\infty$ is inaccessible from $D$ with respect to $X$ if and only if
\begin{equation}\label{e:criterion-accessibility-infty}
\int_{D^1} P_{D^1}(y,0)\, dy <\infty\, .
\end{equation}
\end{prop}
\pf (b) Recall that $\infty$ is inaccessible from $D$ if there exists $x\in D$ such that $\E_x \tau_D<\infty$. Let $r=\max(2|x|,R,1)$ where $R>0$ is the constant from \textbf{C2}$(0, R)$ and \textbf{F2}$(0, R)$. We write
\begin{align*}
\E_x \tau_D &= \int_{B(x,4r)}G_D(x,y)\, dy + \int_{D_{8r}\setminus B(x,4r)}G_D(x,y)\, dy+ \int_{D^{8r}}G_D(x,y)\, dy \\
&=:I+II+III.
\end{align*}
Since $X^D$ is transient, $I=G_D \ind_{B(x,4r)}$ is bounded, hence finite. By \eqref{e:nGassup1-infty} we have that $G_D(x,y)\le c(r)$ for $y\in D_{8r}\setminus B(x,4r)$, hence
$II\le c(r)|D_{8r}|<\infty$.  Since $G_D(x,\cdot)$ is regular harmonic in $D^r$, it follows from \textbf{F2}$(0, R)$ that $G_D(x,y)\asymp P_{D^r}(y,0)$ for all $y\in D^{8r}$. Thus $III\asymp \int_{D^{8r}}P_{D^r}(y,0)\, dy$. Hence, $\E_x \tau_D<\infty$ if and only if $\int_{D^{8r}}P_{D^r}(y,0)\, dy <\infty$. Next,
\begin{eqnarray*}
\int_{D^1}P_{D^1}(y,0)\, dy&=&\int_{D\cap \{1<|y|\le 8r\}}P_{D^1}(y,0)\, dy+ \int_{D^{8r}}P_{D^1}(y,0)\, dy\\
& \le & \int_{\int_{D\cap \{1<|y|\le 8r\}}}P_{B^1}(y,0)\, dy + \int_{D^{8r}}P_{D^1}(y,0)\, dy\\
& =& : IV + V\, .
\end{eqnarray*}
By Proposition 3.1 of \cite{KSVp1}, $P_{B^1}(y,0)\le c_1$ for all $y\in B^1$, 
hence $IV\le c_1 |B_{8r}|<\infty$. Finally, by repeatedly 
applying Lemma 3.9 of \cite{KSVp1} we deduce that
$$
P_{D^r}(y,0)\le P_{D^1}(y,0) \le c_2 P_{D^r}(y,0)\, , \quad y\in D^{8r},
$$
with a constant $c_2>0$ depending on $r$. Thus, $V$ is comparable to $\int_{D^{8r}} P_{D^r}(y,0)\, dy$, proving that $\int_{D^1}P_{D^1}(y,0)\, dy<\infty$ if and only if $\int_{D^{8r}} P_{D^r}(y,0)\, dy<\infty$. This finishes the proof.

\noindent
(a) This can be proved in the similar way, so we omit the proof. \qed

\begin{remark}\label{r:criterion-accessibility}
{\rm
(a) Note that the criterion in Proposition \ref{p:criterion-accessibility} does not depend on $x\in D$. Hence, if $\E_x \tau_D <\infty$ for one $x\in D$, then $\E_x \tau_D <\infty$ for all $x\in D$. Similarly, if $P_D(x,0)<\infty$ for one $x\in D$, then $P_D(x,0)<\infty$ for all $x\in D$.

\noindent
(b) By inspecting the proof of Proposition \ref{p:criterion-accessibility} one can see that it carries over to the case of the process satisfying the assumptions in \cite{KSVp1}.  In particular, $x\mapsto \E_x\tau_D$ (respectively $x\mapsto P_D(x,z_0)$) is either identically infinite or finite for all $x\in D$.
}
\end{remark}

For any open set $V$, let $s_V(x)=\E_x \tau_V$ and let $\omega_V^x=\P_x(X_{\tau_V}\in \cdot)$ be the harmonic measure. By the strong Markov property, for $V\subset D$ we have
$$
s_D(x)=s_V(x)+\int_{D\setminus V}s_D(y)\omega_V^x(dy)\, .
$$
If $\partial V\cap D$ is Lipschitz, it follows from \cite{Szt} that $\omega_V^x(\partial V)=0$ and hence
\begin{equation}\label{e:sVD}
s_D(x)=s_V(x)+\int_{D\setminus V}s_D(y)P_V(x,y)dy\, .
\end{equation}

We now record the following lower bound on the expected exit time from a ball: There exists a constant $c>0$ such that for every $r>0$ and every $x\in \R^d$
\begin{equation}\label{e:exit-ball-lower}
\E_x \tau_{B(x,r)}\ge \frac{c}{\Psi_0(r^{-1})}\, .
\end{equation}
This follows, for example, from the last display in the proof of 
Theorem 2.2 of \cite{CK15} and the proof of Lemma 13.4.2 in \cite{KSV12b}.

Let $\kappa\in (0,1/2]$. Recall that an open set $D$ in $\R^d$ is said to be $\kappa$-fat at $z_0\in \partial D$ if there exists $r_0>0$ such that for every $r\in (0,r_0]$ there exists $A_r\in D$ such that $B(A_r, \kappa r)\subset D\cap B(0,r)$.
An open set $D$ in $\R^d$ is said to be $\kappa$-fat at infinity if there exists $r_0>0$ such that for every $r\ge r_0$ there exists $A_r\in D$ such that $B(A_r, \kappa r)\subset D\cap \overline{B}(0,r)^c$ and $|A_r|<\kappa^{-1}r$, 
cf. Definition 1.3 in \cite{KSV14}.

\begin{prop}\label{p:kappa-fat-infty}
(a) Suppose that {\bf H1} holds. If $D\subset \R^d$ is $\kappa$-fat at $z_0\in \partial D$, then $z_0$ is accessible from $D$ with respect to $X$.

\noindent
(b)  Suppose that {\bf H2} holds.  If $D\subset \R^d$ is $\kappa$-fat at infinity, then infinity is accessible from $D$ with respect to $X$.
\end{prop}
\pf  We prove part (b). The proof of (a) is similar.

Let $A_0=A_{r_0}$ be a point in $D$ such that $B(A_0, \kappa r_0)\subset D\cap \overline{B}(0,r_0)^c$ and $|A_0|<\kappa^{-1}r_0$.
We inductively define the sequence $r_n=4\kappa^{-1}r_{n-1}$, $n\ge 1$, and a sequence of points $A_n=A_{r_n}$ such that $B(A_n,\kappa r_n)\subset D\cap \overline{B}(0,r_n)^c$ and $r_n<|A_n|<\kappa^{-1}r_n$. It is easy to see that the family of balls $(B(A_n, \kappa r_n)_{n\ge 0}$ is pairwise disjoint.

Let $U:=\cup_{n=0}^{\infty}B(A_n, \kappa r_n)$. Then by \eqref{e:sVD} with $V=B(A_0,\kappa r_0)$,
\begin{eqnarray*}
\E_{A_0}\tau_D &\ge & \int_{D\setminus B(A_0, \kappa r_0)}s_D(y)P_{B(A_0,\kappa r_1)}(A_10,y)\, dy \\
&\ge & \int_{D\setminus B(A_0,\kappa r_0)}s_U(y)P_{B(A_0, \kappa r_0)}(A_0,y)\, dy\\
&\ge &\sum_{n=1}^{\infty}\int_{B(A_n, \kappa r_n/2)}s_{B(A_n,\kappa r_n)}(y) P_{B(A_0,\kappa r_1)}(A_0,y)\, dy \\
&\ge & \sum_{n=1}^{\infty}\big(\inf_{y\in B(A_n,\kappa r_n/2)}s_{B(A_n,\kappa r_n)}(y)\big) \int_{B(A_n,\kappa r_n/2)}P_{B(A_0,\kappa r_1)}(A_0,y)\, dy.
\end{eqnarray*}
By \eqref{e:exit-ball-lower}, 
$$
s_{B(A_n,\kappa r_n)}(y)\ge c_1 \Psi_0((\kappa r_n)^{-1})^{-1}\ge c_2 \Psi_0(r_n^{-1})^{-1}
$$ 
for all $y\in B(A_n,\kappa r_n/2)$. Further, if $y\in B(A_n,\kappa r_n)$, then $r_n/2\le |y-A_0|\le 3\kappa^{-1}$. Hence,
by Lemma 3.3 of \cite{KSV12b}, \eqref{e:exit-ball-lower}, \eqref{e:fuku1.2} 
and \eqref{e:fuku1.7gl}, we have that for $y\in B(A_n,\kappa r_n/2)$,
$$
P_{B(A_0,\kappa r_0)}(A_0,y)\ge c_3 \frac{\Psi_0(|y-A_0|^{-1})}{|y-A_0|^d} \Psi_0((\kappa r_0)^{-1})^{-1}\ge c_4 \frac{\Psi_0(r_n^{-1})}{r_n^d}\Psi_0(r_1^{-1})^{-1}\, .
$$
Therefore,
$$
\E_{A_0}\tau_D \ge  \sum_{n=2}^{\infty} c_5 \Psi_0(r_n^{-1})^{-1} \frac{\Psi_0(r_n^{-1})}{r_n^d}\Psi_0(r_1^{-1})^{-1} r_n^d =\infty\, .
$$
By using Remark \ref{r:criterion-accessibility} we see that $\infty$ is accessible form $D$. \qed

In the next result we give a criterion for accessibility of infinity from a thorn-like domain.

Let $f:(2, \infty) \to (0, \infty)$  be a positive non-decreasing function such that  $f(t) \le t$ for all $t>0$ and define
$$
D=D^f:=\{(y_1, \widetilde{y})\in \R^d: y_1>2, |\widetilde{y}|<f(y_1)\}\, .
$$
Here $y=(y_1,\wt{y})$ with $\wt{y}=(y_2,\dots, y_d)\in \R^{d-1}$.
\begin{prop}\label{p:thorn-infty}
Suppose that {\bf H1} and {\bf H2} hold. Then
 infinity is accessible from $D$ if and only if
\begin{equation}\label{e:thorn-infty}
\int_4^{\infty} \frac{\Psi_0(t^{-1})}{\Psi_0(f(t)^{-1})}\, \frac{f(t)^{d-1}}{t^{d}}\, dt=\infty\, .
\end{equation}
\end{prop}

\pf
Assume that the integral in \eqref{e:thorn-infty} is infinite. Fix $x\in D$ and denote $\delta_D(x)$ by $r$. Let $U:=\{(y_1, \widetilde{y})\in \R^d: y_1>4(1+|x|), |\widetilde{y}|<f(y_1)/2\}.$
Since $|x-y| \ge y_1-x_1 \ge 4(1+|x|)-x_1 > 2x_1> 2f(x_1) > r$ for all $y \in U$, we have $U\subset D \setminus B(x, r)$.

Moreover, $B(y, f(y_1/2)/2) \subset D $ for all $y \in U$. In fact, for $z \in B(y, f(y_1/2)/2)$ with $y \in U$ we have
$z_1> y_1- f(y_1/2)/2>y_1/2$, which implies that $z_1> 6$ and $f(y_1/2) \le f(z_1)$. Using the last inequality we see that
$
 |\widetilde{z}| \le  |\widetilde{y}- \widetilde{z}|+ |\widetilde{y}| <f(y_1/2)\le f(z_1)
$.
Thus for $y\in U$,
$$
s_D(y)\ge s_{B(y, f(y_1/2)/2)}(y)\ge \frac{c_1}{\Psi_0((f(y_1/2)/2)^{-1})}\, ,
$$
where the last inequality follows from
\eqref{e:exit-ball-lower}.

Notice that for $y\in U$, $|y|\asymp y_1$.  Thus, since $|z-y| \le  |x|+|z-x|+|y| \le 6 y_1$ for $z \in B(x,r)$,
using $j(y_1)\asymp \frac{\Psi_0(y_1^{-1})}{y_1^d}$ we have
$$
P_{B(x, r)}(x, y)\ge c_2\E_x[\tau_{B(x,r)}] j(y_1) \asymp \frac{\Psi_0(y_1^{-1})}{y_1^d}, \qquad y\in U.
$$
Therefore
\begin{eqnarray*}
s_D(x)&\ge & \int_U s_D(y)P_{B(x, r)}(x, y)dy\\
&\ge& c_3\int^\infty_{4(1+|x|)}\frac{\Psi_0(y_1^{-1})}{\Psi_0((f(y_1/2)/2)^{-1})}\frac{f(y_1/2)^{d-1}}{y_1^d} dy_1\\
&=&c_4 \int_{2(1+|x|)}^{\infty} \frac{\Psi_0(2^{-1}t^{-1})}{\Psi_0(2f(t)^{-1})}\, \frac{f(t)^{d-1}}{t^d}dt\\
&\ge & c_5\int_{2(1+|x|)}^{\infty}\frac{f(t)^{d-1}}{t^{d}}\, \frac{\Psi_0(t^{-1})}{\Psi_0(f(t)^{-1})}dt=\infty\, ,
\end{eqnarray*}
where the last inequality follows from Lemma 1 of \cite{G14}.
Thus $\infty$ is accessible from $D$.

Assume that the integral in \eqref{e:thorn-infty} is finite.
For $r\ge 4$, let $D_r:=D \cap B(0, r)$. Then,  
by Lemma 2.5 and (2.1) in \cite{KSV14b}, and 
considering the infinite cylinder, we get
$\sup_{x_1=r} s_{D_{4r}}(x) \le c_6 \Psi_0(f(4r)^{-1})^{-1}$.
Thus, by \eqref{e:sVD},  we have that for $x\in D$ with $x_1=r$
\begin{align}
s_D(x)=&s_{D_{4r}}(x)+\int_{D\setminus D_{4r}}s_D(y)P_{D_{4r}}(x, y)dy
\label{swqe0}\\
\le &c_6 \Psi_0(f(4r)^{-1})^{-1}+\int_{D\setminus D_{4r}}s_D(y)P_{D_{4r}}(x, y)dy. \label{swqe}
\end{align}
By the argument in the paragraph before Theorem 3.12 in \cite{KSV14a},
 Lemma 5.4 in \cite{KSV12}  is valid for all $r>0$. Hence
$$
P_{D_{4r}}(x, y) \le c_7 s_{D_{4r}}(x) \left(\int_{D \setminus B(0, 2r)} j(|z|)P_{D_{4r}}(z, y)dz+ j(|y|) \right)
$$
Thus using \eqref{swqe0}
\begin{eqnarray*}
\lefteqn{\int_{D\setminus D_{4r}}s_D(y) P_{D_{4r}}(x, y) dy}\\
&\le &c_8s_{D_{4r}}(x)\left( \int_{D\setminus D_{4r}}s_D(y) \int_{D \setminus B(0, 2r)} j(|z|)P_{D_{4r}}(z, y)dzdy\right.\\
&&\qquad\quad \left.+\int_{D\setminus D_{4r}}s_D(y) j(|y|)dy \right)\\
&=&c_8s_{D_{4r}}(x)\left(
 \int_{D \setminus B(0, 2r)} j(|z|)\big(\int_{D\setminus D_{4r}}s_D(y)P_{D_{4r}}(z, y)dy\big)dz\right.\\
 &&\qquad \quad\left.+c_8\int_{D\setminus D_{4r}}s_D(y) j(|y|)dy\right)\\
& \le &c_8s_{D_{4r}}(x)\left(
 \int_{D \setminus B(0, 2r)} j(|z|)s_D(z)dz+\int_{D\setminus D_{4r}}s_D(y) j(|y|)dy\right)\\
&\le &2c_8s_{D_{4r}}(x)
 \int_{D \setminus B(0, 2r)} j(|z|)s_D(z)dz\\
 & \le& 2c_6c_8\Psi_0(f(4r)^{-1})^{-1}
 \int_{D \setminus B(0, 2r)} j(|z|)s_D(z)dz.
\end{eqnarray*}
Applying this to \eqref{swqe}, we get
\begin{align}\label{swqe1}
s_D(x)
\le c_9 \Psi_0(f(4r)^{-1})^{-1}\left(1+\int_{D \setminus B(0, 2r)} j(|z|)s_D(z)dz\right).
\end{align}
Let
$M(r):=\sup_{x_1=r} s_D(x) \Psi_0(f(4r)^{-1})$.
From \eqref{swqe1}, for $r>4$,
\begin{align*}M(r)
&\le c_{10} \left(1+\int_{2r}^\infty \int_{|\wt z| < f(s)} |(s, \wt z)|^{-d} \Psi_0(|(s, \wt z)|^{-1})   M(s) \Psi_0(f(4s)^{-1})^{-1}  d \wt zds\right)\\
&\le c_{11} \left(1+\int_{2r}^\infty  f(s)^{d-1} s^{-d} \Psi_0(s^{-1})   M(s) \Psi_0(f(4s)^{-1})^{-1}  ds\right)\\
&\le c_{11} \left(1+\int_{r}^\infty  f(s)^{d-1} s^{-d} \Psi_0(s^{-1})   M(s) \Psi_0(f(4s)^{-1})^{-1}  ds\right).
\end{align*}
Let
$m(r)=M(1/r)$; by a change of variable we have that for $r<1/4$,
\begin{align*}m(r)
&\le c_{11} \left(1+\int_{1/r}^\infty  f(s)^{d-1} s^{-d} \Psi_0(s^{-1})   M(s) \Psi_0(f(4s)^{-1})^{-1}  ds \right)\\
&= c_{11}\left(1+\int_{0}^r  f(v^{-1})^{d-1} v^{d} \Psi_0(v)   m(v) \Psi_0(f(4v^{-1})^{-1})^{-1}  v^{-2} dv\right).
\end{align*}
By  Gronwall's inequality,  for $r<1/4$,
$$
m(r) \le c_{12} \exp \left(\int_{0}^r  f(v^{-1})^{d-1} v^{d} \Psi_0(v)   \Psi_0(f(4v^{-1})^{-1})^{-1}  v^{-2} dv\right).
$$
Therefore, under the assumption that the integral in \eqref{e:thorn-infty} is finite, we have for all $x \in D$ with $x_1=r>4$,
\begin{align*}
s_D(x) &\le \Psi_0(f(4r)^{-1})
M(r)=\Psi_0(f(4r)^{-1}) m(1/r)
\\
&
\le  c_{12} \Psi_0(f(4r)^{-1}) \exp \left(\int_{r}^\infty  f(s)^{d-1} s^{-d} \Psi_0(s^{-1})   \Psi_0(f(4s)^{-1})^{-1}  ds\right)
\\
&\le  c_{12}\Psi_0(f(16)^{-1})  \exp \left(\int_{4}^\infty  f(4s)^{d-1} s^{-d} \Psi_0(s^{-1})   \Psi_0(f(4s)^{-1})^{-1}  ds\right)\\
&\le  c_{12}\Psi_0(f(16)^{-1})  \exp \left( c_{13}\int_{16}^\infty  f(t)^{d-1} t^{-d} \Psi_0(4t^{-1})   \Psi_0(f(t)^{-1})^{-1}  ds\right)\\
&\le  c_{12}\Psi_0(f(16)^{-1})  \exp \left(c_{14} \int_{16}^\infty  f(t)^{d-1} t^{-d} \Psi_0(t^{-1})   \Psi_0(f(t)^{-1})^{-1}  ds\right)< \infty.
\end{align*}
Here the last inequality follows from Lemma 1 of  \cite{G14}.
Hence infinity is inaccessible. \qed

Suppose that $f(t)=t(\log t)^{-\beta}$, $\beta \ge 0$. Then
\begin{eqnarray*}
I &=&\int_4^{\infty}\frac{\Psi_0(t^{-1})}{\Psi_0(t^{-1}(\log t)^{\beta})} \, \frac{t^{d-1}(\log t)^{-\beta(d-1)}}{t^d}\, dt\\
&=&\int_4^{\infty} \frac{\Psi_0(t^{-1})}{\Psi_0(t^{-1}(\log t)^{\beta})} (\log t)^{-\beta(d-1)}\, \frac{dt}{t}\\
&\ge & c_1\int_4^{\infty} (\log t)^{-\beta(2\delta_4+d-1)} \frac{dt}{t},
\end{eqnarray*}
where the inequality follows from {\bf H2}. When $\beta\le 1/(d-1+2\delta_4)$, the integral above is divergent and hence infinity is accessible. Note that when $\beta >0$, $D$ is not $\kappa$-fat at infinity for any $\kappa\in (0,1/2]$.
Similarly,
$$
I\le c_2 \int_4^{\infty} (\log t)^{-\beta(2\delta_3+d-1)} \frac{dt}{t}\, .
$$
When $\beta > 1/(d-1+2\delta_3)$, the integral above is convergent and hence the infinity is inaccessible.

A result analogous to Proposition \ref{p:thorn-infty} is valid for a finite boundary point.
Let $f:(0,1)\to (0,\infty)$ be a bounded increasing function such that $f(t)\le t$ and define
$$
D_f:=\{x=(x_1,\wt{x}): 0<x_1<1, |\wt{x}|<f(x_1)\}.
$$
\begin{prop}\label{p:thorn-finite}
Assume that {\bf H1} holds. Then the point $0$ is accessible from $D$ if and only if
\begin{equation}\label{e:thorn-finite}
\int_0^{1} \frac{\Psi_0(t^{-1})}{\Psi_0(f(t)^{-1})}\, \frac{f(t)^{d-1}}{t^{d}}\, dt=\infty\, .
\end{equation}
\end{prop}

\subsection{Subordinate Brownian motions}
Let $Y=(Y_t, \P_x)$ be a standard Brownian motion in $\R^d$, and $S=(S_t)$ an independent subordinator with the Laplace exponent $\phi$, $\E[e^{-\lambda S_t}]=e^{-t\phi(\lambda)}$. The subordinate Brownian motion $X=(X_t,\P_x)$ is defined as $X_t=Y(S_t)$. Assume that $\phi$ is a complete Bernstein function with infinite L\'evy measure $\mu$ satisfying the following hypothesis

\smallskip
\noindent
{\bf H}:  There exist constants $\sigma>0$, $\lambda_0 > 0$ and $\delta \in (0, 1]$ such that
\begin{equation*}
  \frac{\phi'(\lambda t)}{\phi'(\lambda)}\leq\sigma\, t^{-\delta}\ \text{ for all }\ t\geq 1\ \text{ and }\
 \lambda \ge \lambda_0\, .
\end{equation*}
When $d \le 2$,  assume that $d+2\delta-2>0$
and  there are $\sigma'>0$ and
\begin{equation}\label{e:new22}
\delta'  \in  \left(1-\tfrac{d}{2}, (1+\tfrac{d}{2})\wedge (2\delta+\tfrac{d-2}{2})\right)
\end{equation}
 such that
\begin{equation}\label{e:new23}
\frac{\phi'(\lambda x)}{\phi'(\lambda)}\geq \sigma'\,x^{-\delta'}\ \text{ for all}\ x\geq 1\ \text{ and }\ \lambda\geq\lambda_0\,;
\end{equation}

Assumption {\bf H}  was introduced and used in \cite{KM} and \cite{KM2}.  It is easy to check that if $\phi$ is a complete Bernstein function satisfying a weak lower scaling condition at infinity
\begin{equation}\label{e:new2322}
a_1 \lambda^{\delta_1}\phi(t)\le \phi(\lambda t)\le a_2 \lambda^{\delta_2}\phi(t)\, ,\qquad \lambda \ge 1, t\ge 1\, ,
\end{equation}
with $a_1, a_2>0$ and $\delta_1, \delta_2\in (0,1)$, then {\bf H}  is automatically satisfied. In that case the process $X$ belongs to the class of isotropic unimodal L\'evy process considered in the previous subsection. The reason for assuming hypothesis {\bf H} here is to cover the case of geometric stable and iterated geometric stable subordinators. Suppose that $\alpha\in (0, 2)$ for $d \ge 2$ and that $\alpha\in (0, 2]$ for $d \ge 3$. A geometric $(\alpha/2)$-stable subordinator is a subordinator with Laplace exponent $\phi(\lambda)=\log(1+\lambda^{\alpha/2})$. Let $\phi_1(\lambda):=\log(1+\lambda^{\alpha/2})$, and for $n\ge 2$, $\phi_n(\lambda):=\phi_1(\phi_{n-1}(\lambda))$. A subordinator with Laplace exponent $\phi_n$ is called an iterated geometric subordinator. It is easy to check that the functions $\phi$
and $\phi_n$ satisfy {\bf H}, but they do not satisfy \eqref{e:new2322}.

The process $X$ clearly satisfies assumption \textbf{A} and {\bf C},
and by symmetry, every semipolar set is polar.  Suppose that $X$ is transient. 
Then it follows from Lemma 5.4 of \cite{KM2} that for all $z_0\in \R^d$, 
 \textbf{C1}$(z_0, R)$, \textbf{F1}$(z_0, R)$, and \eqref{e:nGassup1} (with a uniform constant) hold true.
Moreover, since all Green functions are continuous, 
by Proposition 6.2 of \cite{KSVp1},
$\lim_{x\to z}G_D(x,y)=0$ for every regular boundary point $z$ of $\partial D$.
Therefore the conclusions of Corollary \ref{c:main-mb0} hold true.

Suppose now that $X$ is an (iterated) geometric $\alpha$-stable process with $0<\alpha <2$. Then $X$ satisfies condition {\bf H2} from the previous subsection 
(see Example 5.1 of \cite{KSVp1}) and 
by the same arguments we conclude that Theorem \ref{t:main-mb1} is true.

\section{Minimal thinness is a local property   }\label{s:mt-local}

The purpose of this section is to establish several results analogous to those 
in Section 9.5 of \cite{AG} and to conclude 
that minimal thinness is a local property.

The setting is the following: $(\X,d,m)$ is a metric measure space with countable base as before. Since bounded closed sets are compact, the topology of $\X$ is locally compact. Let $X=(X_t,\P_x)$ be a Hunt process in $\X$ satisfying
Assumption \textbf{A}.
The cone of excessive functions with respect to $X$ is denoted by $\SS(X)$. We assume that $(\X, \SS(X))$ is a balayage space in the sense of \cite{BH}.
Let $D\subset \X$
 be an open set, $X^D$ the killed process and $\SS(X^D)$ the cone of
excessive function with respect to $X^D$.
By Proposition V.1.1 and Proposition VI.3.20 of \cite{BH},
$(D,\SS(X^D))$ is also a balayage space in the sense of \cite{BH}.
In particular, all functions
in $\SS(X^D)$ are lower semi-continuous (l.s.c.) Moreover, by definitions
and results from 
 p.~94 and Lemma III 1.2  of \cite{BH}, 
bounded harmonic functions 
on $D$ with respect to $X^D$ are
continuous. Since we will be interested only in $X^D$, all notions
defined below are relative to $X^D$.

For any (numerical) function $f:D\to (-\infty,\infty]$ we define its lower semi-continuous (l.s.c.) regularization $\wh{f}$ by
$$
\wh{f}(x)=f(x)\wedge\left(\liminf_{y\to x}f(y)\right)\, .
$$
Then $\wh{f}$ is the largest l.s.c.~function dominated by $f$: $\wh{f}\le f$.
We remark that in this section the hat  $\ \widehat{ }\ $ denotes the l.s.c~regularization and \emph{not} the notions related to the dual process.
For a Borel set $A\subset D$, let $S_A=\inf\{t\ge 0:\, X_t\in A\}$ be the debut of $A$,
and $T_A=\inf\{t>0:\, X_t\in A\}$ the hitting time of $A$. For $u\in \SS(X^D)$, the reduced function of $u$ on $A$ is defined as 
(see p.243 of \cite{BH}):
\begin{eqnarray*}
R_u^A&=&\inf\{v\in \SS(X^D):\, v\ge u \textrm{ on }A\}\\
&=&\inf\{v\in \SS(X^D):\, v\le u, v=u \textrm{ on }A\}\, .
\end{eqnarray*}
Its l.s.c.~regularization $\wh{R}_u^A:=\wh{R_u^A}$ is called the balayage of $u$ on $A$. Then $\wh{R}_u^A\in \SS(X^D)$. The probabilistic interpretations of the reduced function and the balayage 
are (cf. VI.3 of \cite{BH})
$$
R_u^A(x)=\E_x[u(X_{S_A})]\, ,\qquad \wh{R}_u^A(x)=\E_x[u(X_{T_A})]\, .
$$

We have the following properties of $R_u^A$ and $\wh{R}_u^A$: $R_u^A=u$ on $A$, $\wh{R}_u^A\le R_u^A\le u$ 
(p.243 of \cite{BH}), 
$\wh{R}_u^A=R_u^A$ on $A^c$ (Proposition VI.2.3 of \cite{BH}), 
$\{\wh{R}_u^A<R_u^A\}$ is
semipolar (Proposition VI.5.11 of \cite{BH}), hence polar by \textbf{A}.

\medskip
Let $u:D\to [0,\infty)$ be continuous and harmonic in $D$ 
with respect to $X^D$, $E\subset D$ an open set, and $w:E\to [0,\infty)$ harmonic in $E$ with respect to $X^E$ such that
$w\le u-R_u^{D\setminus E}$. We set $w\equiv 0$ on $D\setminus E$.

\begin{lemma}\label{l:w-harmonic}
For every bounded open set $U\subset \overline{U}\subset D$, it holds that
$$
w(x)=\E_x[w(X_{\tau_{U\cap E}})]\, , \quad x\in U\cap E\, .
$$
\end{lemma}
\pf We first show that there exists a polar set $N\subset \partial E\cap D$ such that for every $z\in (\partial E\cap D)\setminus N$,
\begin{equation}\label{e:w-to-0}
\lim_{x\to z, x\in E}w(x)=0\, .
\end{equation}
Note that $R_u^{D\setminus E}(x)=\wh{R}_u^{D\setminus E}(x)$ for all $x\in E$. Hence by continuity of $u$ and lower semi-continuity of $\wh{R}_u^{D\setminus E}$,
\begin{align*}
\limsup_{x\to z, x\in E}w(x)&\le \limsup_{x\to z, x\in E}\big(u(x)-\wh{R}_u^{D\setminus E}(x)\big)\\
&=u(z)-\liminf_{x\to z, x\in E}\wh{R}_u^{D\setminus E}(x) \le u(z)-\wh{R}_u^{D\setminus E}(z)\, .
\end{align*}
Let $N=\partial E\cap D\cap \{\wh{R}_u^{D\setminus E}<R_u^{D\setminus E}\}$. Then $N$ is polar and it follows from the last display that for all $z\in (\partial E\cap D)\setminus N$ we have
$$
\limsup_{x\to z, x\in E}w(x)\le u(z)-R_u^{D\setminus E}(z)=0\, .
$$

For each $n\ge 1$ define $U_n:=\{x\in U\cap E:\, d(x,E^c)>\frac1n\}$. Then $U_n$ is bounded and open in $E$, $U_n\subset \overline{U}_n\subset U\cap E$, $\overline{U}_{n+1}\subset U_n$, and $U\cap E=\cup_{n=1}^{\infty}U_n$. By harmonicity of $w$, for any $x\in U\cap E$ and $n$ large enough,
\begin{eqnarray*}
w(x)&=&\E_x[w(X^E_{\tau_{U_n}})]=\E_x[w(X^E_{\tau_{U_n}}):\, \tau_{U_m}=\tau_{U\cap E} \textrm{ for some }m\ge 1]\\
& &+\E_x[w(X^E_{\tau_{U_n}}):\, \tau_{U_m}<\tau_{U\cap E} \textrm{ for all }m\ge 1]\, .
\end{eqnarray*}
Since $w$ is dominated by $u$ which is continuous on $D$, it is bounded on the relatively compact set $U$.
Hence by the dominated convergence theorem and \eqref{e:w-to-0},
\begin{eqnarray*}
\lefteqn{\lim_{n\to \infty}\E_x[w(X^E_{\tau_{U_n}}):\, \tau_{U_m}<\tau_{U\cap E} \textrm{ for all }m\ge 1]}\\
&=&\E_x\left[\lim_{n\to \infty}w(X^E_{\tau_{U_n}})\I_{(X^E_{\tau_{U\cap E}-}\in (\partial E\cap D)\setminus N)}:\, \tau_{U_m}<\tau_{U\cap E} \textrm{ for all }m\ge 1\right]=0\, .
\end{eqnarray*}
Further,
\begin{eqnarray*}
\lefteqn{\lim_{n\to \infty}\E_x[w(X^E_{\tau_{U_n}}):\, \tau_{U_m}=\tau_{U\cap E} \textrm{ for some }m\ge 1\}}\\
&=&\E_x\left[\lim_{n\to \infty}w(X^E_{\tau_{U_n}}):\, \tau_{U_m}=\tau_{U\cap E} \textrm{ for some }m\ge 1\} \right]\\
&=&\E_x\left[w(X^E_{\tau_{U\cap E}}):\, \tau_{U_m}=\tau_{U\cap E} \textrm{ for some }m\ge 1\} \right]\\
&=&\E_x[w(X^E_{\tau_{U\cap E}}):\, \tau_U<\tau_E]\, .
\end{eqnarray*}
This proves the lemma. \qed

\begin{lemma}\label{l:v-supermedian}
Let
\begin{equation}\label{e:v-definition}
v(x):=\left\{\begin{array}{cl}
w(x)+R_u^{D\setminus E}(x)\, , & x\in E\\
u(x)\, , &x\in D\setminus E\, .
\end{array}\right.
\end{equation}
For every bounded open set $U\subset \overline{U}\subset D$ it holds that
\begin{equation}\label{e:v-supermedian}
\E_x[v(X^D_{\tau_U})]\le v(x)\, , \qquad x\in U\, .
\end{equation}
\end{lemma}
\pf We first note that $v\le (u-R_u^E)+R_u^E=u$ in $E$, and clearly $v=u$ in $D\setminus E$. Hence, if $x\in U\cap (D\setminus E)$, then $\E_x[v(X^D_{\tau_U})]\le \E_x[u(X^D_{\tau_U})]=u(x)=v(x)$.

Assume now that $x\in U\cap E$. Since $R_u^{D\setminus E}=u$ on $D\setminus E$, we have
\begin{eqnarray*}
\E_x[v(X^D_{\tau_U})]&=&\E_x[v(X^D_{\tau_U});\, X^D_{\tau_U}\in E]+\E_x[v(X^D_{\tau_U});\, X^D_{\tau_U}\in D\setminus E]\\
&=& \E_x[w(X^D_{\tau_U});\, X^D_{\tau_U}\in E]+\E_x[R_u^{D\setminus E}(X^D_{\tau_U});\, X^D_{\tau_U}\in E]\\
& &+ \E_x[R_u^{D\setminus E}(X^D_{\tau_U});\, X^D_{\tau_U}\in D\setminus E]\\
&=&\E_x[w(X^D_{\tau_U});\, X^D_{\tau_U}\in E] + \E_x[R_u^{D\setminus E}(X^D_{\tau_U})]=:A+B\, .
\end{eqnarray*}
Next, by using that $w=0$ on $D\setminus E$, and the fact that $X^D_t=X^E_t$ for all $t<\tau_E$,
\begin{eqnarray*}
A&=&\E_x[w(X^D_{\tau_U});\, X^D_{\tau_U}\in E, \tau_U<\tau_E]+\E_x[w(X^D_{\tau_U});\, X^D_{\tau_U}\in E, \tau_E<\tau_U]\\
&=&\E_x[w(X^E_{\tau_U})]+\E_x[w(X^D_{\tau_U});\,  \tau_E<\tau_U]\\
&=&\E_x[w(X^E_{\tau_{U\cap E}})]+\E_x[w(X^D_{\tau_U});\,  \tau_E<\tau_U]\\
&=&w(x)+\E_x[w(X^D_{\tau_U});\,  \tau_E<\tau_U]=:w(x)+A_2\, .
\end{eqnarray*}
In the last line we used Lemma \ref{l:w-harmonic}. We split $B$ into two parts:
$$
B=\E_x[R_u^{D\setminus E}(X^D_{\tau_U});\,  \tau_E<\tau_U]+\E_x[R_u^{D\setminus E}(X^D_{\tau_U});\, \tau_U\le \tau_E]:=B_1+B_2\, ,
$$
and combine $B_1$ with $A_2$:
\begin{eqnarray*}
A_2+B_1&=&\E_x[w(X^D_{\tau_U});\,  \tau_E<\tau_U]+
\E_x[R_u^{D\setminus E}(X^D_{\tau_U});\,  \tau_E<\tau_U]\\
&=&\E_x[(w+R_u^{D\setminus E})(X^D_{\tau_U});\,  \tau_E<\tau_U]\\
&\le &\E_x[u(X^D_{\tau_U});\,  \tau_E<\tau_U]\\
&=&\E_x\left[\E_{X^D_{\tau_E}}\left(u(X^D_{\tau_U})\right);\, \tau_E<\tau_U\right] \\
&=&\E_x[u(X^D_{\tau_E});\, \tau_E<\tau_U ].
\end{eqnarray*}
In the penultimate line we used the strong Markov property at time $\tau_E$, and in the last line harmonicity of $u$ (note that $X^D_{\tau_E}\in U\setminus E$ on $\tau_E<\tau_U$).

Finally, for $B_2$ we use that $N:=\{\wh{R}_u^{D\setminus E}\neq R_u^{D\setminus E}\}$ is polar, hence $\P_x(X^D_{\tau_U}\in N)=0$. Therefore, by using that $\wh{R}_u^{D\setminus E}(y)=\E_y[u(X^D_{\tau_E})]$ in the second line, and the strong Markov property in the third,
\begin{eqnarray*}
B_2&=&\E_x[\wh{R}_u^{D\setminus E}(X^D_{\tau_U});\, \tau_U\le \tau_E]\\
&=&\E_x\left[\E_{X^D_{\tau_U}}\left(u(X^D_{\tau_E})\right);\, \tau_U\le \tau_E\right]\\
&=&\E_x[u(X^D_{\tau_E});\, \tau_U \le \tau_E]\, .
\end{eqnarray*}
Putting everything together we get
\begin{eqnarray*}
\E_x[v(X^D_{\tau_U})]&=& w(x)+A_1+B_1+B_2\\
& \le& w(x)+ \E_x[u(X^D_{\tau_E});\, \tau_E<\tau_U]+\E_x[u(X^D_{\tau_E});\, \tau_U \le \tau_E]\\
&=&w(x)+\E_x[u(X^D_{\tau_E})]=w(x)+R_u^{D\setminus E}(x)\\
&=&w(x)+\wh{R}_u^{D\setminus E}(x)=v(x)\, .
\end{eqnarray*}
In the last line we used that $\wh{R}_u^{D\setminus E}= R_u^{D\setminus E}$ on $E$. \qed

\begin{lemma}\label{l:v-regularization}
Let $v$ be defined by \eqref{e:v-definition} and let $\wh{v}(x):=\liminf_{y\to x}v(y)$ be its lower semi-continuous regularization. Then $\wh{v}$ is excessive for $X^D$. Moreover, $\wh{v}\le u$ and there exists a polar set $N\subset \partial E\cap D$ such that $\wh{v}=u$ on $(D\setminus E)\setminus N$.
\end{lemma}
\pf First note that $\wh{v}\le v$ on $D$. Let $U\subset \overline{U}\subset D$ be open. Define
$$
\wt{v}(x):=\E_x[\wh{v}(X^D_{\tau_U})]\, ,\quad x\in D\, .
$$
By the proof of Lemma III.1.2 in \cite{BH}, 
$\wt{v}$ is lower semi-continuous in $U$. Moreover, by Lemma \ref{l:v-supermedian},
$$
\wt{v}(x)=\E_x[\wh{v}(X^D_{\tau_U})]\le \E_x[v(X^D_{\tau_U})]\le v(x)\, , \qquad x\in U\, .
$$
Hence, by lower semi-continuity of $\wt{v}$ in $U$, for every $x\in U$,
$$
\wt{v}(x)\le \liminf_{y\to x, y\in U}\wt{v}(y)\le \liminf_{y\to x, y\in U}v(y) =\wh{v}(x)\, .
$$
This proves that
\begin{equation}\label{e:v-superharmonic}
\E_x[\wh{v}(X^D_{\tau_U})]\le \wh{v}(x)\, \qquad \textrm{ for all } x\in U\, .
\end{equation}

Now, for any open $U\subset \overline{U}\subset D$, let $H_U(x,dy)=\P_x(X^D_{\tau_U}\in dy)$. Then the family $H_U(x, \cdot)$ (over all relatively compact open $U\subset D$) forms a family of harmonic kernels, 
cf. Chapter II of \cite{BH}. In the notation of \cite{BH}, 
\eqref{e:v-superharmonic} means that $\wh{v}\in {^*{\mathcal H}^+}(D)$. 
By Corollary III.2.1 of \cite{BH}, the latter family is equal to 
$\SS(X^D)$. Hence, $\wh{v}$ is excessive with respect to $X^D$.

Clearly, $\wh{v}\le v\le u$ on $D$.
Recall that $v=u$ on $D\setminus E$.
Let $z\in \partial E\cap D$.  Then
\begin{eqnarray*}
\liminf_{x\to z, x\in E}v(x)&=&\liminf_{x\to z, x\in E}(w(x)+R_u^{D\setminus E}(x))\\
&\ge &\liminf_{x\to z, x\in E}w(x)+\liminf_{x\to z, x\in E}R_u^{D\setminus E}(x)\\
&\ge & \liminf_{x\to z, x\in E}w(x)+ \wh{R}_u^{D\setminus E}(z)
\end{eqnarray*}
(since $\wh{R}_u^{D\setminus E}$ is the l.s.c.~regularization of $R_u^{D\setminus E}$).
By \eqref{e:w-to-0}, $\liminf_{x\to z, x\in E}w(x)=0$ for all $z\in (\partial E\cap D)\setminus N_1$
with $N_1$ being a polar set.
Also, $\wh{R}_u^{D\setminus E}=R_u^{D\setminus E}$ except on a polar set $N_2$. By setting $N=N_1\cup N_2$, we see that for all $z\in (\partial E\cap D)\setminus N$,
$$
\liminf_{x\to z, x\in E}v(x)\ge R_u^{D\setminus E}(z)=u(z)\, .
$$
Clearly, for all $z\in \partial E\cap D$,
$$
\liminf_{x\to z, x\in D\setminus E}v(x)=\liminf_{x\to z, x\in D\setminus E}u(x)\ge u(z)\, .
$$
Together the last two displays give that for all $z\in (\partial E\cap D)\setminus N$,
$$
\wh{v}(z)=\liminf_{x\to z}v(x)\ge u(z)\, .
$$
\qed

We note that for every $g:D\setminus E\to [0,\infty)$ the function $x\mapsto \E_x[g(X_{\tau_E})]=\wh{R}_g^{D\setminus E}(x)$ is harmonic in $E$ with respect to $X^D$. Since $R_g^{D\setminus E}=\wh{R}_g^{D\setminus E}$ on $E$, it follows that $R_g^{D\setminus E}$ is harmonic in $E$ with respect to $X^D$.

\medskip
In what follows, $\partial_M D$ denotes the Martin boundary of $D$ with respect to $X^D$, $\partial_m D$ the minimal Martin boundary, and $D_M=D\cup \partial_M D$ the Martin space (with the Martin topology). For $z\in \partial_m D$ let $M_D(\cdot, z)$ be the Martin kernel (based at $x_0\in E$). 
Then $M_D(\cdot, z)$ is continuous and harmonic in $D$ with respect to $X^D$.
We recall that $E\subset D$ is \emph{minimally thin} in $D$ at $z\in \partial_m D$ with respect to $X^D$ if $\wh{R}_{M_D(\cdot, z)}^E\neq M_D(\cdot, z)$.

\begin{prop}\label{p:h-minimal}
Let $E\subset D$ be an open set in $D$, $z\in \partial_m D$ such that $z$ is in the closure of $E$ in $D_M$. Assume that $D\setminus E$ is minimally thin at $z$ in $D$ with respect to $X^D$. Let
$$
h(x):=M_D(x,z)-R_{M_D(\cdot,z)}^{D\setminus E}(x)\, ,\qquad x\in E\, .
$$
Then $h$ is a minimal harmonic function in $E$ with respect to $X^E$.
\end{prop}
\pf
We first prove that $h$ is harmonic with respect to $X^E$. Let $U\subset \overline{U}\subset E$ be relatively compact open in $E$. Then
\begin{eqnarray*}
\lefteqn{\E_x\left[h(X^E_{\tau_U}) \right]
=\E_{x}\left[M_D(X^E_{\tau_U} ,z) \right]-
\E_x\left[ \E_{X^E_{\tau_U}} \left[M_D\left(X^D_{S_{D\setminus E}}, z \right)   \right] \right]}\\
&= &\E_x\left[M_D(X^D_{\tau_U} ,z) \right]-\E_x\left[M_D(X^D_{\tau_U} ,z);\, \tau_U=\tau_E \right]\\
&&\qquad-
\E_x\left[ \E_{X^D_{\tau_U}} \left[M_D\left(X^D_{S_{D\setminus E}}, z \right)   \right];\, \tau_U < \tau_E \right]\\
&= &\E_x\left[M_D(X^D_{\tau_U} ,z) \right]-\E_x\left[M_D(X^D_{\tau_E} ,z);\,  \tau_U=\tau_E \right]\\
&&\qquad -\E_x \left[M_D\left(X^D_{S_{D\setminus E}}, z \right);\, \tau_U < \tau_E \right]\\
&= &\E_x\left[M_D(X^D_{\tau_U} ,z) \right]-\E_x\left[M_D(X^D_{D\setminus E} ,z);\,  \tau_U=S_{D\setminus E}\right]\\
&&\qquad -
\E_{x} \left[M_D\left(X^D_{D_{D\setminus E}}, z \right);\, \tau_U < S_{D\setminus E} \right]\\
&=&M_D(x,z)-\E_x \left[M_D\left(X^D_{S_{D\setminus E}}, z \right)   \right]=h(x).
\end{eqnarray*}

Now suppose that $w:E\to[0,\infty)$ is harmonic in $E$ with respect to $X^E$ and $w\le h$. Define $v$ analogously to \eqref{e:v-definition} by
$$
v(x):=\left\{\begin{array}{cl}
w(x)+R_{M_D(,\cdot,z)}^{D\setminus E}(x)\, , & x\in E\\
M_D(x,z)\, , &x\in D\setminus E\, ,
\end{array}\right.
$$
and let $\wh{v}$ be its l.s.c.~regularization. 
By Lemma \ref{l:v-regularization}, 
$\wh{v}\in \SS(X^D)$, $\wh{v}\le M_D(\cdot,z)$ on $D$, and $\wh{v}=M_D(\cdot, z)$ on $(D\setminus E)\setminus N$, $N$ polar. By the Riesz decomposition,
$$
\wh{v}=aM_D(\cdot,z)+G_D \mu\, ,
$$
where $0\le a\le 1$ and $\mu$ is a measure charging no polar set (since $\wh{v}$ is locally bounded, the same holds for $G_D\mu$, hence $\mu$ cannot charge polar sets). Note that $\wh{v}=w+R_{M_D(\cdot, z)}^{D\setminus E}$ on $E$.
The function $R_{M_D(\cdot, z)}^{D\setminus E}$ is harmonic in $E$ with respect to $X^D$. By assumption, $w$ is harmonic in $E$ with respect to $X^E$, and hence harmonic in $E$ with respect to $X^D$ (we extend $w=0$ on $D\setminus E$). Therefore, $\wh{v}$ is harmonic in $E$ with respect to $X^D$ which implies the same for $G_D\mu=\wh{v}-aM_D(\cdot, z)$.

Recall that $D\setminus E$ is thin at $y\in D$ if and only if $\wh{R}_{G_D(\cdot,y)}^{D\setminus E}\neq G_D(\cdot, y)$
(this can be proved along the same lines as the corresponding proof for minimal thinness, 
cf. Proposition 6.2 of \cite{KSV14b}).
Let
$$
A=\{y\in \partial E\cap D:\, \wh{R}_{G_D(\cdot,y)}^{D\setminus E}\neq G_D(\cdot, y)\}\, .
$$
By Proposition VI.5.12 of \cite{BG}, $A$ is polar, and hence $\mu(A)=0$.

Now consider $R_{G_D\mu}^{D\setminus E}$. This function is harmonic in $E$ with respect to $X^D$. Moreover, $R_{G_D\mu}^{D\setminus E}\le G_D\mu$ on $D$. Hence, $G_D\mu -R_{G_D\mu}^{D\setminus E}\ge 0$ and is harmonic in $E$ with respect to $X^D$. Note that $G_D\mu -R_{G_D\mu}^{D\setminus E}=0$ on $D\setminus E$. Hence, $G_D\mu -R_{G_D\mu}^{D\setminus E}$ is harmonic in $E$ with respect to $X^E$. On the other hand, for $x\in E$,
\begin{eqnarray*}
\lefteqn{G_D\mu(x)-R_{G_D\mu}^{D\setminus E}(x)=G_D\mu(x)-\wh{R}_{G_D\mu}^{D\setminus E}(x)}\\
&=&\int_D G_D(x,y)\, \mu(dy)-\int_D\wh{R}_{G_D(\cdot,y)}^{D\setminus E}(x)\, \mu(dy)\\
&=&\int_D\left[G_D(x,y)-\wh{R}_{G_D(\cdot,y)}^{D\setminus E}(x)\right]\, \mu(dy)\\
&=&\int_E\left[G_D(x,y)-\wh{R}_{G_D(\cdot,y)}^{D\setminus E}(x)\right]\, \mu(dy)+\int_A\left[G_D(x,y)-\wh{R}_{G_D(\cdot,y)}^{D\setminus E}(x)\right]\, \mu(dy)\\
&=&\int_E G_E(x,y)\, \mu(dy)=G_E\mu(x)\, .
\end{eqnarray*}
In the last line we used that $\mu(A)=0$ and the formula for the Green function of $X^E$: $G_E(x,y)=G_D(x,y)-\E_x[G_D(X_{\tau_E},y)]$. This shows that $G_D\mu -R_{G_D\mu}^{D\setminus E}$ is at the same time harmonic in $E$ (with respect to $X^E$) and the potential of the measure $\mu_{|E}$. Hence, it is identically zero in $E$, that is, $G_D\mu =R_{G_D\mu}^{D\setminus E}$ in $E$, hence in $D$.

Since $v$ and $\wh{v}$ differ at most on a polar set, and $v=M_D(\cdot,z)$ on $D\setminus E$, we see that $G_D\mu=(1-a)M_D(\cdot, z)$ outside a polar set. Therefore
$$
(1-a)R_{M_D(\cdot,z)}^{D\setminus E}=R_{G_D\mu}^{D\setminus E} =G_D\mu\, .
$$
Hence, on $E$ we have
\begin{align*}
w&=\wh{v}-R_{M_D(\cdot,z)}^{D\setminus E}=aM_D(\cdot,z)+(1-a)R_{M_D(\cdot,z)}^{D\setminus E}-R_{M_D(\cdot,z)}^{D\setminus E}\\
&=a\left(M_D(\cdot,z)-R_{M_D(\cdot,z)}^{D\setminus E}\right)=ah\,
,
\end{align*}
which completes the proof. \qed

\begin{remark}
{\rm
The assumption that $D\setminus E$ is minimally thin at $z$ in $D$ with respect to $X^D$ is used to conclude that $h\neq 0$.
If $D\setminus E$ is not minimally thin at $z$ in $D$ with respect to $X^D$, then $R_{M_D(\cdot,z)}^{D\setminus E}=M_D(\cdot, z)$.
}
\end{remark}

\begin{prop}\label{p:convergence}
Let $E\subset D$ be an open set in $D$, $z\in \partial_m D$ such that $z$ is in the closure of $E$ in $D_M$. Assume that $D\setminus E$ is minimally thin at $z$ in $D$ with respect to $X^D$ and let
$$
h(x):=M_D(x,z)-R_{M_D(\cdot,z)}^{D\setminus E}(x)\, ,\qquad x\in E\, .
$$
Let $\zeta=\zeta(z)$ be the Martin boundary point of $E$ associated with the minimal harmonic function $h$. Assume that $(x_n)_{n\ge 1}$ is a sequence of points in $E$ that converges to $z$ in $D_M$ and also
\begin{equation}\label{e:liminf-positive}
\liminf_{n\to\infty}\frac{G_E(x_0,x_n)}{G_D(x_0,x_n)}>0\, .
\end{equation}
Assume further that for every subsequence $(x_{n_k})$, $G_E(\cdot, x_{n_k})/G_E(x_0, x_{n_k})$ converges to a harmonic function with respect to $X^E$.
Then $(x_n)_{n\ge 1}$ converges to $\zeta$ in $E_M$ (Martin space of $E$).
\end{prop}
\pf Let $(x_n)_{n\ge 1}$ be a sequence in $E$ converging to $z$ in $D_M$ and such that \eqref{e:liminf-positive} holds. Assume that $(x_n)$ does not converge to $\zeta$ in $E_M$. This implies that there exists a subsequence $(x_{n_k})$ with the property that $G_E(\cdot, x_{n_k})/G_E(x_0, x_{n_k})$ converges in $E$ to a function $u:E\to [0,\infty)$ such that $u\neq h/h(x_0)$. By assumption, $u$ is harmonic with respect to $X^E$. It follows from \eqref{e:liminf-positive}, that by choosing a further subsequence (if necessary) we can arrange that
$$
\lim_{k\to \infty}\frac{G_E(x_0, x_{n_k})}{G_D(x_0, x_{n_k})}=a>0\, .
$$
Therefore, on $E$ we have that
$$
\lim_{k\to \infty}\frac{G_E(\cdot, x_{n_k})}{G_D(x_0,x_{n_k})}=au\, .
$$
Since $G_E(\cdot, y)=G_D(\cdot,y)-R_{G_D(\cdot,y)}^{D\setminus E}$, and since
$$
\frac{R_{G_D(\cdot,y)}^{D\setminus E}(x_0)}{G_D(x_0,y)}=R_{M_D(\cdot,y)}^{D\setminus E}(x_0)
$$
(which easily follows from the probabilistic representation of the reduced function), we get by use of Fatou's lemma in the last line that
\begin{eqnarray*}
a u(x)&=& \lim_{k\to \infty}\left(\frac{G_D(x,x_{n_k})}{G_D(x_0,x_{n_k})}-\frac{R_{G_D(\cdot,x_{n_k})}^{D\setminus E}(x)}{G_D(x_0,x_{n_k})}\right)\\
&=&\lim_{k\to \infty}\left(M_D(x,x_{n_k})-R_{M_D(\cdot,x_{n_k})}^{D\setminus E}(x)\right)\\
&\le & M_D(x,z)-R_{M_D(\cdot,z)}^{D\setminus E}(x)=h(x)\, .
\end{eqnarray*}
Since $u$ is harmonic for $X^E$, it follows from Proposition \ref{p:h-minimal} that $au$ is proportional to $h$. Since $u(x_0)=1$, that would imply $u=h/h(x_0)$ which contradicts the assumption. \qed

If $F\subset E\subset D$, and $v\in \SS(X^E)$, let $^E R_v^F$ denote the reduced function of $v$ on $F$ with respect to $X^E$.

\begin{lemma}\label{l:reduits}
Let $F\subset E\subset D$, $u\in \SS(X^D)$, and define $v:=u-R_u^{D\setminus E}$. Then $v\in \SS(X^E)$ and
\begin{equation}\label{e:reduits}
^E R_v^F=R_u^{D\setminus F}-R_u^{D\setminus E}\, .
\end{equation}
\end{lemma}
\pf Since the excessiveness implies that $u(x) \ge \E_x [u(X^D_{S_{D\setminus E}})]$, $v$ is non-negative. If $x \in E$, by the strong Markov property,
\begin{eqnarray*}
\E_x\left[v(X^U_t)\right]&=&\E_x\left[u(X^E_t)\right] -\E_x\left[ \E_{X^E_t} \left[u\left(X^D_{S_{D\setminus E}} \right)   \right] \right]\\
&=&\E_x\left[u(X^D_t)\right]-\E_x\left[u(X^D_t): t \ge \tau_E\right]\\
&&\qquad -\E_x\left[ \E_{X^D_t} \left[u\left(X^D_{S_{D\setminus E}} \right)   \right]:  t < \tau_E \right]\\
&=&\E_x\left[u(X^D_t)\right]-\E_x\left[u(X^D_t): t \ge \tau_E\right] -\E_x\left[ u\left(X^D_{S_{D\setminus E}} \right)   :  t < \tau_E \right]\\
&=&\E_x\left[u(X^D_t)\right]-\E_x\left[ u\left(X^D_{t \vee \tau_E} \right)\right].
\end{eqnarray*}
 By the excessiveness of $u$ for $X^D$, $\E_x\left[u(X^D_t)\right] \le u(x)$ and
$$ \E_x\left[ u\left(X^D_{t \vee \tau_E} \right) \right] \ge \E_x\left[ u\left(X^D_{\tau_E} \right)   \right].$$
Thus
$$
\E_x\left[v(X^E_t)\right]
\le u(x) -\E_x \left[u\left(X^D_{S_{D\setminus E}} \right)   \right].
$$
Moreover
\begin{align*}
\lim_{t \downarrow 0} \E_x\left[v(X^E_t)\right]
&= \lim_{t \downarrow 0}\E_x\left[u(X^D_t)\right] - \lim_{t \downarrow 0}\E_x\left[ u\left(X^D_{t \vee \tau_E} \right)\right]\\
& =u(x) -\E_x \left[u\left(X^D_{S_{D\setminus E}} \right)   \right].
\end{align*}

Note that for $x \in E$,
\begin{align*}
&\E_x \left[u\left(X^D_{S_{D\setminus F} } \right ) \right]\\
=&
\E_x \left[u\left(X^D_{S_{E\setminus F}} \right ); X_{S_{D\setminus F}} \in E  \right]+
\E_x \left[u\left(X^D_{S_{D\setminus F} } \right ); S_{D\setminus F}=S_{D\setminus E}   \right]\\
=&\E_x \left[u\left(X^E_{S_{E\setminus F}} \right )  \right]+
\E_x \left[u\left(X^D_{S_{D\setminus E} } \right ); S_{D\setminus F}=S_{D\setminus E}   \right]\\
=&\E_x \left[u\left(X^E_{S_{E\setminus F}} \right )  \right]+
\E_x \left[u\left(X^D_{S_{D\setminus E} } \right )  \right]-
\E_x \left[u\left(X^D_{S_{D\setminus E} } \right ); S_{D\setminus F}<S_{D\setminus E}   \right].
\end{align*}
By strong Markov property,
$$
\E_x \left[u\left(X^D_{S_{D\setminus E} } \right ); S_{D\setminus F}<S_{D\setminus E}   \right]=
\E_x \left[\E_{X^E_{S_{E\setminus F}}} \left[u\left(X^D_{S_{D\setminus E}} \right)   \right]\right ] .
$$
Thus
\begin{align*}
&\E_x \left[u\left(X^D_{S_{D\setminus F} } \right ) \right]\\
=&\E_x \left[u\left(X^D_{S_{D\setminus E}  }  \right)\right]+\E_x \left[u\left(X^E_{S_{E\setminus F}} \right )  \right]
-\E_x \left[\E_{X^U_{S_{E\setminus F}}} \left[u\left(X^D_{S_{D\setminus E}} \right)   \right]\right ].
\end{align*}
Therefore
$$
\E_x \left[u\left(X^D_{S_{D\setminus F} } \right ) \right]
=\E_x \left[u\left(X^D_{S_{D\setminus E}  }  \right)\right]+\E_x \left[v\left(X^E_{S_{E\setminus F}} \right )  \right],
$$
which is \eqref{e:reduits}. \qed

\begin{prop}\label{p:min-thin}
Let $E\subset D$ be an open set in $D$, $z\in \partial_m D$ such that $z$ is in the closure of $E$ in $D_M$. Assume that $D\setminus E$ is minimally thin at $z$ in $D$ with respect to $X^D$ and let
$$
h(x):=M_D(x,z)-R_{M_D(\cdot,z)}^{D\setminus E}(x)\, ,\qquad x\in E\, .
$$
Let $\zeta=\zeta(z)$ be the Martin boundary point of $E$ associated with the minimal harmonic function $h$. Let $F\subset E$. Then $F$ is minimally thin at $\zeta$ in $E$ with respect to $X^E$ if and only if $F$ is minimally thin at $z$ in $D$ with respect to $X^D$.
\end{prop}
\pf The set $F$ is minimally thin at $\zeta$ with respect to $X^E$ if and only if $^E R_h^F\neq h$. By Lemma \ref{l:reduits} (with $F$ replaced by $E\setminus F$ and $u=M_D(\cdot, z)$),
$$
^E R_h^{E\setminus F}=R_{M_D(\cdot, z)}^{D\setminus (E\setminus F)}-R_{M_D(\cdot, z)}^{D\setminus E}=R_{M_D(\cdot, z)}^{(D\setminus E)\cup F}-R_{M_D(\cdot, z)}^{D\setminus E}\, .
$$
Since $h=M_D(\cdot, z)-R_{M_D(\cdot,z)}^{D\setminus E}$, we see that $^E R_h^F\neq h$ if and only if $R_{M_D(\cdot, z)}^{(D\setminus E)\cup F}\neq M_D(\cdot, z)$. The last condition is equivalent to $(D\setminus E)\cup F$ being minimally thin at $z$ in $D$ with respect to $X^D$. Since $D\setminus E$ is not minimally thin at $z$, the latter is equivalent to $F$ being minimally thin at $z$ in $D$ with respect to $X^D$. \qed

\begin{remark}{\rm
Proposition \ref{p:min-thin} does not depend on Proposition \ref{p:convergence}.
}
\end{remark}

\bigskip

Let $D\subset \X$ be an open unbounded set.
Suppose $E$ is an open subset of $D$ such that for some $R>0$ it holds that
$D\cap \overline{B}(z_0,R)^c=E\cap \overline{B}(z_0,R)^c$.
Assume that $\infty$ is accessible both from $E$ and from $D$.
Assume that the assumptions \textbf{A},  \textbf{C}, \textbf{C2}$(z_0, R)$ and \textbf{F2}$(z_0, R)$ for $X$ and  $\wh{X}$ are satisfied.
By Theorem \ref{t:main-mb1}
there is only one Martin boundary point of $E$ associated with $\infty$, say $\infty^E$, and this point is minimal, $\infty^E\in \partial_m E$. In the same way, there is only one Martin boundary point of $D$ associated with $\infty$, say $\infty^D$, and this point is also minimal, $\infty^D\in \partial_m D$.
Hence, the concept of minimal thinness at $\infty$ of a set $F\subset E$ makes sense with respect to both $X^E$ and $X^D$. In fact, we have the following result.

\begin{thm}\label{t:min-thin-infty}
Suppose that  \textbf{A},  \textbf{C}, \textbf{C2}$(z_0, R)$ and
\textbf{F2}$(z_0, R)$ for $X$ and $\wh{X}$  hold true.
Let $D\subset \X$ be an unbounded open set,
and let $E$ be an open subset of $D$ such that for some $R>0$ it holds that $D\cap \overline{B}(z_0,R)^c=E\cap \overline{B}(z_0,R)^c$. Assume that $\infty$ is accessible from $E$ and from $D$. Suppose that $F\subset E$. Then $F$ is minimally thin at $\infty$ with respect to $X^E$ if and only if it is minimally thin at $\infty$ with respect to $X^D$.
\end{thm}
\pf
Let $x_0\in E$ and choose $r_0 > 2(d(x_0, E^c)\wedge R)$. For every $r\in (0,r_0)$, both  $G_D(x_0, \cdot)$ and $G_E(x_0, \cdot)$ are regular harmonic  in $D\cap \overline{B}(z_0,r)^c$ with respect to $X$ and vanish in $\overline{B}(z_0,r)^c\setminus D$. Let $x_1\in E\cap \overline{B}(z_0,2r)^c=D\cap \overline{B}(z_0,2r)^c$ be fixed. By the boundary Harnack principle,
$$
\frac{G_E(x_0,x)}{G_D(x_0,x)}\ge c^{-1} \frac{G_E(x_0,x_1)}{G_D(x_0,x_1)}\, , \qquad \textrm{for all }x\in D\cap \overline{B}(z_0,8r)^c\, .
$$
This implies that
\begin{equation}\label{e:min-thin-infty}
\liminf_{E\ni x\to \infty}\frac{G_E(x_0,x)}{G_D(x_0,x)} >0\, .
\end{equation}
Let $M_D(\cdot, \infty)=M_D(\cdot, \infty^D)$, respectively $M_E(\cdot, \infty)=M_E(\cdot, \infty^E)$, be the Martin kernels at $\infty$ for $D$, respectively $E$. Define
$$
h(x)=M_D(x,\infty)-R_{M_D(\cdot, \infty)}^{D\setminus E} (x).
$$
By Proposition \ref{p:h-minimal}, $h$ is a minimal harmonic function with respect to $X^E$. Let $\zeta\in \partial_m E$ be the minimal Martin boundary point of $E$ corresponding to $h$.
Let $(x_n)_{n\ge 1}$ be a sequence of points in $E$ converging to $\infty$. By \eqref{e:min-thin-infty},
$$
\liminf_{n\to \infty}\frac{G_E(x_0,x_n)}{G_D(x_0,x_n)} >0\, .
$$
Also, it follows from Lemma \ref{l:limit-exists} and Theorem \ref{t:main-mb1}(b)
that, for every subsequence $(x_{n_k})$, $G_E(\cdot, x_{n_k})/G_E(x_0,x_{n_k})$ converges to the harmonic function $M_E(\cdot, \infty)$. It follows from Proposition \ref{p:convergence} that $(x_n)_{n\ge 1}$ converges to $\zeta$ in the Martin topology of $E_M$. Thus, $\zeta\in \partial_m E$ is
associated to $\infty$.
By uniqueness, $\zeta=\infty^E$ and therefore $h=M_E(\cdot, \infty^E)=M_E(\cdot,\infty)$.
The claim of the theorem now follows from Proposition \ref{p:min-thin}. \qed

\begin{remark}
{\rm
Suppose that $\infty$ is accessible from $E$. Since $G_E(x,w)\le G_D(x,w)$, $x,w\in E$, implies that $\E_x \tau_E\le \E_x \tau_D$ for $x\in E$, we see that $\E_x \tau_D=\infty$ for all $x\in E$. If the assumptions of \cite{KSVp1} are satisfied, it follows from Remark \ref{r:criterion-accessibility} that $\infty$ is also accessible from $D$.
}
\end{remark}

One can similarly prove the following theorem saying that minimal thinness is a local property at a finite boundary point.

\begin{thm}\label{t:min-thin-local}
Suppose that  \textbf{A},  \textbf{C}, \textbf{C1}$(z_0, R)$ and
\textbf{F1}$(z_0, R)$ for $X$ and $\wh{X}$ hold true.
Let $D\subset \X$,
$z_0\in \partial D$, and let $E$ be an open subset of $D$ such that for some $R>0$ it holds that $D\cap B(z_0,R)=E\cap B(z_0,R)$. Assume that $z_0$ is accessible from $E$ and from $D$. Suppose that $F\subset E$. Then $F$ is minimally thin at $z_0$ with respect to $X^E$ if and only if it is minimally thin at $z_0$ with respect to $X^D$.
\end{thm}

\bigskip
\noindent
{\bf Acknowledgements:} Part of the research for this paper was done
during the visit of Renming Song and Zoran Vondra\v{c}ek to Seoul National University
from May 24 to June 8 of 2015.
They thank the Department of
Mathematical Sciences of Seoul National University for the hospitality.
We also thank the referee for very helpful comments on the first version of this
paper.

\end{doublespace}
\bigskip

\vspace{.1in}
\begin{singlespace}


\end{singlespace}

\
\vskip 0.1truein

\parindent=0em

{\bf Panki Kim}

Department of Mathematical Sciences and Research Institute of Mathematics,

Seoul National University, Building 27, 1 Gwanak-ro, Gwanak-gu Seoul 08826, Republic of Korea

E-mail: \texttt{pkim@snu.ac.kr}

\bigskip

{\bf Renming Song}

Department of Mathematics, University of Illinois, Urbana, IL 61801,
USA

E-mail: \texttt{rsong@math.uiuc.edu}

\bigskip

{\bf Zoran Vondra\v{c}ek}

Department of Mathematics, University of Zagreb, Zagreb, Croatia, and \\
Department of Mathematics, University of Illinois, Urbana, IL 61801,
USA

Email: \texttt{vondra@math.hr}

\end{document}